\documentclass[a4page]{article}
\usepackage[english]{babel}
\usepackage[pdftex]{color,graphicx,hyperref}
\usepackage{amsmath, amssymb}
\usepackage{amsthm}
\usepackage{fancyhdr}

\allowdisplaybreaks 

\usepackage[left=3.5cm,right=2.5cm,top=2.5cm,bottom=2.5cm,includeheadfoot]{geometry}

\usepackage{array}
\usepackage{booktabs}
\usepackage{tabularx}
\usepackage{longtable}

\usepackage{graphicx}
\graphicspath{{bilder/}}

\usepackage{amsmath}
\usepackage{amssymb}
\usepackage{amsfonts}
\usepackage{amsthm}

\usepackage{color}
\usepackage{xcolor}

\usepackage{rotating}

\usepackage{lscape}

\usepackage[multiple]{footmisc}

\usepackage[pdftex]{hyperref}

\usepackage{setspace}

\usepackage{caption}

\usepackage{pdfpages}

\usepackage{url}

\usepackage{setspace}

\usepackage{lipsum}

\usepackage{fancyhdr}
\def \qed{\hfill$\Box$}

\newcommand{\e}{\mathbb{E}} 
\newcommand{\E}{\mathbb{E}}
\newcommand{\var}{\mathbb{V}{\rm ar}} 
\newcommand{\p}{\mathbb{P}} 
\renewcommand{\L}{\mathbf{L}} 
\renewcommand{\P}{\mathbb{P}} 
\newcommand{\R}{\mathbb{R}} 
\newcommand{\N}{\mathbb{N}} 
\newcommand{\ZZ}{\mathbb{Z}} 

\newcommand{\1}{\mathbf{1}} 
    \def\independenT#1#2{\mathrel{\setbox0\hbox{$#1#2$}%
    \copy0\kern-\wd0\mkern4mu\box0}} 

\newcommand{\vecx}[1]{\mathbf{\underline{#1} } } 

\newcommand{\cB}{\mathcal{B}}
\newcommand{\cC}{\mathcal{C}}
\newcommand{\cE}{\mathcal{E}}
\newcommand{\cF}{\mathcal{F}}
\newcommand{\cG}{\mathcal{G}}
\newcommand{\cH}{\mathcal{H}}
\newcommand{\cK}{\mathcal{K}}

\newcommand{\cN}{\mathcal{N}}

\newcommand{\cQ}{\mathcal{Q}}

\newcommand{\cT}{\mathcal{T}}
\newcommand{\cX}{\mathcal{X}}

\newcommand{\Om}[1]{\Omega^{(#1)} }
\newcommand{\F}[2]{\cF^{(#1)}_{#2}}
\newcommand{\Prob}[1]{\P^{(#1)}}
\newcommand{\CE}[1]{\E^{(#1)}}

\newcommand{\X}[2]{X^{(#1)}_{#2}}

\newcommand{\Samp}[1]{\cX^{(#1)} }
\newcommand{\emeas}[1]{\nu^{(#1)}_{M}}
\newcommand{\emeasi}[2]{\nu^{(#1)}_{#2}}

\newcommand{\DW}[2]{\Delta W^{(#1)}_{#2}}
\newcommand{\yM}[2]{y^{(M)}_{#1}(#2)}
\newcommand{\zM}[2]{z^{(M)}_{#1}(#2)}
\newcommand{\byM}[2]{\bar y^{(M)}_{#1}(#2)}
\newcommand{\bzM}[2]{\bar z^{(M)}_{#1}(#2)}
\newcommand{\ObsM}[2]{S^{(M)}_{#1}(#2)}
\newcommand{\ObsMk}[3]{S^{(#1,M)}_{#2}(#3)}

\newcommand{\Obsk}[3]{S^{(#1)}_{#2}(#3)}

\newcommand{\psim}{\psi^{(M)}}
\newcommand{\psimk}[1]{\psi^{(#1,M)}}
\newcommand{\psik}[1]{\psi^{(#1)}}

\newcommand{\plx}[3]{p^{(#2)}_{#1}(#3) }

\newcommand{\pl}[2]{p_{#1}^{(#2)} }
\newcommand{\LinSpace}[1]{\cK_{#1}}
\newcommand{\LinSpacek}[2]{\cK^{(#1)}_{#2}}
\newcommand{\BasDim}[1]{K_{#1}}
\newcommand{\BasDimk}[2]{K^{(#1)}_{#2}}

\newcommand{\OLS}{{\bf OLS}}

\newcommand{\HX}{{$\bf (A_X)$}}
\newcommand{\HXp}{{$\bf (A'_X)$}}
\newcommand{\HXpp}{{$\bf (A''_X)$}}
\newcommand{\HXppp}{{$\bf (A'''_X)$}}

\newcommand{\HK}{{$\bf (A_{\cK})$}}

\newcommand{\Cz}[1]{C_{z,#1}}

\newcommand{\thetaL}{\theta_L}

\newcommand{\thetaP}{\theta_\Phi}

\newcommand{\lawnorm}[2]{\|#1 \|_{#2,\infty} }
\newcommand{\empinorm}[2]{\| #1\|_{#2,M} }

\newcommand{\fM}{f^{(M)} }

\newtheorem{definition}{Definition}[section]
\newtheorem{theorem}[definition]{Theorem}
\newtheorem{proposition}[definition]{Proposition}
\newtheorem{lemma}[definition]{Lemma}
\newtheorem{corollary}[definition]{Corollary}
\theoremstyle{definition} \newtheorem{remark}[definition]{Remark}
\theoremstyle{definition} \newtheorem*{remark*}{Remark}

\theoremstyle{theorem} \newtheorem{alg}{Algorithm}

\newcommand \trace{{\rm Tr}}

\newcommand\cpib{{R_\pi}}

\newcommand {\thetal}{{\theta_L}}

\makeatletter
\def\timenow{\@tempcnta\time
\@tempcntb\@tempcnta
\divide\@tempcntb60
\ifnum10>\@tempcntb0\fi\number\@tempcntb
:\multiply\@tempcntb60
\advance\@tempcnta-\@tempcntb
\ifnum10>\@tempcnta0\fi\number\@tempcnta}
\makeatother

\newcommand{\Ht}{$\bf (A_{\pi })$}
\newcommand{\Htp}{$\bf (A'_{\pi })$}

\newcommand{\Hg}{$\bf (A_{\Phi })$}

\newcommand{\Hf}{$\bf (A_{f })$}

\newcommand{\tk}[2]{t^{(#1)}_{#2}}
\newcommand{\pik}[1]{\pi^{(#1)}}
\newcommand{\Y}[2]{y^{ (#1)}_{ #2}}
\newcommand{\y}[3]{y^{ (#1)}_{ #2}(#3)}

\newcommand{\yn}[2]{y^{ (#1)}_{ #2}}

\newcommand{\Z}[2]{z^{(#1)}_{ #2}}
\newcommand{\z}[3]{z^{(#1)}_{ #2}(#3)}

\newcommand{\zn}[2]{z^{(#1)}_{ #2}}

\newcommand{\Del}[1]{\Delta^{( #1)} }
 
\newcommand{\DL}[2]{\Delta L^{ (#1)}_{ #2} }

\newcommand{\Tk}[2]{T^{(#1) }_{ #2}}

\newcommand{\normM}[3]{\|#3 \|_{#1,#2,M}}

\newcommand{\tily}[2]{\tilde y^{(#1)}_{#2} }
\newcommand{\tilz}[2]{\tilde z^{(#1)}_{#2} }

\newcommand{\Xmarkov}[2]{ V^{(#1)}_{#2} }

\title{Multilevel approximation of backward stochastic differential equations 
}
\date{\today}

\author{
\begin{minipage}{0.5\textwidth}
\begin{centering}
D. Becherer 
\footnote{Authors acknowledge support from German Schience Fondation DFG, Berlin Mathematical School and {\sc Matheon}}
\\
Institut f\"ur Mathematik \\
Humboldt Universit\"at zu Berlin \\
Unter den Linden 6 \\
D-10099 Berlin, Germany \\
\end{centering}
\end{minipage}
\begin{minipage}{0.5\textwidth}
\begin{centering}
P. Turkedjiev\footnote{Corresponding author: \href{mailto:turkedjiev@cmap.polytechnique.fr}{\tt turkedjiev@cmap.polytechnique.fr}.
Research supported by the Chair Financial Risks of the Risk Foundation, the FiME Laboratory, and the Chair Finance and Sustainable Developement, under the aegis of Louis Bachelier Finance and Sustainable Growth laboratory, a joint initiative with Ecole Polytechnique. 
}\\
Centre de Math\'ematiques Appliqu\'ees\\
Ecole Polytechnique and CNRS\\
Route de Saclay\\
F 91128 Palaiseau cedex, France \\
\end{centering}
\end{minipage}
}

\begin{document}
\maketitle

\begin{abstract}
We develop a multilevel approach  to compute approximate solutions to backward differential equations (BSDEs).
The  fully implementable algorithm of our  multilevel scheme constructs 
sequential martingale control variates 
along a sequence of refining time-grids
 to  reduce statistical approximation errors in an adaptive and generic way.
We provide an error analysis with
 explicit 
and non-asymptotic error estimates for the multilevel scheme
under general conditions on the forward process and the BSDE data.
It is shown that the multilevel approach can reduce the computational complexity
to achieve precision $\epsilon$, ensured by error estimates, 
essentially by one order (in $\epsilon^{-1}$)
in comparison to established methods,
which is substantial. 
 Computational examples support the  validity of the theoretical analysis, demonstrating 
efficiency improvements in practice.

\end{abstract}

\section{Introduction}

The concept of Multilevel Monte Carlo has been introduced by \cite{gile:08} as a {simulation} method for the  efficient computation of linear expectations $E[\Phi(X)]$ 
of functions $\Phi$ of diffusion processes $X$;
 see also  \cite{Heinrich01}.  Multilevel Monte Carlo (MLMC) {is an active research area, evolving in} many directions; For instance, 
 \cite{Dereich11} studies the case where $X$ is the solution to a L\'evy-diven stochastic differential equation and \cite{BelomestnySD11}
 develops a multilevel approach for the problem of optimal stopping, where an expectation $E[\Phi(X_{\tau})]$ is maximized over a family of stopping times $\tau$.

{Our} paper develops a {novel} multilevel approximation algorithm  {for solutions} {to} backward stochastic differential equations, which can be seen as  a non-linear generalization 
of the probabilistic Feynman-Kac representation for linear expectations of diffusions, with many applications in 
optimal control and mathematical finance, see e.g.\ \cite{elka:peng:quen:97}.  
{To this end, we consider}  backward stochastic differential equations (BSDEs) 
 of the form
\begin{equation}
\label{eq:BSDE}
Y_t = \Phi(X_T) + \int_t^T f(s,X_s,Y_s,Z_s) ds - \int_t^T Z_s dW_s + N_{t,T}\,,\quad t\in[0,T],
\end{equation}
on a  filtered probability space $(\Omega, \cF_T, (\cF_t)_{t\in[0,T]},\P)$,  
satisfying the usual conditions  with finite horizon $T<\infty$  and a $q$-dimensional Brownian motion $W$.
The {\em terminal condition}
 $\Phi$ is a deterministic function 
satisfying some standard conditions (see Section~\ref{section:ass}), 
{while $X$ is} 
an exogenously given Markov process with 
 fixed initial value $X_0=x_0$, 
and $(N_{t,T})_{0\le t\le T}$ is a martingale orthogonal to {$W$}.
A solution to \eqref{eq:BSDE} is a {suitable} pair $(Y,Z)$ of $\R\times(\R^q)^\top$-valued processes.
Typically, BSDEs cannot be solved explicitly and one uses discrete time approximation.
Fixing a time-grid $\pi:= \{0 = t_0 , \ldots , t_N = T\}$, let $(X_i)_{0\le i \le N}$ be a suitable  discrete time approximation of $(X_t)_{0\le t \le T}$ on $\pi$.
We will build upon analysis in \cite{gobe:turk:13a}
on  so-called multi-step forward dynamical programming (MDP) equations
\begin{align}
\displaystyle Y_i & := \E[ \Phi(X_N) + \sum_{j=i}^{N-1} f_j(X_j,Y_{j+1}, Z_j)(t_{j+1} - t_j) | \cF_{t_i} ] , \quad\text{and} \label{eq:MDP} \\
\displaystyle \quad (t_{i+1} - t_i) \times Z_i& := \E[ (W_{t_{i+1}} - W_{t_i})^\top(\Phi(X_N) + \sum_{j=i+1}^{N-1} f_j(X_j,Y_{j+1}, Z_j)(t_{j+1} - t_j) ) | \cF_{t_i} ];
\nonumber 
\end{align}
{this process} $(Y_i,Z_i)$ ($i=N,\ldots,0$)  is {called} 
the {\em discrete BSDE} solution.
Further, we make use of a known splitting technique to decompose 
 the discrete BSDE 
 into the sum  $(Y,Z)= (y + \bar y, z + \bar z)$  of a the components of a  system of two (discrete) BSDEs given by
\begin{align}
\label{eq:lin}
 y_i  &:= \E[ \Phi(X_N) | \cF_{t_i} ]   \quad \text{and}\quad  (t_{i+1} - t_i) \times z_i := \E[ (W_{t_{i+1}} - W_{t_i})^\top \Phi(X_N) | \cF_{t_i} ] 
\end{align} 
{ for $i=N,\ldots, 0$, and, likewise, }     
\begin{align}
\label{eq:nonl}
 \displaystyle \bar y_i & := \E\left[  \sum_{j=i}^{N-1} f_j(X_j,y_{j+1} + \bar y_{j+1}, z_j + \bar z_j)(t_{j+1} - t_j)\,\Big|\, \cF_{t_i}\right]  \quad\text{and} \\
\nonumber \displaystyle  (t_{i+1} - t_i) \times \bar z_j & := \E\left[ (W_{t_{i+1}} - W_{t_i})^\top \sum_{j=i}^{N-1} f_j(X_j,y_{j+1} + \bar y_{j+1}, z_j + \bar z_j)(t_{j+1} - t_j) \Big| \cF_{t_i} \right].
\end{align}

We call the system (\ref{eq:lin} - \ref{eq:nonl}) for $(y,z)$ and $(\bar y,\bar z)$ the splitting scheme.
By adding together the equations \eqref{eq:lin} and \eqref{eq:nonl}, one recovers the original discrete BSDE \eqref{eq:MDP}.
In order to solve the system (\ref{eq:lin} - \ref{eq:nonl}), one must first solve for $(y,z)$ and then use that solution to solve $(\bar y, \bar z)$.
In general, one must approximate the conditional expectation operator to obtain a fully implementable algorithm, and for this we will make use of Monte Carlo least-squares regression, a method initiated in the BSDE context by \cite{gobe:lemo:wari:05}.
Our method will be to develop a novel multi-grid algorithm, which we term the multilevel algorithm, in order to efficiently approximate $(y,z)$, 
then to use the so-called least-squares multistep-forward dynamical programming algorithm (LSMDP) \cite{gobe:turk:13a} to approximate $(\bar y, \bar z)$.
In this paper, we focus on the error between the solution of (\ref{eq:lin} - \ref{eq:nonl}) and our fully implementable scheme; this is in the spirit of \cite{lemo:gobe:wari:06}\cite{gobe:turk:13a}. 
{The present paper is} not concerned with {the} error {from}  using {time-}discret{izing} 
schemes, like \eqref{eq:MDP}, to 
{approximate} \eqref{eq:BSDE};
for 
{analysis of this error},
{one can refer to extensive research}
\cite{zhan:04,bouc:touz:04,bouc:elie:08,gobe:laba:07,gobe:makh:10,imke:dosr:10,rich:11,hu:nual:song:11,chas:cris:14,%
geis:geis:gobe:12,turk:13b,dosr:lion:szpr:13,chas:rich:13}.

Thanks to improved regularity properties, the LSMDP algorithm  can be solved much more efficiently for \eqref{eq:nonl} than for \eqref{eq:MDP}, particularly in high dimension.
We show in Section \ref{section:conclusion} that, typically, resolving \eqref{eq:nonl} with LSMDP incurs a complexity of 
$O(\varepsilon^{-4-d/2} \ln(\varepsilon^{-1}+1)^d)$, 
whereas solving \eqref{eq:MDP} with LSMDP incurs a complexity of 
$O(\varepsilon^{-4-d} \ln(\varepsilon^{-1}+1)^d)$, 
where $\varepsilon$ is the precision and $d$ is the dimension of $X$.
The multilevel algorithm sequentially builds approximations of $(y,z)$ on a refining sequence of dyadic time-grids $\{ \pik k \ : \ k\ge0\}$ and takes the form of an adaptive martingale control variates algorithm:
assuming we have already constructed the the solution to \eqref{eq:lin} on the time-grid $\pik {k-1} := \{ 0 = \tk {k-1}0, \ldots , \tk {k-1}{2^{k-1}} = T\}$, which we denote $(\yn {k-1}{}, \zn {k-1}{})$, we use it to construct the solution on $\pik k := \{ 0 = \tk {k}0, \ldots \tk {k-1}{2^{k}} = T\}$ as follows:
\begin{align}
\displaystyle
 \yn ki & := \E[ \Phi(\X k{2^k}) - \sum_{j = \alpha( i) + 1}^{2^{k-1} -1} \zn {k-1}j (W_{\tk {k-1}{j+1}} - W_{\tk {k-1}j} ) | \cF_{\tk ki}], 
 \label{eq:ml}
\\
 \displaystyle
(\tk k{i+1} - \tk ki) \times \zn ki & := \E[ (W_{\tk {k}{i+1}} - W_{\tk {k}i} )^\top \big( \Phi(\X k{2^k}) -\yn ki - \sum_{j = \alpha( i) + 1}^{2^{k-1} -1} \zn {k-1}j (W_{\tk {k-1}{j+1}} - W_{\tk {k-1}j} )\big) | \cF_{\tk ki}],
\nonumber
\end{align}
for $\alpha( i) := \max \{0\le j \le 2^{k-1} \ : \ \tk {k-1}j \le \tk ki \} $ and $i \in \{ 0,\ldots,2^k-1\}$.
In order to solve \eqref{eq:ml}, one must first solve $\yn ki$, then $\zn ki$, then iterate the procedure for $i-1$; once one has solved $\yn k0$ and $\zn k0$, one may proceed to the time-grid $\pik {k+1}$.
Observe that, under the conditional expectation operator, the multilevel scheme \eqref{eq:ml} matches the discrete BSDE \eqref{eq:lin} whenever $\pi = \pik k$.
However, when the conditional expectation is replaced by the Monte Carlo least-squares operator, 
the multilevel formulation suffers from substantially less variance than the LSMDP formulation.
Indeed, we demonstrate in Section \ref{section:ML:err an} that the  complexity of the multilevel algorithm is typically
$
O(\varepsilon^{-2 - d} \ln(\varepsilon^{-1}+1)) 
$
whereas the LSMDP algorithm for \eqref{eq:lin} incurs a complexity of 
$
O(\varepsilon^{-3 - d} )
$; 
here, $\varepsilon$ is the precision and $d$ is the dimension of $X$.
We see that we have an order one improvement (up to log terms) 
by using the multilevel algorithm, which is substantial.
The overall complexity of using the splitting scheme with multilevel to approximate (\ref{eq:lin} - \ref{eq:nonl}) is 
\[
O(\varepsilon^{-2 - d} \ln(\varepsilon^{-1}+1)) + O(\varepsilon^{-4-d/2} \ln(\varepsilon^{-1}+1)^d)
\]
which should be compared to the complexity of the LSMDP scheme $O(\varepsilon^{-4-d} \ln(\varepsilon^{-1}+1)^d)$ for \eqref{eq:MDP}, so we see  a substantial overall gain in the complexity of the algorithm.
The reduction of the complexity is largely because one needs to generate fewer simulations of the process $X$.
This has the secondary  effect that it reduces the memory needed  to run the algorithm.
Since we are typically working with high dimensional problems (e.g., $d\ge 5$), the memory usage is typically very high, therefore reducing the memory usage is extremely important for practical implementation.

To conclude the introduction, we summarise the novelty of our results and compare them to the existing literature.
The majority of this paper is dedicated to the analysis of the multilevel scheme, which is, to the best of our knowledge, the first adaptive multi-grid algorithm for variance reduction in the approximation of BSDEs.
We make great efforts to keep make our results applicable in high generality and give examples of many  situations  {of interest} in which our assumptions are valid, see Section \ref{section:ass} - including some instances of discontinuous or path-dependent Markov process $X$.
The applicability of our results is not exclusive to the examples we present.
We  mention that the  splitting scheme (\ref{eq:lin} - \ref{eq:nonl}) has been studied in the literature in the past.
In the continuous time setting, it was used by \cite{gobe:makh:10}\cite{turk:13b} in order to determine regularity properties of the BSDE \eqref{eq:BSDE} when the filtration $(\cF_t)_{t \ge 0}$ is  generated by the Brownian motion $W$.
It has been used in the discrete time setting  by \cite{bend:stei:12} to design a numerical scheme based on a martingale basis technique, and in \cite{gobe:turk:13a} for {a} proxy scheme for BSDEs.
The {common} idea for both {techniques}  is  {to make efficient use of some a-priori knowledge} for  good approximation of the  {solution} $(y,z)$ {to the } discrete BSDE \eqref{eq:lin}, 
{if such is available e.g. though analytic knowledge of a suitable martingale basis or of approximate PDE solutions.} 
{In comparison},  the multilevel scheme does { not} {require such a-priori knowledge}. 
We demonstrate that obtaining the approximation of $(y,z)$ is the most expensive part {in} approximating the splitting scheme {\em without} multilevel; 
{this explains the overall efficiency gains that can be obtained by a (generic) multilevel aproximation of $(y,z)$.}  
We present explicit error estimates for our algorithm and demonstrate in a quantitative manner that we are able to obtain substantial complexity improvements; we also  provide numerical examples to corroborate these claims.
In the latter part of this paper, we determine explicit error estimates for the splitting scheme with multilevel used to approximate the solution of (\ref{eq:lin} - \ref{eq:nonl}); $(\bar y, \bar z)$ is approximated with an LSMDP scheme once $(y,z)$ is computed with the multilevel scheme.

We also use results  on the improved regularity \cite{turk:13b} of $(\bar y, \bar z)$ {from} \eqref{eq:nonl}, compared to $(Y,Z)$ {from} \eqref{eq:MDP}, to demonstrate that one can obtain better 
complexity {because one can choose a lower dimensional regression basis}.

{\bf Organization of the paper:}
Section \ref{section:notation} provides 
{some} notation used 
{within} paper.
 In Section \ref{section:ass}, we state the assumptions to be used throughout the paper, and 
{give several  examples to show that these assumptions permit  }
our algorithm to apply
in high generality.
In Section \ref{section:MC}, we present the multilevel scheme and compute explicit error estimates for the fully implementable scheme. 
These error estimates are then used to perform a 
complexity analysis that demonstrates a theoretical improved efficiency of the multilevel scheme compared to the LSMDP scheme.
In Section \ref{section:nonzero}, we perform the error analysis for the LSMDP scheme used to approximate $(\bar y, \bar z)$ \eqref{eq:nonl} given that $(y,z)$ \eqref{eq:lin} has been computed using the multilevel scheme.
Finally, {in} Section \ref{section:conclusion}, 
{we present} 
a complexity analysis of the splitting scheme with multilevel,
 and compare it to the LSMDP scheme 
 with
and without splitting 
 to demonstrate the efficiency gains.
Numerical example{s are included} throughout to demonstrate the improved efficiency of the fully{-implementable} multilevel algorithm in 
{actual computations}.

\subsection{Notation and conventions}
\label{section:notation}

For a given probability space $(\Omega,\cF,\P)$ and sub-$\sigma$-algebra $\cG$, 
we write $\L_2(\cG)$ for the space of $\cG$-measurable, square integrable random variables.
Constants are always understood to be finite and non-negative.
Filtrations in continuous time are taken to satisfy the usual condition of right continuity and completeness. 
Markov processes and semimartingales  in continuous time are taken to have cadlag paths. 
Inequalities between random variables (cadlag processes) are understood to hold  almost everywhere with respect to $\P$ ($\P\otimes dt$).
For any vector or matrix $V$, we denote its transpose by $V^\top$.
The usual Euclidian norm on some $\R^n$ (or $\R^{m\times n}$)
is denoted by  $|\cdot|$.
For any functions $f : \R^k \to \R^n$, the supremum norm is denoted by $|f(\cdot)|_\infty := \sup_{x\in\R^k} |f(x)|$.
For $L \ge 0$ and $l\in\N$, we define the truncation function $\cT_L : \R^l \to \R^l$ by $\cT_L(x) := (-L\vee x_1 \wedge L, \ldots, -L\vee x_l \wedge L)$
The multilevel approach is working along a refining sequence of time-grids and a sequence of approximating processes
evolving on those.
To this end, we introduce the following notation.
For each $k\ge0$, we denote by $\pik k := \{\tk k0, \ldots , \tk k{2^k}\}$ a time-grid with $2^k$ time-points, $\Del k_i := \tk k{i+1} - \tk ki$ the $(k+1)$-th time increment
and $\DW ki := W_{\tk k{i+1}} - W_{\tk ki}$ the $k$-th level $(i+1)$-th Brownian increment.
To deal with the the referencing of time indicies between the multiple time grids $\pik k$,  we define functions $\alpha : \{0,\ldots, 2^{k}\} \rightarrow \{0, \ldots, 2^{k-1}\}$  by
\begin{equation}
\alpha(i) := \alpha^{(k)}(i):=\max \{ j \in \{0,\ldots , 2^{k-1} \} \ : \ \tk {k-1} j \le \tk ki \}
\label{eq:alpha}
\end{equation}
To ease notation, we simply write $\alpha$ for $\alpha^{(k)}$ whenever the level $k$ is clear from the context. 
For the $\sigma$-algebras
$\F ki := \cF_{\tk ki}$ ($i=0,\ldots, 2^k$) we denote the respective conditional expectations by 
\(   
\E^k_i[\cdot] := \E[\cdot|\F ki]\,.          
\)
A stochastic process $X=(X_t)_{t\in [0,T]}$ that is piecewise constant with nodes at the time points of $\pik k$ 
we call discrete and write $ \X ki:= X_{\tk ki}$. In straightforward way, any $(\F ki)$-adapted process can be seen as a cadlag process in continuous time.
 We say that $X$ is $(\F ki)$-adapted if $X_i$ is $\F ki$-measurable  for each $i \in \{0,\ldots, 2^k\}$,
and call it  an $(\F ki)$-martingale if it is a martingale in the filtration $(\F ki)$.
Finally, a $\pik k$-Markov chain is a {discrete  process which is a  $(\F ki)$-Markov} chain.

\section{Assumptions}
\label{section:ass}
In this section, we state the conditions on the Markov processes/chains $X$, the time-grids $\pik k$, the terminal function $\Phi$ and the driver $f$ for the paper. 
{As in \cite{gobe:turk:13a}, we strive for a high level of generality under which the subsequent analysis is valid, 
to make results applicable to as wide a class of problems as possible. This includes but should not be restricted to the concrete 
examples of relevant practical problems, that are detailed in Section \ref{section:ass:egs} to explain and illustrate our general assumptions, which might appear overly abstract at first sight.}   
Section \ref{section:0:general} derives elementary consequences for representations and a-priori estimates that 
will be useful 
in what follows.
{Some additional assumptions that will be required only for the analysis in Sections \ref{section:nonzero} and \ref{section:conclusion}
will be detailed in
Section \ref{ass:two}}.

\subsection{General  assumptions}
\label{ass:gen}
The following conditions will hold throughout the entirety of this paper.
\begin{enumerate}
\item[\Hg]  The function $\Phi : \R^d\rightarrow \R$  
is measurable and is uniformly bounded by $C_\Phi$. 

\item[\HX]
There is a family of $({\cal F}_t)$-Markov processes $\{ (\X{t,x}{s})_{s\in[t,T]} \ : \ (t,x) \in [0,T]\times \R^d \}$ which share the same (possibly time-inhomogenous) Markov dynamics, in the sense that
 for the same semigroup of contraction operators $P_{t,s}$ ($t\le s$) acting on bounded measurable functions $h:\R^d\to \R$  it holds 
\begin {align}
P_{t,s} h(x)= E[h(\X{t,x}{s})] \quad \text{and}\quad P_{t,s} h(\X{u,x}{t})= E[h(\X{u,x}{s}|{\cal F}_t] 
\quad \text{for $0\le u\le t\le s\le T$.}
\end{align} 
We let $\X{t,x}{s} =x$ for $s \in [0,t]$.
This family satisfies the following properties:
\begin{itemize}
\item[(i)]
the Markov process $X$  is in this family and satisfies $X = \X {0,x_0}{}$;
\item[(ii)]
there is a  constant $C_X$ such that, for all $x,x'\in\R^d$ and $s \in [t,T]$,
\[
\E[ | \Phi(\X{t,x}T) -\Phi(\X{t,x'}T)|^2] \le C_X |x - x'|^2 \ \text{ and } \ \E[|\X {t,x}s - x|^2] \le C_X (t-s);
\]
\item[(iii)]
there exist deterministic functions $u : [0,T] \times \R^d \to \R$ and $v: [0,T] \times \R^d \to (\R^q)^\top$,  measurable,
such that, for any  $(t,x)$ in $[0,T]\times \R^d $,
the square integrable (bounded) martingale $\yn {t,x}{s}=E[\Phi(\X {t,x}{T})\,|\, {\cal F}_s]$ ($s\in [0,T]$)
and the predictable integrand $\zn{t,x} {}$ from the It\^o martingale representation
\begin{equation}
\label{eq:bsde:lin:tx}
 \yn{t,x}{s} = \Phi(\X{t,x}T) - \int_s^T \zn{t,x}{r} dW_r\,,\quad s\in[0,T],
 \end{equation}
 admit versions 
\(
\yn{t,x}{s} = u(s,\X{t,x}s) \quad \text{and} \quad \zn{t,x}{s} = v(s,\X{t,x}s).
\)
 \item[(iv)]
there exist constants $\theta \in (0,1]$ and $C_X \ge 0$ such that, for all $t\in[0,T)$ and $x\in\R^d$, $|v(t,x)|$ is bounded  by $C_X/(T-t)^{(1-\theta)/2}$. 
Moreover, the functions $u(t,\cdot)$ and $v(t,\cdot)$ are Lipschitz continuous with Lipschitz constants $C_X/(T-t)^{(1-\theta)/2}$ and $C_X/(T-t)^{1-\theta/2}$ respectively.

\end{itemize}

\item[\Ht] The set of time-grids $\{\pik k \ : \ k \ge 0\}$ satisfies 
\begin{itemize}
\item[(i)]
{$\pik{k+1}$ is refinement of $\pik k$; }

\item[(ii)]
there exists a  constant $C_X$ such that $\max_{0\le i \le 2^k-1} \Del k_i \le C_X 2^{-k}$;

\item[(iii)] 
there exists a  constant $c_X$ such that $\min_{0\le i \le 2^k-1} \Del k_i \ge c_X 2^{-k}$;

\item[(iv)]
recalling the processes $( \yn{t,x}{}, \zn{t,x}{})$ solving \eqref{eq:bsde:lin:tx}, there is a  constant $C_X$ such that for all $k\ge0$, $\tk ki \in \pik k$ and $x\in\R^d$,
\[
\sum_{j =0}^{2^k-1} \E[ \int_{\tk kj}^{\tk k{j+1}} |\zn {\tk ki, x}t - \tilz {k,i,x}{j}|^2 dt ] \le C_X 2^{-k},  \text{ where } \tilz {k,i,x}j := {1 \over \Del k_i} \E^k_i[ \int_{\tk kj}^{\tk k{j+1}} \zn {\tk ki, x}t dt].
\]
\end{itemize}

\item[\HXp]
There is a family of  $\R^d$-valued, $\pik k$-Markov chains  $\{\X{k, i,x}{} \ : \ x\in\R^d, 0\le i \le 2^k, \ \X{k,i,x}j =x \ \forall \ j \le i \}$
satisfying the properties:
\begin{itemize}
\item[(i)]
recalling the parameter $\theta$ from \HX(iv), there is a  constant $C_X$ such that, 
for all $x\in\R^d$, $k \ge 0$ and $i \in\{0,\ldots,2^k\}$,  $\E[|\Phi(\X{\tk ki,x}T) - \Phi(\X{k,i,x}{2^k}) |^2] \le C_X 2^{-k}$ {and $\E[| \Phi(\X{k,2^k-1,x}{2^k}) - \Phi(x) |^2] \le C_X 2^{-2\theta k}$};

\item[(ii)]
for all $k\ge 0$ and $j \in \{0,\ldots, 2^k\}$ and $l \in \{j,\ldots,2^k\}$, there exist a $\cG^{(k)}_{j}\otimes\cB(\R^d)$-measurable functions $\Xmarkov k{j,l}  : \Omega \times \R^d \rightarrow \R^d$
with a $\sigma$-algebra $\cG^{(k)}_{j} \subset \cF_T$  independent of $\cF_{t_j}$  and containing $\sigma(\DW kr  :  r \ge j)$,   such that $\X{k,i,x}{l} = \Xmarkov{k}{j,l} (\X{k,i,x}{j})$;

\item[(iii)]
there exists a  constant $C_X$ such that, for all $x\in\R^d$, $k \ge 1$, $i \in \{0,\ldots ,2^k\}$, and $j \in\{i,\ldots,2^k\}$
\[
\E[| \X {k,i,x}j - \X {k-1,\alpha(i),x}{\alpha(j)}|^2] \le C_X 2^{-k}.
\]
\end{itemize}
{For brevity}, we denote the process $(\X {k,0,x_0}{i})_{0\le i \le 2^k}$ by $(\X k{i})_{0\le i \le 2^k}$.

\end{enumerate}

\subsection{Properties of the discrete BSDE derived from the assumptions} 
\label{section:0:general}

This section  
{collects some} elementary results that {follow}
directly from the general conditions 
{in} Section \ref{ass:gen}
{and are useful} in Sections \ref{section:MC} -- \ref{section:conclusion}.
Recall the family of Markov processes 
from \HX. 
{For every  $i \in\{ 0,\ldots ,2^k\}$ and $x\in\R^d$,} the Kunita-Watanabe decomposition 
guarantees
{(unique)}
existence of a 
pair 
 $(\tily{k}{},\tilz{k}{})=(\tily{k,i,x}{},\tilz{k,i,x}{}) $ of square integrable, $\pik k$-adapted processes
and a square integrable $(\F ki)$-martingale $L^{(k)}$ such that
\begin{equation}
\tily kl = \Phi(\X{\tk ki,x}T) - \sum_{j=l}^{2^k-1} \tilz kj \DW kj  - \sum_{j=l}^{2^k-1} \DL kj, \qquad l \ge i,
\label{eq:disc:bsde:0}
\end{equation}
where $\DL ki := L^{(k)}_{i+1} - L^{(k)}_i$, $L^{(k)}_0 =0$,  and $L^{(k)}$ is {(strongly) orthogonal} to $W^{(k)}$ in the sense that
{$\left(W^{(k)}_j L^{(k)}_j\right)_{0\le j\le N}$ is an  $(\F ki)$-martingale, i.e.}
{\(\E^k_j[\DW kj \DL kj] = 0 \quad \text{for all } j$. 
}
We will determine an explicit representation for $(\tily k{}, \tilz k{}, L^{(k)})$ 
in terms the solution of the continuous time BSDE $(\yn {\tk ki,x}{},\zn {\tk ki,x}{})$ {from \HX(iii)}
and the conditional expectation $\E^k_l[\cdot]$ {in the next lemma.} 
{This permits one to establish}
important a-priori bounds on the 
{process } $\tilz kl$ {by Corollary \ref{cor:bd:z}}.
\begin{lemma}
\label{lem:disc:bsde:0}
For any time-grid $\pik k$, $k\ge0$, $i \in \{0,\ldots , 2^k\}$, $l \ge i$, and $x\in\R^d$ {holds}
\begin{align}
\tily kl =  &\E^k_l[\Phi(\X{\tk ki,x}T)] \quad \text{and} \quad \Del k_l \tilz kl = \E^k_l[ (\DW kl)^\top\Phi(\X{\tk ki,x}T)], 
\label{eq:bsde:0:cdn} \\
\DL kl  =  &\int_{\tk kl}^{\tk k{l+1}} (\zn {\tk ki,x}s - \tilz kl)dW_s,
\label{eq:l:inc:0} \\
\Del{k}_l\tilz kl   =  &\E^k_l[ \int_{\tk kl}^{\tk{k}{l+1}} \zn{\tk ki, x}t dt ] .
\label{eq:z:integral}
\end{align}
\end{lemma}
{\bf Proof.}
{Equalities in \eqref{eq:bsde:0:cdn}  are}
{well known and easily obtained by} 
taking conditional expectations 
{of} \eqref{eq:disc:bsde:0} 
{itself or in product with} $(\DW kl)^\top$.
{Equality $\yn{\tk ki,x}{\tk kl} =\E^k_l[\Phi(\X{\tk ki,x}T)]$} 
implies
\begin{align*}
\Y kl = \yn {\tk ki, x}{\tk kl}  = \Phi(\X{\tk ki,x}T) - \sum_{j=l}^{2^k-1} \tilz kj \DW kj - \sum_{j=l}^{2^k-1} \int_{\tk kj}^{\tk k{j+1}} (\zn {\tk ki,x}s - \tilz kj) dW_s
\label{eq:lem:bsde:0:1}
\end{align*}
{and thereby $\DL kj=\int_{\tk kj}^{\tk k{j+1}} ( \zn {\tk ki,x}s  - \tilz kj) dW_s$.} 
{Equality \eqref{eq:z:integral}  follows easily. Indeed,} 
multiplying with $(\DW ki)^\top$ and taking conditional expectation $\E^k_l[\cdot]$
{yields}  
{\(0 = \Del{k}_l \tilz kl
- \E^k_l[ \int_{\tk kl}^{\tk{k}{l+1}} \zn {\tk ki,x}t  dt ]\).}
Moreover,
\[
0 =   \E^k_l[\Phi(\X{\tk ki,x}T) (\DW kl)^\top] 
      - \E^k_l[(\DW ki)^\top \int_{\tk kl}^{\tk k{l+1}} \zn {\tk ki,x}t dW_t] 
   = \Del{k}_l \tilz kl  
      - \E^k_l[ \int_{\tk kl}^{\tk{k}{l+1}} \zn {\tk ki,x}t  dt ].
\]
which proves \eqref{eq:z:integral}.
\qed

{Naturally, one obtains} {a Markov representation for the solutions of the discrete BSDEs (\ref{eq:lin} - \ref{eq:nonl}).

\begin{lemma}
\label{lem:markov:bsde}
For all $k\ge 0$ and $j \in \{0,\ldots, 2^k-1\}$, there exist deterministic functions $\yn kj : \R^d \to \R$ and $\zn kj : \R^d \to (\R^q)^\top$ such that 
\[
\displaystyle \E^k_j [\Phi(\X{k,i,x}{2^k})] = \y kj{\X{k,i,x}j} \quad \text{and} \quad {1 \over \Del k_j} \E^k_j [(\DW kj)^\top\Phi(\X{k,i,x}{2^k})] = \z kj{\X{k,i,x}j} 
\]
for all $i \in \{0,\ldots, 2^k\}$ and $x\in\R^d$.
Moreover, there exist deterministic functions
 $\bar y^{(k)}_j : \R^d \to \R$ and $\bar z^{(k)}_j : \R^d \to (\R^q)^\top${, with $\bar y^{(k)}_{2^k}(\cdot)= 0$,} { for $ j\in \{0,\ldots, 2^k-1\}$}  
\begin{align*}
& \E^k_j[  \sum_{l=j}^{N-1} f_l(\X{k,i,x}l,\y k{l+1}{\X{k,i,x}{l+1}} 
      + \bar y^{(k)}_{l+1}(\X{k,i,x}{k+1}),  \z kl {\X{k,i,x}l} + \bar z^{(k)}_k(\X{k,i,x}l))\Del k_l] 
       = \bar y^{(k)}_j(\X{k,i,x}j)    \\
 & \quad \text{and} \quad  \E^k_j [{(\DW kj)^\top  \over \Del k_j} \left(\sum_{l=j+1}^{N-1} f_l(\X{k,i,x}l,\y k{l+1}{\X{k,i,x}{l+1}} + \bar y^{(k)}_{l+1}(\X{k,i,x}{l+1}), \z kl {\X{k,i,x}l} + \bar z^{(k)}_l(\X{k,i,x}l))\Del k_l \right)] \\ 
 &  \qquad \qquad= \bar z_j(\X{k,i,x}j) \qquad\text{for all $i \in \{0,\ldots, 2^k\}$ and $x\in\R^d$.}
\end{align*}
\end{lemma}
}
{This follows directly from \HXp(ii) by routine conditioning arguments from measure theory, like {\cite[Lemma 4.1]{gobe:turk:13a}}.}

Finally, we present almost sure {absolute} bounds, uniform in $x$, for the functions $\y kix$ and $\z kix$.
{Such} bounds are crucial for {and repeatedly used} 
in Section~\ref{section:MC}.

\begin{corollary}
\label{cor:bd:z}
There exists a constant $C_{X}$ such that, for all $k \ge 0$, $i \in \{0,\ldots , 2^k\}$, $l \ge i$, and $x\in\R^d$,
 {$|\tilz {k}l| = |\tilz {k,i,x}l|$}
 is bounded {(a.s.)} 
 by 
\begin{equation}
\Cz{k,l} := {C_X\over (T- \tk kl)^{(1-\theta)/2} } 
\end{equation}
This implies that, for all $k \ge 0$, the function $\z ki\cdot$ is absolutely bounded by $\Cz{k,i}$.
Further, there exists { a constant} $C_y $
independent of $k$ and $i$, such that 
{$|\y ki\cdot|$} is bounded 
by $C_y$ (a.s.).
\end{corollary}
{\bf Proof.}
Recall the version $v(t,\X{\tk ki,x} t)$ of $ \zn {\tk ki,x}t $ from \HX(iii), and the absolute bound on the function $x\to v(t,x)$ from \HX(iv).
Using the representation \eqref{eq:z:integral}, it follows that
\begin{align*}
&|\tilz ki|  = {1 \over \Del{k}_i} \left| \E^k_i[ \int_{\tk ki}^{\tk{k}{i+1}} v(t,\X{\tk ki,x} t)  dt ] \right| 
 \;\le\; {C_X \over \Del{k}_i}  \int_{\tk ki}^{\tk{k}{i+1}}  {dt \over (T-t)^{(1-\theta)/2} }  \\
& = {2 C_X  (T-\tk ki) \over \Del k_i(1+\theta) (T-\tk ki)^{(1- \theta)/2} } - {2 C_X  (T-\tk k{i+1}) \over \Del k_i(1+\theta)  (T-\tk k{i+1})^{(1-\theta)/2}} 
\;\le\; {2 C_X  \Del k_i \over \Del k_i(1+\theta) (T-\tk ki)^{(1- \theta)/2} } .
\end{align*}
In order to obtain bounds for $\z kix $, we use additionally the condition \HXp(i) in order to obtain
\[
|\z kix | \le |\z kix - \tilz ki| + |\tilz ki| \le  \sqrt{\frac{1}{\Del k_i} \E[| \Phi(\X{k,i,x}{2^k} ) - \Phi(\X{\tk ki, x}T)|^2] } +  |\tilz ki| \le C_X +  |\tilz ki| 
\]
We mildly abuse notation by stating $|\z kix | \le C_X (T-\tk ki)^{-(1-\theta)/2}$ to 
{ease} 
notation.
The bound on $\yn ki$ is immediate from the boundedness of $\Phi(\cdot)$ in \Hg.
\qed

\subsection{Additional assumptions and properties for  Section \ref{section:nonzero} and \ref{section:conclusion}}
\label{ass:two}

The following conditions are used  in Sections \ref{section:nonzero} and \ref{section:conclusion}.

\begin{enumerate}
\item[\Hf] For every $i \in\{0,\ldots,2^k\}$, the driver $f_i : \R^d \times \R \times( \R^q)^\top \rightarrow \R$ 
is $\cB(\R^d)  \otimes \cB(\R) \otimes \cB((\R^q)^\top)$-measurable and satisfies the following properties:
\begin{itemize}
\item[({i)}]
for all $x\in\R^d$, 
$(y,z) \mapsto f_i(x,y,z)$ is Lipschitz continuous with $\tk ki$-dependant Lipschitz constant:
there exist constants $L_f>0$ finite and $\thetaL \in (0,1]$ such that, for all  $(y,z),(y',z') \in  \R\times(\R^q)^\top$,
\begin{equation*}
|f_i(x,y,z) - f_i(x,y',z')| \le{ L_f \{ |y-y'| + |z-z'| \} \over (T-\tk ki)^{(1-\thetaL)/2} };
\end{equation*}
\item[(ii)] 
${|}f_i(x,0,0) {|}$, is {uniformly} bounded by a 
 constant $C_f
 $: $|f_i(x,0,0)|\le C_f$ for all $x$;

\item[(iii)] recalling the functions $u$ and $v$, and the parameter $\theta$, from \HX(iii), there are deterministic functions $U : [0,T] \times \R^d \to \R$ and $V : [0,T] \times \R^d \to (\R^q)^\top$ such that there is a version of $(Y,Z)$ of the solution of the BSDE \eqref{eq:BSDE} satisfying 
\[
Y_t - u(t,X_t) = U(t,X_t), \qquad Z_t - v(t,X_t) = V(t,X_t)
\]
for all $t$ 
{(a.s.)}. 
Moreover, there exists a 
constant $C_X$ such that
{$|U(t,x)|$ is bounded by $C_X(T-t)^{\theta + \thetaL}$,}
 {and} $U(t,\cdot)$ is Lipschitz continuous with Lipschitz constant $C_X (T-t)^{ - (1-\thetaL-\theta)/2}$,
 and $V(t,\cdot)$ is Lipschitz continuous with 
 Lipschitz constant $C_X(T-t)^{-\{1-(\thetaL+\theta)/2\}}$ 
 {at} $t \in [0,T)$.
\end{itemize}

\item[\HXpp]
Recall the functions $\bar y^{(k)}_j : \R^d \to \R$ and $\bar z^{(k)}_j : \R^d \to (\R^q)^\top$   from Lemma \ref{lem:markov:bsde},
and the functions $U : [0,T] \times \R^d \to \R$ and $V : [0,T] \times \R^d \to (\R^q)^\top$ from \Hf(iii).
For all $x \in \R^d$ and $\tk kj \in \pik k$, there exists a constant $C_X$ such that 
\begin{align*}
&\max_{j\le i \le 2^k-1} \e^k_j[ |   \bar y^{(k)}_i (\X {k , j, x}i) - U(\tk ki , \X {\tk kj,x}{\tk ki} )|^2 ] \; \le\; C_X 2^{-k} , \\
 & \text{{and}}\quad \sum_{i=0}^{2^k-1} \e[ |   \bar z^{(k)}_ i(\X ki) - V(\tk ki, X_{\tk ki}) |^2 ] \Del k_i \;\le \;C_X 2^{-k}.
\end{align*}

\item[\Htp]  Recalling the parameter $\thetaL$ from \Hf(i), the time-grids
{$\pik k := \{0=t_0 < \ldots < \tk k{2^k} = T \}$, $ \ k\ge 0$,} are such that
\begin{align}
\label{eq:delta}
C_{\pik k}  :=\sup_{i<2^k}\frac{\Del k_i }{(T- \tk ki)^{1-\thetal} }\; \longrightarrow \;0 \quad \text{as }k\rightarrow +\infty, \quad\\
\label{eq:R}
\text{  {and}}\quad \limsup_{k \rightarrow \infty} \ R_{\pik k}  < +\infty, \quad \text{where } R_{\pik k} {:=}
\sup_{0 \le i \le 2^k-2} \frac{\Del k_i}{\Del k_{i+1}}. 
\end{align}

\end{enumerate}

{We now use \HXpp \ to prove  {a-priori} 
bounds for the functions $\bar y^{(k)}_j (x)$ and $ \bar z^{(k)}_j (x)$ from Lemma \ref{lem:markov:bsde},
 and consequently for the processes $\bar y^{(k)}_j $ and $\bar z^{(k)}_j $, 
 {similar to} Corollary \ref{cor:bd:z}.
These bounds will be crucial for constructing algorithms and obtaining error estimates in Section \ref{section:nonzero}.
}

\begin{lemma}
\label{cor:bd:bar:z}
There exists a  constant $C_X$
such that, for all $k \ge 0$, $i \in \{0,\ldots , 2^k - 1\}$, $x\in\R^d$ we have 
\begin{equation*}
|\bar y^{(k)}_i(x) | \le C_{y,i} := {C_X    (T-\tk ki)^{(\thetaL + \theta)/2} }, \qquad
| \bar z^{(k)}_i (x) | \le \Cz{i} :=   {C_X \over (T-\tk k{i})^{(1-\{ \thetaL\vee\theta + \theta\})/2}  } . 
\end{equation*}
\end{lemma}
{\bf Proof.}
Recall 
the function $U(t,x)$ from Assumption \Hf.
The bound on $\bar y^{(k)}_i(x)$ is obtained from the trivial decomposition
\(
|\bar y^{(k)}_i(x))| =  |\bar y^{(k)}_i(x) - U(\tk ki,x) | + | U(\tk ki, x)| .
\)
Then, using the bounds from \Hf(iii) and \HXpp(ii) (with $j = i$), the result on $| \bar y^{(k)}_i(x) |$ follows.
{By} a mild abuse of notation, {we replace } $ C_X (T + T^{\theta + \thetaL} )$ by $C_X$
{to simplify notation}.

{For} the bound on $| \bar z^{(k)}_i (x) | $, we recall that $\Del k_i  \bar z^{(k)}_i (x) =  \E ^k_i[  \bar y^{(k)}_{i+1} (\X {k,i,x}{i+1})\DW ki ]$
{and} treat the cases  $i =2^k-1$ and $i < 2^k-1$ separately.
For $i =2^k-1$, 
the Cauchy-Schwarz inequality and \HXp(i) {yield}
\[
(\Del k_{2^k-1})^2 \left|\bar z^{(k)}_{2^k-1} (x)\right|^2 = \left|\E\left[\DW k{2^k-1} \{ \Phi(\X {k,2^k-1,x}{2^k} ) - \Phi(x) \} \right]\right|^2 \le C_X (\Del k_{2^k-1})^{1+2\theta}
\]
 implying 
 $|\bar z^{(k)}_{2^k-1}(x) |^2 \le C_X (T- \tk k{2^k-1})^{-1 + 2\theta}$, as required.
For $i<2^k-1$, 
apply{ing} Cauchy-Schwarz 
{yields}
\begin{align*}
&(\Del k_i)^2 |\bar z^{(k)}_i (x)|^2  = \left| \E ^k_i\left[ \left\{ \bar y^{(k)}_{i+1}(\X {k,i,x}{i+1}) \pm U(\tk k{i+1}, \X{\tk ki, x}{\tk k{i+1}} ) - U(\tk k{i+1},x) \right\} \DW ki \right] \right|^2 \\
&\qquad \le 2 \Del k_i  \E ^k_i\left[ \big| \bar y^{(k)}_{i+1}(\X {k,i,x}{i+1}) - U(\tk k{i+1}, \X{\tk ki, x}{\tk k{i+1}} )   \big|^2 \right]   
  + 2 \Del k_i \E^k_i\big[|U(\tk k{i+1}, \X{\tk ki, x}{\tk k{i+1}} ) - U(\tk k{i+1},x)|^2\big].
\end{align*}
Using assumption \HXpp(ii) for the first term and the Lipschitz continuity of $U$ in \Hf(iii) for the second yields 
\begin{align*}
(\Del k_i)^2 |\bar z^{(k)}_i(x) |^2 
&\le 2 (\Del k_i )^2 C_X 
+ { 2 (\Del k_i )^2 C_X  \over  (T-\tk k{i+1})^{1-\{ \thetaL + \theta\}}   }
\le 2 (\Del k_i )^2 C_X 
+ { 2 (\Del k_i )^2 C_X  \over  (T-\tk k{i})^{1-\{ \thetaL + \theta\}}   },
\end{align*}
where one exchanges $(T-\tk k{i+1})$ by $(T-\tk k{i})$
from \Htp.
By mild abuse of notation, we {re}write $C_X :=
C_X\vee \sqrt {2C_X (1 + T^{1-\{\thetaL + \theta\}})}$ to simplify the result.
\qed

\subsection{Examples satisfying the general assumptions}
\label{section:ass:egs}

{This section details explicit examples of processes, time-grids and functions to illustrate and explain the conditions from Section~\ref{ass:gen}.}

\noindent
$\blacktriangleright${\bf Assumption \HX.}
Property (ii) is a Lipschitz continuity property of the payoff $\Phi(\X{t,x}T)$ with respect to the initial value $x$.
It is satisfied if (a) $\Phi$ is locally Lipschitz
 continuous,  i.e.\  for some constant $c$ and $l\in [0,\infty)$ {holds} $|\Phi(x) - \Phi(x')| \le c (|x|^l + |x'|^l) |x - x'|$ {for all $x,x'$};
 or (b) H\"older continuous with H\"older exponent greater than or equal to $1/2$, and,
 in both cases (a) and (b), $\X{t,x}{}$ solves an SDE (which may have jumps) with Lipschitz continuous coefficients.
{A lower H\"older regularity in case (b) would lower the convergence rate of the numerical scheme in Theorems \ref{thm:0:z:er}, \ref{thm:spec:basis}, cf.\ Remark \ref{rem:gen:hold}.}

Property (iii) is a classical property of Markovian BSDEs.
It is satisfied when $\X{t,x}{}$ solves an SDE with deterministic 
({Markovian}),
 Lipschitz continuous coefficient {functions} in a Brownian filtration \cite[Theorem 4.1]{elka:peng:quen:97} 
or in 
{a} L\'evy filtration \cite[Proposition 4]{nual:scho:01}.
The property is also known to hold in the setting where $\X{t,x}s$ is of the form $(S^{(t,x)}_{r_1\wedge s},\ldots, S^{(t,x)}_{r_l\wedge s}, S^{(t ,x)}_s )$, where $S^{(t,x)}$ is the solution of an SDE  in the Brownian filtration and $t \le r_1<\ldots < r_l \le T $, see \cite{ma:zhan:02, geis:geis:gobe:12}.

There are two important instances where one can show Property (iv) to be valid.
Firstly, suppose that $\X {t,x}{}$ solves an SDE with deterministic ({Markovian}), bounded and continuously differentiable coefficients,
 whose partial derivatives are bounded and H\"older continuous  
 {with} diffusion coefficient 
 {being} uniformly elliptic. 
Then $v(t,x) = (\sigma(t,x) \nabla_x u(t,x))^\top$ and the boundedness of $|v(t,x)|$ 
follows from classical gradient bounds {of} parabolic PDEs \cite{frie:64};
  {and}  $\theta$ is equal to the H\"older exponent of $\Phi$.
Secondly,
if $\Phi$ is locally Lipschitz continuous, this result holds with $\theta =1$ if $\X {t,x}{}$ solves an SDE with deterministic ({Markovian}), Lipschitz continuous coefficients having linear growth;
the path dependant setting $\X {t,x}s = (S^{(t,x)}_{r_1\wedge s},\ldots, S^{(t,x)}_{r_l\wedge s}, S^{(t ,x)}_s )$ - where $S^{(t,x)}$ solves an SDE with deterministic ({Markovian}), Lipschitz continuous coefficients having linear growth - is also valid in this setting.

\noindent $\blacktriangleright${\bf Assumption \Ht.} 
{
Condition (i) is  to ensure that $\tk {k-1}{\alpha(i)+1} > \tk ki$ for all $k$ and  $i$; 
later, when we introduce condition  \HXppp \ in Theorem \ref{thm:0:z:er}, this condition becomes crucial.
It is satisfied by the time grids with points $\tk ki := T - T(1 - i/2^k)^{1/\beta}$ for any $\beta \in (0,1]$, which includes the uniform time-grid.
}
{
Condition (iv) is the most complex of the requirements, and has been studied extensively in recent years.
Let us first consider the case of the Brownian filtration.
Then, condition (iv) is satisfied for (locally) Lipschitz continuous $\Phi$ and uniform time-grids if $\X {\tk ki, x}{}$ is the solution of an SDE 
with deterministic ({Markovian}), Lipschitz continuous coefficients 
{of} linear growth \cite{zhan:04,rich:12}
 local Lipschitz continuity is {meant as described in}
 \HX.
For H\"older continuous (fractionally smooth) $\Phi$, it is satisfied by the time-grids with points $\tk ki := T - T(1 - i/2^k)^{1/\beta}$
if $\beta $ is less than the H\"older (fractional smoothness) exponent of $\Phi$, and
$\X {\tk ki,x}{}$ solves a continuous SDE with deterministic ({Markovian}), bounded 
and twice continuously differentiable coefficients $b(t,x)$ for the drift and $\sigma(t,x)$ for the volatility,
 whose partial derivatives are bounded and H\"older continuous, and $\sigma$ is uniformly elliptic \cite{gobe:makh:10};
note that the time-grid also satisfies properties (i)-(iii), cf.\ \cite[Lemma 5.3]{turk:13b} for {a} proof of (ii).
We remark that this rate of convergence may not be optimal, cf.\ 
\cite{gobe:laba:07}\cite{laba:07}.
The path dependent setting $\X {\tk ki,x}s = (S^{(\tk ki,x)}_{r_1\wedge s},\ldots, S^{(\tk ki,x)}_{r_l\wedge s}, S^{(\tk ki ,x)}_s )$ 
with $\Phi$ fractionally smooth and $S^{(t,x)}$ {being} the solution of an SDE with bounded, twice differentiable coefficients, 
whose partial derivatives are bounded and H\"older continuous, 
also satisfies the condition 
{if suitable} time grids are used; 
 cf.\  \cite{geis:geis:gobe:12}.
In a filtration generated by a L\'evy process, \cite{bouc:elie:08} showed that the uniform time grid was sufficient to have this property
 if the terminal condition is of the form $\Phi(X_T)$, for $\X {\tk ki,x}{}$ solving an SDE with Lipschitz continuous coefficients
{of} linear growth, and $\Phi$ is Lipschitz continuous.
}

\noindent $\blacktriangleright${\bf Assumption \HXp.}
Condition (i) is a ``good-approximation" criterion for the Markov process by the Markov chain. 
It is satisfied if $\Phi$ is (locally) Lipschitz continuous
and $\X{\tk ki,x}{}$ solves an SDE (with jumps) whose coefficients are deterministic ({Markovian}), 
Lipschitz continuous and have linear growth;  
$\X{k,i,x}{}$ may be the Euler scheme approximation of $\X{\tk ki,x}{}$ on the time-grid $\pik k$.
{However,  if  the terminal condition has  a lower regularity,} the Euler scheme 
{might} not satisfy this condition;
 for example,   in the case where $\Phi$ has only bounded variation, see \cite[Theorem 5.4]{avik:09}.
Higher order approximation schemes may be required for H\"older exponent $\theta$ less than $1$. 

Condition (ii) is slightly stronger requirement on the Markov chain than the basic definition; 
it is satisfied by   most approximation schemes for SDEs, including the Euler scheme.

Condition (iii) is a typical estimate required in multilevel Monte Carlo type approximation schemes for SDEs,
{cf.}\ \cite{gile:08,gile:szpr:12,gile:szpr:13a} and references therein. 
It is a property satisfied, for instance, by the Euler scheme for an SDE with deterministic ({Markovian}), Lipschitz continuous coefficients
of linear growth.
{If a convergence rate for the Markov chains 
{ were lower},
a lower rate of convergence of the BSDE multilevel scheme would be obtained; see Remark \ref{rem:coup}.
}

{
\noindent $\blacktriangleright${\bf Assumption \Hf.} 
The condition \Hf(i) of  Lipschitz continuous driver is  standard in the literature \cite{lemo:gobe:wari:06,bend:denk:07,bouc:elie:08,bria:laba:14} for $\thetaL =1$, and has more recently been extended to the setting $\thetaL<1$ \cite{gobe:turk:13a,gobe:turk:13b,turk:13b}. The {case} $\thetaL <1$ allows to treat some cases of quadratic BSDEs \cite{gobe:turk:13a}.

The Lipschitz continuity of $V(t,\cdot)$ condition  \Hf(iii)  is available, for example, from \cite[Corollary 4.3]{turk:13b}. 
This result is for a Brownian filtration, where $\X{t,x}{}$ solves an SDE with deterministic ({Markovian}), bounded, twice differentiable coefficients whose partial derivatives are bounded and H\"older continuous, and whose volatility matrix is uniformly elliptic. 
The Lipschitz continuity of the function $V(t,\cdot)$ is equal to $\lim_{\varepsilon \to 0} \phi(t,\varepsilon,\thetaL,\theta) = C_X/(T-t)^{1-(\thetaL -\theta)/2}$ for all $t\in[0,T)$; note that $\theta$ is denoted $\thetaP$ in that work.
The convergence requires the condition  that $\thetaL + \theta \ge 1$. 
The estimate of the Lipschitz constant of $U(t,\cdot)$ comes from the standard result that there exists a constant independent of $(t,x)$ such that $|\nabla_x U(t,x)| \le C |V(t,x)|$ for all $(t,x)$, and that $|V(t,\cdot)|_\infty \le C_X /(T-t)^{(1-\thetaL-\theta)/2}$  \cite[Corollary 2.13]{turk:13b}.
Likewise, the almost sure bound on $U(t,x)$ is available in \cite[eq. (3.12)]{turk:13b}.
}

{
\noindent $\blacktriangleright${\bf Assumption \HXpp.}
This condition is a discretization property; it has been proved under quite  general conditions in for example, \cite[Section 3]{gobe:makh:10}\cite[Section 3]{turk:13b}. Both references treat $X$ which solves an SDE with deterministic ({Markovian}), bounded, twice differentiable coefficients whose partial derivatives are bounded and H\"older continuous, and whose volatility matrix is uniformly elliptic;
the terminal conditions can be fractionally smooth.
}

\noindent $\blacktriangleright${\bf Assumption \Htp.} 
The additional conditions on the time grid \Hf(iv) are required to obtain a-priori estimates in the discrete setting; see \cite[Prop. 3.2]{gobe:turk:13a}.
Although the conditions seem  abstract, they are in fact satisfied by the time-grids given as examples in \Ht; see \cite[Lem.B.1]{turk:13b}.

\section{Multilevel least-squares Monte Carlo scheme}
\label{section:MC}
In this section, we construct approximations the functions $\y kix$ (and $\z kix$) {in Lemma \ref{lem:markov:bsde}}
for each level $k$ of the multilevel algorithm and each time-point $i$ of the grid $\pik k$.
The approximating functions are denoted by
\begin{equation}
\yn {k,M}i : \R^d \rightarrow \R, \qquad\text{respectively } \zn {k,M}i: \R^d \rightarrow (\R^q)^\top.
\label{eq:MC:proxys}
\end{equation}

The multilevel algorithm uses Monte Carlo least-squares regression to approximate conditional expectations
{and will be described in} Section \ref{section:ML:algorithm}.
{We apply ordinary least-squares regression, as in} \cite{gobe:turk:13a},
{whose terminology we recall } in Section \ref{section:ML:prelims}.
 {In contrast to} \cite{gobe:turk:13a}, 
 {we have} to pay {extra} attention to the
{novel} multilevel structure 
of the algorithm.
In Section \ref{section:ML:err an}, a comprehensive error analysis
gives an upper bound for the global error
\begin{equation}
\max_{0\le i \le 2^k -1} \E[|\y k i{\X ki} - \y{k,M}i{\X ki}|^2] + \sum_{i=0}^{2^k-1} \E[|\z k i{\X ki} - \z{k,M}i{\X ki}|^2] \Del{k},
\label{eq:ML:str er} 
\end{equation}
{and shows how it depends on}
  {numerical} parameters
 (the number of Monte Carlo simulations, the {choice} of basis functions)  and global error on the level $k-1$ of the algorithm. 
This error analysis 
{enables us}, 
{in} Section~\ref{section:ML:err an},  
to calibrate the numerical parameters of the {multilevel} algorithm 
 and to compare 
 {its complexity to that of } 
 {alternative} algorithms.

\subsection{Preliminaries}
\label{section:ML:prelims}
This section
introduces ordinary least-squares regression (OLS)
 to approximate the conditional expectation operator in the multilevel scheme. 
We will build on {a general but versatile}
 Definition~\ref{def:ls}
  for OLS
  to express 
our 
algorithms concisely.  
 OLS admits an elementary
 theory
 (see Proposition~\ref{prop:ls:reg:properties}),
 that enables
  a general (distribution-free) but tight error analysis 
in Section \ref{section:ML:err an}.

\begin{definition}[Ordinary least-squares regression]\label{def:ls}
For $l,l'\geq 1$ and for probability spaces $(\tilde \Omega,\tilde \cF,\tilde \P)$ and $(\R^l,\cB(\R^l),\nu)$,
let $S$ be a $\tilde\cF\otimes \cB(\R^l)$-measurable
$\R^{l'}$-valued function such that  $S(\omega,\cdot)$ {is in} $ \L_2(\cB(\R^l),\nu)$ for $\tilde \P$-a.e. $\omega \in \tilde \Omega$,
and $\LinSpace{}$ a linear subspace of $\L_2(\cB(\R^l),\nu)$,
 spanned by {some (finite or countable)
 set of} deterministic $\R^{l'}$-valued functions $\{p_k(.)\,:\, k\geq 1\}$.
The least squares approximation of $S$ in the (closure $\bar{\LinSpace{}}$of) space $\LinSpace{}$ with respect to $\nu$ is the ($\tilde \P\times\nu$-a.e.) unique,
$\tilde \cF\otimes\cB(\R^l)$-measurable function
\begin{equation}
\label{eq:mc:mls}
 S^\star(\omega, \cdot) :=\arg\inf_{\phi \in \LinSpace{} }\int |\phi(x) - S(\omega,x) |^2 \nu(dx) = \arg\min _{\phi \in \bar{\LinSpace{}} }\int |\phi(x) - S(\omega,x) |^2 \nu(dx).
\end{equation}
We say that $S^\star$ solves $\OLS(S,\LinSpace{},\nu)$.

{On the other hand, suppose that $\nu_M=M^{-1} \sum_{m=1}^M \delta_{\cX^{(m)}}$ is a discrete probability measure 
 on  $(\R^l,\cB(\R^l))$,
where $\delta_x$ is the Dirac measure on $x$ and $\cX^{(1)},\ldots,\cX^{(M)}: \tilde \Omega \rightarrow \R^l$ are   i.i.d. random variables.
For an $\tilde\cF\otimes \cB(\R^l)$-measurable $\R^{l'}$-valued function $S$ such that $\big|S\big(\omega,\cX^{(m)}(\omega)\big)\big| < \infty$ for any $m$ and $\tilde \P$-a.e. $\omega \in \tilde\Omega$,
the least squares approximation of $S$ in the space $\LinSpace{}$ with respect to $\nu_M$ is the ($\tilde \P$-a.e.) unique,
$\tilde \cF\otimes\cB(\R^l)$--measurable function 
\begin{equation}
\label{eq:mc:ls:empi}
 S^\star(\omega, \cdot) := \arg\inf_{\phi \in \LinSpace{} }
\frac 1M \sum_{m=1}^M |\phi\big(\cX^{(m)}(\omega) \big) - S \big(\omega,\cX^{(m)}(\omega) \big) |^2 .
\end{equation}
We say that $S^\star$ solves $\OLS(S,\LinSpace{},\nu_M)$.
}
\end{definition}

In order to 
{explain} the 
{computational obstacles to be addressed, let us}
first
express the Markov functions given \HXp(ii) in terms of an algorithm involving OLS { by using  Definition \ref{def:ls}}. 
\begin{alg}
\label{alg:mlce}
Initialize by setting $\y 01 \cdot := \Phi(\cdot)$, $\y 00\cdot := \E[\Phi(\X 01)]$ and $T \z 00\cdot := \E[W_T \Phi(\X 01)]$.
Recursively for $k \ge 1$, assume that $(\y {k-1}{}\cdot, \z {k-1}{}\cdot)$ have already been computed, set $\y k{2^k}\cdot = \Phi(\cdot)$, and, 
for any $i \in \{0,\ldots , 2^k-1\}$, let $\LinSpace {l'}$ be the space  $\L_2(\cB(\R^d), \P \circ (\X ki)^{-1} ; (\R^{l'})^\top )$
{for $l'\in\{1,q\}$, and}
\begin{align}
\label{eq:multilevel:cdnexp:ols}
\left.
\begin{array}{l}
\Y ki (\cdot) \text{ solves }  \OLS ( \Obsk k{{Y},i} {\vecx{x}, \vecx{\bar x},\vecx{w}} , \LinSpace 1 , \nu_{k})  \text{ {for}} \\
\\
\qquad \Obsk k{{Y},i} {\vecx x,  \vecx{\bar x},\vecx w}  : =\Phi(x_{2^k})  - \sum_{j =\alpha(i) +1}^{2^{k-1} - 1} \Z {k-1}j (\bar x_{j}) (w_{2j} + w_{2j +1}),
\\ \\
  \Z ki (\cdot )  \text{ solves }  \OLS( \Obsk k{Z,i} {\vecx x,  \vecx{\bar x},\vecx{w}} ,\LinSpace q , \nu_{k})  \text{ {for}} \\
  \\
\displaystyle \qquad \Obsk k{Z,i} {\vecx x,  \vecx{\bar x},\vecx w} : =   { w_i \over \Del k_{i}} \left( \Obsk k{Y,i} {\vecx x,  \vecx{\bar x},\vecx w} - \y ki{x_i} \right),
\end{array}
\right\}
 \end{align}
{for}
$\vecx x = (x_0, \ldots , x_{2^k}) \in \R^{(2^k+1) \times d}$, $\vecx {\bar x} = (\bar x_{0} ,\ldots, \bar x_{2^{k-1}}) \in \R^{(2^{k-1}+1)\times d}$, 
$\vecx{w} = (w_0, \ldots, w_{2^{k}  -1}) \in \R^{2^{k}  \times q}$, and
 $\nu_{k}$ {being} the law of $( \X k0, \ldots, \X k{2^k},  \X {k-1}0, \ldots, \X {k-1}{2^{k-1}}, \DW {k}{0},\ldots,\DW{k}{2^{k}-1})$.
\end{alg}
{\bf Intuition for Definition \ref{def:ls}. }
In Algorithm \ref{alg:mlce} above, we {are using}
Definition \ref{def:ls} 
{with respect to the (theoretical) law instead of the empirical measure (as in Algorithm~\ref{alg:mlmc})}. 
Here,
$l = (2^k+1) \times d+2^{k}\times q$ and $l' = 1$ (resp. $q$).
The function $S(\cdot)$ is given by $\Obsk k{Y,i}\cdot$ (resp. $\Obsk k{Z,i}\cdot$), which is 
deterministic, 
 {hence there is no need for a} probability space $(\tilde \Omega,\tilde \cF ,\tilde \P)$ {here}.
Finally, the measure $\nu$ is the law of the trajectories of the Markov chain $\X k{}$ and the Brownian increments $\DW k{}$, i.e. $\nu = \nu_{k}$.

Algorithm \ref{alg:mlce} {in this form is not really implementable, but } illustrates 
two  computational issues
{ that we are going to} overcome with the empirical least-squares regression algorithm in 
Section \ref{section:ML:algorithm} below: 
firstly, the linear space $\LinSpace 1$ (resp. $\LinSpace q$) is usually infinite dimensional,
{which is infeasible for actual} 
computations;
secondly, 
{generic actual} computation of the integrals \eqref{eq:mc:mls}  
{is hindered in general by the fact that}
 { the  law $\nu_{k}$ may not be available in explicit terms}.

\subsection{Fully implementable algorithm}
\label{section:ML:algorithm}

To avoid 
{regression onto} possibly infinite dimensional space{s} $\LinSpace 1$
 {and} $\LinSpace q$ as in Algorithm \ref{alg:mlce}, we regress on 
 predetermined (user defined) finite dimensional subspace{s},
defined as linear spans of {finite sets of} so-called basis functions:
\begin{definition}[Finite dimensional approximation spaces]
\label{def:fin dim:space}
For each $k \ge 1$ and $i \in \{0,\ldots, 2^k -1\}$, define finite-dimensional functional 
linear spaces
of dimension
$\BasDimk k{Y,i}$ (resp. $\BasDimk k{Z,i}$) by
\begin{equation*}
\begin{cases}
  \LinSpace{Y,k,i} := {\rm span}\{p_{Y,k,i,1},\ldots, p_{Y,k,i,\BasDim{Y,k,i} } \} \text{ for } p_{Y,k,i,j} :\R^d \rightarrow \R \text{ s.t. }
\e[|p_{Y,k,i,j}(\X ki)|^2] < +\infty,\\[2mm]
 \LinSpace{Z,k,i} := {\rm span}\{p_{Z,k,i,1},\ldots, p_{Z,k,i,\BasDim{Z,k,i} } \} \text{ for } p_{Z,k,i,j} :\R^d \rightarrow \R^q \text{ s.t. }
\e[|p_{Z,k,i,j}(\X ki)|^2] < +\infty.
\end{cases}
\end{equation*}
The {\em minimal} error {afforded} by these approximation spaces is denoted
\begin{equation*}
{\Tk{Y,k}{1,i}  := \inf_{ \phi \in \LinSpacek k{Y,k,i} }\e \Big[  |\phi(\X ki) - \y ki{\X ki}|^2 \Big],}
\qquad
{\Tk{Z,k}{1,i}  := \inf_{ \phi \in \LinSpacek k{Z,k,i} }\e \Big[  |\phi(X_i) - \z ki{X_i}|^2 \Big].}
\end{equation*}
\end{definition}

To avoid integration with respect to 
{some (computationally inaccessible)} law $ \nu_i$, as in Algorithm \ref{alg:mlce}, 
{ the next Algorithm~\ref{alg:mlmc} will use simulation to approximate it by} the empirical measure.
\begin{definition}[Simulations and empirical measures]
\label{def:sims and empi}
For $k\ge 0$, generate $M_k\geq 1$ independent copies ({\em simulations})
{$\cC_k := \{( \DW{k,m}{},\X{k,m}{}, \X{k-1,m}{} ) \; : \  m = 1,\dots,M_k \}$}
 of the trajectories of the Markov chains 
and the Brownian increments $(\DW k{} , \X k{}, \X{k-1}{})$.
Denote by $\emeas k$ the empirical probability measure of the $\cC_k$-simulations, 
i.e.
\begin{equation*}
\emeas k = \frac 1{M_k} \sum_{m=1}^{M_k} \delta_{( \X {k,m}0, \ldots, \X {k,m}{2^k},  \X {k-1,m}{0}, \ldots, \X {k-1,m}{2^{k-1}},  \DW {k,m}0, \ldots, \DW {k,m}{2^k - 1}) }.
\end{equation*}
Denote by $\Samp m$ the concatenation of the trajector{ies} of the Markov chains $\X {k,m}{}$, $\X {k-1,m}{}$ and the Brownian increments $\DW {k,m}{}$, 
i.e.

\[\Samp m := (\X{k,m} 0 ,\ldots , \X {k,m}{2^k}, \X{k-1,m}0 ,\ldots , \X {k-1,m}{{2^{k-1}}}, \DW{k,m}0 , \ldots, \DW {k,m}{2^k-1}).\]
\end{definition}
\noindent{\bf Notation and assumptions for the simulations.}
 {Each} $\cC_k$ forms a \emph{cloud of simulations}.
Without loss of generality, up to a generation of extra simulations, we assume 
$M_k \geq \max_{0 \le i \le 2^k-1} \BasDimk k{Y,i} \vee \BasDimk k{Z,i}$.
Furthermore, 
{let} the  clouds of simulations $(\cC_k \ : \ k \ge 0)$ 
{be} independently generated.
All  clouds are defined on {one} probability space $(\Om M, \F M{} , \Prob M)$.
{To construct the probability space that supports the analysis of our algorithm, we simply extend the previous probability space} 
supporting   $(\DW k{}, \X k{})_{k\ge 0}$, 
which serves as a generic element for {any single} simulations,
{by passing to the usual} product space $(\bar \Omega, \bar \cF, \bar \P)=(\Omega, \cF,\P) \otimes (\Om M, \F M{} , \Prob M)$.
{To simplify} notation, we  write $\P$ (resp. $\E$) {instead of} $\bar \P$ (resp. $\bar \E$). 

{In the sequel, we will frequently use} conditioning 
 to integrate {only} with respect to a specific cloud {of simulations}, rather than to take global expectation; 
the  following $\sigma$-algebras will be used for this.
\begin{definition}
\label{def:sim cnd}
For every $k\ge 0$ and  $i\in\{0,\ldots,2^k-1\}$, define the $\sigma$-algebras
\begin{equation*}
\F *k  :=   \sigma(\cC_{k},\dots,\cC_0), \qquad
\F {M}{k,i}  :=  \F *{k-1} \vee \sigma(\X{k,m}{j}, \ \X{k-1,m}{\alpha(j)} \; : \; 1 \le m \le M_i, \ 1 \le j \le i)
\end{equation*} 
and let 
$\CE *_{k}[\cdot] $ (resp. $\CE M_{k,i}[\cdot] $) be the conditional expectation with respect to $\F *k$ (resp. $\F M{k,i}$).
\end{definition}

{Now, we are in position to formulate} a fully implementable algorithm:
\begin{alg}
\label{alg:mlmc} 
Initialize by setting $\y {0,M}1{\cdot} = \Phi(\cdot)$ and 
\[
\y {0,M}0{\cdot} = {1 \over M_0} \sum_{m=0}^{M_0} \Phi(\X {0,m}1 )
\quad \text{and}\quad \z {0,M}0\cdot = {1 \over M_0} \sum_{m=0}^{M_0} \frac{\Phi(\X {0,m}1 ) \DW{0,m}0 }{T}.
\]
Recursion for $k \ge 1$: Assume that $(\y {k-1,M}{}{\cdot}, \z{k-1,M}{}{\cdot} )$ have {already} been computed. Set $\y {k,M}{N}{\cdot} := \Phi(\cdot)$, and, for each $i\in \{0,\dots,2^k-1\}$,  compute first $\y {k,M} i\cdot$ and then $\z {k,M} i\cdot$ as follows: 
\begin{equation}
\label{eq:ml:functions}
\y {k,M} i\cdot := \cT_{C_y}\big(\psimk k_{Y,i}(\cdot) \big)  \quad
\text{and} \quad \z {k,M} i\cdot = \cT_{{\Cz {k,i} }}\big( \psimk k_{Z,i}(\cdot) \big),
\end{equation}
where 
the bounds $C_y$ and $\Cz {k,i}$ are  from Corollary \ref{cor:bd:z},
 the truncation functions $\cT_{C}(\cdot)$ are  {from} Section \ref{section:notation}, and
\begin{align}
\left.
\begin{array}{l}
\displaystyle \psimk k_{Y,i}(\cdot) \text{ solves }  \OLS ( \ObsMk k{Y,i} { \vecx x, \vecx{\bar x}, \vecx w}  , \LinSpacek k{Y,i} , \emeas k)  , 
\\
\displaystyle \ObsMk k{Y,i} {\vecx x, \vecx {\bar x}, \vecx w}  := \Phi(x_{2^k})  - \sum_{j =\alpha(i) +1}^{2^{k-1} - 1} \Z {k-1,M}j (\bar x_{j}) (w_{2j} + w_{2j+1}), \\
\\
\displaystyle \psimk k_{Z,i}(\cdot)  \text{ solves }  \OLS( \ObsMk k{Z,i} {\vecx x, \vecx{\bar x},  \vecx w} , \LinSpacek k{Z,i} , \emeas k)  , \\
\\
\displaystyle \quad \ObsMk k{Z,i} {\vecx x, \vecx {\bar x}, \vecx w} : =    {w_i \over \Del k_{i}} \left( \ObsMk k{Y,i} {\vecx x, \vecx{\bar x}, \vecx w} - \y {k,M}i{x_i} \right),
\end{array}
\right\}
\label{eq:PsiM}  
 \end{align}
for
 $\vecx x = (x_0, \ldots , x_{2^k}) \in \R^{(2^k +1) \times d}$, $\vecx {\bar x} = (\bar x_{0} ,\ldots, \bar x_{2^{k-1}}) \in \R^{(2^{k-1}+1)\times d}$, 
$\vecx{w} = (w_0, \ldots, w_{2^{k}  -1}) \in \R^{2^{k} \times q}$.
\end{alg}

{\bf Intuition for Definition \ref{def:ls}. }
In Algorithm \ref{alg:mlmc} above, we are clearly in the empirical measure setting 
of Definition \ref{def:ls}. 
Here
$l = (2^k+1) d \times (2^{k-1}+1) d \times 2^{k}  q$ and, for each $m \in \{ 1, \ldots, M_k\}$, the $\R^l$-valued random variable $\Samp m$ is the trajectory of the Markov chains and Brownian increments as given in Definition \ref{def:sims and empi}.
$\nu_M$ is the empirical measure $\emeas k$.
The probability space $(\tilde \Omega , \tilde \cF , \tilde \P)$ is $(\Omega, \F M{}  , \P)$, 
i.e. the space generated all the sample clouds $\{\cC_k \  : \ k\ge 0\}$. 
The random function $S(\cdot)$ is the sample dependent function $\ObsMk k{Y,i}\cdot$ (resp. $\ObsMk k{Z,i}\cdot$), which is clearly $\F *{k-1} \otimes \cB(\R^{l})$-measurable.

{Actual} computation of OLS in the  Algorithm \ref{alg:mlmc} {uses}
numerical linear algebra \cite{golu:vanl:96}.

\subsection{Error analysis}
\label{section:ML:err an}
In this section, we determine upper bounds for the global error of Algorithm \ref{alg:mlmc}
\begin{equation}
\label{eq:glob:er}
\bar\cE(k) := \max_{0\le i \le 2^k-1} \bar \cE(Y,k,i) + \sum_{i=0}^{2^k-1} \bar \cE(Z,k,i) \Del k_i
\end{equation}
on each level $k\ge0$, for local error terms given by $\bar \cE(Y,k,i) := \E[|\y k{i}{\X{k}i} - \y{k,M}{i}{\X{k}i}|^2]$ and  $\bar \cE(Z,k,i) := \E[|\z k{i}{\X{k}i} - \z{k,M}{i}{\X{k}i}|^2]$.
In order to do so, it will suffice to find upper bounds for the error terms 
\begin{equation}
\left.
\begin{array}{rcl}
\cE(Y,k,i) & := & \E[ \frac{1}{M_k}\sum_{m=1}^{M_k}|\y k{i}{\X{k,m}i} - \y{k,M}{i}{\X{k,m}i}|^2], \\
\cE(Z,k,i) & := & \E[ \frac{1}{M_k}\sum_{m=1}^{M_k}|\z k{i}{\X{k,m}i} - \z{k,M}{i}{\X{k,m}i}|^2] 
\end{array}
\right\}
\label{eq:empi:error}
\end{equation}
{thanks to} the relationship 
in 
Proposition~\ref{prop:eq:y:z:err:deco:M} 
({similar to} \cite[Prop.4.10]{gobe:turk:13a}):
\begin{proposition}\label{prop:eq:y:z:err:deco:M} 
For each $k\in\{0,\ldots,\kappa\}$ and $i \in \{0,\dots , 2^k-1\}$, we have
\begin{align}
  &\bar \cE(Y,k,i)
 \le 2 \cE(Y,k,i)  
+  {2028 (\BasDimk k{Y,i}+1) C_y^2 \log(3M_k) \over M_k},
\nonumber
\\
&\bar \cE(Z,k,i) \le 2\cE(Z,k,i)  
+  {{2028 (\BasDimk k{Z,i}+1) q \Cz{k,i}^2 \log(3M_k) } \over  M_k};
\nonumber
\end{align}
we recall that $C_y = C_\Phi$ and $\Cz{k,i} = C_X^2/(T-\tk ki)^{(1-\theta)/2}$ from Corollary \ref{cor:bd:z}.
\end{proposition}
Since $\y{k,M}{i}{\cdot}$ and $\z{k,M}{i}{\cdot}$ is computed with the samples $\X {k,m}i$, 
which are also used in the empirical norm inside the expectation of $\cE(Y,k,i)$ and $\cE(Z,k,i)$,
 it turns out that the error analysis of $\cE(Y,k,i)$ and $\cE(Z,k,i)$ is more tractable than that of $\bar\cE(Y,k,i)$ and $\bar\cE(Z,k,i)$, 
and 
{ an important aim for our analysis will be} to find upper bounds for these terms;
 Proposition \ref{prop:eq:y:z:err:deco:M} then allows us to compute upper bounds the $\bar\cE(Y,k,i)$ and $\bar\cE(Z,k,i)$
 from  $\cE(Y,k,i)$ and $\cE(Z,k,i)$ and a correction in terms of the number of basis functions, the number of simulations, the time-grid, and the almost sure bounds $C_y$ and $\Cz{k,i}$. 
It turns out that the correction term is of the same order as one of the error terms in the estimate of $\cE(Z,k,i)$, up to the $\ln(M_k)$ term; see Theorems \ref{thm:0:z:er} and \ref{thm:spec:basis}.
Therefore, the impact of the correction terms on the convergence rate of the global error is essentially the same as the impact of the terms $\cE(Y,k,i)$ and $\cE(Z,k,i)$.
The proof of Proposition \ref{prop:eq:y:z:err:deco:M} is analogous to the proof of \cite[Proposition 4.10]{gobe:turk:13a},
{as} the latter involves only almost sure bounds and general concentration of measure inequalities (\cite[Proposition 4.9]{gobe:turk:13a}).
{Therefore,} we provide no proof here {but refer} 
to that paper.
The correction terms in Proposition \ref{prop:eq:y:z:err:deco:M} have an interpretation as the error due to {\em interdependence} between the cloud used to construct $(\y {k,m}{}{\cdot},\z {k,m}{}{\cdot})$ and the sample used for the empirical norm.

It will be convenient to use the following notation of random norms  in  subsequent analysis; the norms are random because their values depend on the samples of Definition \ref{def:sims and empi} and no global expectation is taken.
\begin{definition}
Let
$\varphi: \Om M \times \R^d \rightarrow \R$ or $\R^q$ be $\F M{} \otimes \cB({\R^d})$-measurable.
For each $k\ge 0$ and $i\in\{0,\ldots,2^k-1\}$, define the random norms
\begin{align}
\lawnorm{\varphi}{k,i}^2  := \int_{\R^d} |\varphi(x)|^2\ \P \circ (\X ki)^{-1} (dx) \quad \text{{and}}\quad
\empinorm{\varphi}{k,i}^2  := \frac{1}{M_i} \sum_{m=1}^{M_i} |\varphi(\X{k,m}i)|^2.
\nonumber
\end{align}
\end{definition}
The norm $\lawnorm \cdot{k,i} $ makes use of the law of $\X ki$, whereas $\empinorm\cdot{k,i}$ makes use of the empirical measure of the samples $\{\X{k,m}i \ : \ m = 1, \ldots , M_k\}$.
Indeed, the error terms \eqref{eq:empi:error} can be written $\cE(Y,k,i) = \E[ \empinorm{\y k{i}{\cdot} - \y{k,M}{i}{\cdot}}{k,i}^2]$ and  $\cE(Z,k,i) = \E[\empinorm{\z k{i}{\cdot} - \z{k,M}{i}{\cdot}}{k,i}^2]$.
Moreover, it follows from the tower law that
\[
\bar \cE(Y,k,i) := \E[\lawnorm{\y k{i}{\cdot} - \y{k,M}{i}{\cdot}}{k,i}^2] \quad \text{{and}}\quad \bar \cE(Z,k,i) = \E[\lawnorm{\z k{i}{\cdot} - \z{k,M}{i}{\cdot}}{k,i}^2].
\]
We come to the main results of this paper, the error propagation of Algorithm \ref{alg:mlmc}.
Two theorems are presented based on different assumptions.
The proofs of the two theorems are very similar in that they are based on a common error decomposition technique.
For this reason, we prove them simultaneously 
{and explain} where the proofs differ; 
the proofs are lengthy and 
deferred to Section \ref{section:main:proof}.

\begin{theorem}
\label{thm:0:z:er} 
In addition to the general assumptions, assume also 
\begin{enumerate}
\item[\HXppp]
{For any time point $t\in\pik {k_{1}} \cap \pik {k_{2}}$ that belongs to}  
two time-grids {$\pik {k_{1}}$ and $\pik {k_{2}}$}
 {for some $k_1,k_2$}, 
it{ holds} that $\X {k_1}i = \X {k_2}j$.
\end{enumerate}
{Then, f}or every $k\ge 0$, $i\in\{0,\ldots,2^k-1\}$, 
the error term $\cE(Y,k,i)$ is bounded above by
\begin{align}
& 
 { 4 \times 2^{-k} \BasDimk k{Y,i} { \delta}\over   M_k} \left\{   3 C_X^2 + (2+q)\right\} 
 + 
 {2 \BasDimk k{Y,i} \over  M_k }   \sum_{j = \alpha(i) +1}^{2^{k-1}-1} \left| \z {k-1}{j}{\cdot} -\z {k-1,M}{j}{ \cdot} \right|_\infty^2 \Del{k-1}_j  
+\Tk {Y,k}{1,i}
\label{eq:MC:y:er}
\end{align}
and the error term $\cE(Z,k,i)$ is bounded by
\begin{align}
\displaystyle
 & { {  12 \delta \BasDimk k{Z,i}  (2 +5T^{1-\theta} )C_X^2 \over c_X M_k (T -\tk ki)^{1-\theta} } 
  + 4  \delta {(2+q) \BasDimk k{Z,i} C_{X} \over c_X M_k} }
 \nonumber\\
& 
{+
{2 \BasDimk k{Z,i} \over \Del k_ i M_k }   \Big\{  \left|\y {k}{i}{\cdot} - \y {k,M}{i}{\cdot} \right|_\infty^2 + \sum\limits_{j = \alpha(i) +1}^{2^{k-1}-1} \left| \z {k-1}{j}{\cdot} -\z {k-1,M}{j}{ \cdot} \right|_\infty^2 \Del{k-1}_j  \Big\} 
  +\Tk {Z,k}{1,i}
}
\label{eq:MC:z:er}
\end{align}
\end{theorem}

\begin{remark*}
The Assumption \HXppp \ is trivially valid
if $\X k{}$ 
{can be taken as the finite dimensional marginals of $X$, provided those are available in closed form, like for instance for (geometric) Brownian motion.}
If one 
{is computing} $\X k{}$ with an Euler scheme, for example, one would fix a maximal level, $\kappa$ say, 
and 
{could} obtain $\X k{}$ by running the Euler scheme {once on the finest} time-grid $\pik \kappa$ 
{and selecting} only the values associated $\pik k$ {for every $k \le \kappa$}.
We  remark that the assumption is not necessary for Theorem \ref{thm:spec:basis}.

\end{remark*}

\begin{remark}
\label{rem:best/worst:gen}
The error bounds in Theorem \ref{thm:0:z:er} above are
{not easy to apply}  
because 
{it appears difficult to} quantify the terms in the norms $|\cdot|_\infty$ {more explicitly};
these norms are stronger than the norms used to quantify the error $\cE(\cdot,k,i)$, and we have no precise estimates for them.
It seems to be difficult to 
{replace} the use of this strong norm in general, cf.\ \cite{benz:gobe:13} who 
{obtain} estimates using the same norm when using general basis functions.
However, we can plug the absolute $y$- and $z$-bounds
 Algorithm \ref{alg:mlmc} into  (\ref{eq:MC:y:er},\ref{eq:MC:z:er}) to obtain a rough upper bound on the error
\begin{align*}
& \cE(Y,k,i)  \le \Tk {Y,k}{1,i} + {  16 \times  2^{-k  } \BasDimk k{Y,i} \over M_k}
+{8 (C_X^2 \vee C_\Phi^2 ) T^\theta \BasDimk k{Y,i} \over  \theta M_k },
\\
&  \cE(Z,k,i)  \le    \Tk {Z,k}{1,i} + { 16(2+T^{1-\theta}) \BasDimk k{Z,i} \over c_X (T-\tk ki)^{1-\theta}  M_k } 
 +
{8  (C_{X}^2\vee C_\Phi^2) (1 + T^\theta \theta^{-1}) \BasDimk k{Z,i} \over   \Del k_ i M_k } .
\end{align*}
This is the ``worst-case" error estimate in the sense that we assume the terms $ \left| \y {k-1}{j}{\cdot} -\y {k-1,M}{j}{ \cdot} \right|_\infty$ and $ \left| \z {k-1}{j}{\cdot} -\z {k-1,M}{j}{ \cdot} \right|_\infty$ are maximal.

We use this estimate for a crude comparison to the usual least-squares multistep forward dynamical programming 
(LSMDP) scheme {\cite{gobe:turk:13a}}, i.e. Algorithm \ref{alg:mlmc} with the correction terms from level $k-1$ removed.
For every $k \ge 0$ and $i \in\{0,\ldots ,2^k -1\}$, the {corresponding} error estimates for the LSMDP algorithm (when the same time-grid and the same basis functions are used) are
\begin{equation}
\cE_{\text{MDP}}(Y,k,i) \le \Tk {Y,k}{1,i} + {C_{\Phi}^2 \BasDimk k{Y,i} \over M_{k} } \quad \text{and}\quad
\cE_{\text{MDP}}(Z,k,i) \le \Tk {Z,k}{1,i} + {C_{\Phi}^2 \BasDimk k{Z,i} \over \Del k_{i} M_{k} }.
\label{eq:er:MDP}
\end{equation}
We see that the dependence on the number of basis functions $\BasDimk k{\cdot,i}$, the time increment $\Del k_i$, and the number of simulations $M_k$ is the same for both the multilevel and the MDP scheme, although the constants may differ.
In this setting, the behavior of the both algorithms with respect to each of these parameters 
{might} be the same.
We {emphasize, however,} that this is  {a rather rough} ``worst case scenario", in which the approximations of  $\zn {k-1}{}$ are as bad as 
{absolute a-priori bounds would permit}.
We now turn 
to the opposite extreme, a ``best case" scenario, in which the $|\cdot|_\infty$-terms  are negligible;
studying (\ref{eq:MC:y:er},\ref{eq:MC:z:er}), we see these terms are counter-balancing the negative impact of the time-increment $\Del k_i$. 
This leads to ``best-case" error estimates
\begin{align*}
\cE_{\text{best}}(Y,k,i) & \le \Tk {Y,k}{1,i} + {  16 \times  2^{-k  } \BasDimk k{Y,i} \over   M_k}
 \quad \text{and}\quad  \cE_{\text{best}}(Z,k,i)  \le    \Tk {Z,k}{1,i} + { 2 \BasDimk k{Z,i} \over c_X   M_k } ,
\end{align*}
{motivating}, in particular,
an improvement in the dependence on 
$\Del k_i$, which no longer appears in the denominator.
\end{remark}

Remark \ref{rem:best/worst:gen}{, despite its  crude quantitative nature, is}  encouraging { as a  first} comparison between the multilevel  algorithm and the LSMDP algorithm.
 {A main obstacle for}  more precise statements was that 
{error bounds were given in terms of} very strong norms in the setting  {with a}  general basis.
We 
{next provide} a more precise comparison to the LSMDP 
{for} a specific choice of basis.

\begin{theorem}
\label{thm:spec:basis} 
In addition to the general assumptions, assume that
\begin{enumerate}
\item[\HK] 
{
the basis functions are  indicator functions, i.e. $p_{\cdot,k,i,j} := \1_{A_{\cdot,k,i,j} }$ on {disjoint} sets $A_{\cdot,k,i,j}$;
moreover, there exists $\delta \ge 1$ such that either $\P(\X ki \in A_{\cdot,k,i,j}) \ge  1 / (\delta K)$. 
}
\end{enumerate}
For every $k\in\{1,\ldots,\kappa\}$, $i\in\{0,\ldots,2^k-1\}$, 
the error term $\cE(Y,k,i)$ is bounded above by
\begin{align}
& { 4 \times 2^{-k} \BasDimk k{Y,i} { \delta}\over   M_k} \left\{   3 C_X^2 + (2+q)\right\} 
 + 
 {2\BasDimk k{Y,i} {\delta} \over  M_k }   \sum_{j = \alpha(i) +1}^{2^{k-1}-1} \bar \cE(Z,k-1,j) \Del{k-1}_j  
+\Tk {Y,k}{1,i}
\label{eq:spec:y:er}
\end{align}
and the error term $\cE(Z,k,i)$ is bounded by
\begin{align}
& 
 {  12 \delta \BasDimk k{Z,i}  (2 +5T^{1-\theta} )C_X^2 \over c_X M_k (T -\tk ki)^{1-\theta} } 
  + 4  \delta {(2+q) \BasDimk k{Z,i} C_{X} \over c_X M_k} 
  +\Tk {Z,k}{1,i} \nonumber \\
 & +
{2 \BasDimk k{Z,i} { \delta} \over \Del k_ i M_k }  \left\{  {8qC_X^2  T^\theta \theta^{-1} \Del k_i \over c_X}  
+ \bar \cE(Y,k,i) +   2kq \ln(2)  \sum_{j = \alpha(i) +1}^{2^{k-1}-1} \bar \cE(Z,k-1,j) \Del{k-1}_j  \right\}   .
\label{eq:spec:z:er}
\end{align}
\end{theorem}
The assumption \HK \  would be satisfied, e.g., if the sets $A_{\cdot,k,i,j}$ 
{have equal probability $1/K$ (thus $\delta=1$)} under the 
measure $\p \circ (\X ki)^{-1}$.
Moreover, if $\X ki$ has a density $\phi_X(x)$ that is bounded from below {away from zero} on a compact $A \subset \R^d$, it follows that $\P(\X ki \in H_j) \ge \min_{x\in A} \phi_X(x) \int_{H_j} {1} dx$ for all sets $H_j \subset A$, 
so 
{a partition of $A$ into sets} $A_{\cdot,k,i,j}$,  satisfying $\int_{A_{\cdot,k,i,j}} {1} dx = const/ \BasDimk k{\cdot,i}$ for all $j =1,\ldots ,\BasDimk k{\cdot,i}$ would form a basis  satisfying \HK.

\begin{remark*}
{Let us note that}  Theorem \ref{thm:0:z:er} {differs from} 
Theorem \ref{thm:spec:basis}
 {in}  that {the latter does not require}
assumption \HXppp  \
{and uses a weaker norm than $|\cdot|_{\infty}$}.
We  believe that estimates in the weaker norm and  {moreover relaxing} the assumption \HXppp  \ may 
{hold true} for a class of basis functions 
 {beyond} \HK. 
 {Computational} examples later indeed
indicate that multilevel benefits prevail beyond the assumptions under which complexity gains are proven subsequently.
\end{remark*}

{Comparing the error bounds from Theorem \ref{thm:0:z:er} and Theorem \ref{thm:spec:basis}, one sees that  the terms given in the norm $|\cdot|_\infty$ in the former theorem have  been replaced by the equivalent terms in the weaker norm $\bar \cE(\cdot)$ in the latter.
From theorem \ref{thm:spec:basis}, the upper bound for the global error \eqref{eq:glob:er} { of the multilevel scheme} on the level $k$ is expressed in terms of the number of time-steps, the time-increments, the number of basis functions, the bias of the basis functions, the number of simulations, and the global error \eqref{eq:glob:er}
 on {the previous} level $k-1$ (i.e. the terms $\bar \cE(Z,k-1,j)$).
}
For the remainder of this section, the order notation $O(\cdot)$ will be used: {we write} $g(y)$ {is in} $O(y)$
 if there exists a constant $C$,
 {not depending on} 
 the level $k$, such that 
 {$\limsup_{y\downarrow 0} g(y)/y \le  C$}.
We set the numerical parameters -- the basis functions and the number of simulations -- of the multilevel algorithm 
so that the global error $\bar\cE(k)$ attains a precision level $O(\varepsilon)$ for 
$\varepsilon >0$.
We use this calibration to compute complexity and to compare the multilevel algorithm to the LSMDP scheme more precisely than in Remark \ref{rem:best/worst:gen}.

\begin{remark}
\label{rem:perid}
Our theoretical complexity analysis below applies error estimates  from Theorem~\ref{thm:spec:basis} and hence requires its assumption to hold; in particular, the basis is to satisfy \HK . 
Furthermore, it is required that the basis is such that the approximations errors $\Tk {Y,k}{1,i}$ and $\Tk {Z,k}{1,i}$ are of order $O(\varepsilon)$ for basis dimension $\BasDimk k{Y,i}$ resp. $\BasDimk k{Z,i}$ as stated in {\bf Choice of basis} below.
In combination, these assumptions appear restrictive, but computational examples later will indicate empirical multilevel benefits beyond these assumptions.    
{On the other hand, there is a class of examples in which the required assumptions are satisfied, and we exhibit this class for the remainder of this remark.}
Suppose that the solutions $x\mapsto (\y kix, \z kix)$ are periodic, that is, there exists 
$\lambda:=(\lambda_1,\dots,\lambda_d)\in \R^d_+$ such that
$(\y kix, \z kix) =  (\y ki{x + n \lambda } ,\z ki{x + n\lambda } )$ for all $n \in \ZZ$ and $x\in\R^d$.
As an example one can think of $\X ki = W_{\tk ki}$ and $\Phi(x) = \sin(\beta \cdot x)$, whence $\lambda_i = 2\pi/\beta_i$.
More generally, one can consider any periodic terminal condition $\Phi$, and $\X ki$ to be the marginals of the solution to a stochastic differential equation whose coefficient functions have the same periodicity as $\Phi$.
For every $t_0\in (0,T]$ and $\kappa \le k$, we assume that the marginals $ \X \kappa j $ have a density bounded from below by $c(t_0)>0$  (independent of $\kappa$)  in the domain $D :=\otimes_{i=1}^d [-\frac{\lambda_i}2, \frac{\lambda_i}2]$ for all $j$ such that  $\tk \kappa j > t$.
This property is satisfied if $\X ki$ were the marginal of  the solution to a stochastic differential equation whose generator is uniformly elliptic; 
hence the marginal density is bounded from below by a Gaussian density \cite{koha:03b}.
Let $\{ B_{k,i,1}, \ldots , B_{k,i,\BasDimk k{i} } \}$ be a  hypercube partition of $D$, and define the basis functions $p_{\eta, k,i,j}(x)$ ($j =1,\ldots,\BasDimk ki$, $\eta = Y,Z$) to be the indicator functions $\1_{A_{k,i,j}}(x)$ on the sets 
$\displaystyle A_{k,i,j} := \bigcup_{n\in\ZZ} \{x + n\lambda \ : \ x \in B_{i,j} \}$.
Then,
$ \displaystyle
\P(A_{k,i,j}) \ge \P(B_{k,i,j}) \ge {c(t_0) \mu(D) / \BasDimk ki},
$
where $\mu$ is the Lebesgue measure;
hence, the condition \HK \ is satisfied with $\delta = 1/(c(t_0) \mu(D))$.
$\delta$ may therefore be considered a constant with respect to the level $k$ and the precision level so long as one considers the global error on the interval $[t_0,T]$.
\end{remark}

{
\begin{remark*}
{As an alternative case to the periodic one outlined in Remark \ref{rem:perid},
one may also think about a forward process $X$ which is a diffusion reflected within some compact domain, such as to ensure its density being bounded away from zero, so that one could argue similarly as in the periodic case above. 
To make this ansatz rigorous, however, would require $\L_2$-regularity properties like \Ht(ii) to hold for reflected diffusions $X$.  We are not aware of such results being yet available.
}
\end{remark*}
}

The  error bounds  (\ref{eq:spec:y:er} -- \ref{eq:spec:z:er})  of Theorem \ref{thm:spec:basis} 
show that a sufficient criterion to achieve an global error $O(\varepsilon)$ is to ensure that each of the terms in the sums in (\ref{eq:spec:y:er} -- \ref{eq:spec:z:er}) 
is bounded 
by $O(\varepsilon)$, {and} we use this {criterion} to develop a 
calibration procedure.
Furthermore, assume that the assumptions of Remark \ref{rem:perid} hold, namely that the basis functions satisfy \HK \ and periodicity.

\noindent
{\bf $\blacktriangleright$ Choice of basis.}
{We first choose a basis satisfying \HK \ so that $\Tk {Y,k}{1,i}$ and $\Tk {Z,k}{1,i}$ are bounded above by $O(\varepsilon)$ for all $i$.
{
Let $\LinSpacek k{Y,i} = \LinSpacek k{Z,i}$ and the set $\{B_{k,i,j} \ : \ j=1,\ldots , \BasDimk k{\cdot, i} \}$ be the uniform hypercubes on the set $D$.
}
Thanks to \Ht(iv), \HXp, it is  sufficient for the boundedness of  $\Tk{Y,k}{1,i}$, $\Tk{Z,k}{1,i}$
to ensure that
\begin{align*}
\max_{0\le i \le 2^k-1} {\min}_{\phi \in \LinSpacek k{Y,i}} & \E[ |u(\tk ki, \X ki) - \phi(\X ki) |^2  ]  
 + \sum_{i = 0}^{2^k -1} {\min}_{\phi \in \LinSpacek k{Z,i}}\E[ |v(\tk ki, \X ki) - \phi(\X ki) |^2 ] \le O(\varepsilon).
\end{align*}
We assume (for simplicity) that $\theta =1$.
Thanks to  \HX(iv), the Lipschitz constant of $v(t,\cdot)$ is equal to  $O((T-t)^{-1/2})$,
 so it suffices to set the hypercube diameter at time $\tk ki$ equal to $\sqrt{T-\tk ki} O(\sqrt{\varepsilon})$,
{whence the dimension of the basis is $\BasDimk k{Z,i} = (T-\tk ki)^{-d/2}O(\varepsilon^{-d/2})$.}
Since the Lipschitz constant of $u(\tk ki, \cdot)$ is $O(1)$,  it follows that $\Tk{Y,k}{1,i} \le O(\varepsilon)$ with the same basis.
}

{
\noindent
{\bf $\blacktriangleright$ Number of simulations.}
The choice of basis fixes the number of basis functions $\BasDimk k{\cdot,i} = \BasDimk k{\cdot,i} (\varepsilon)$. 
We  choose $M_k = \max_i O(\varepsilon^{-1}  k \BasDimk k{Z,i}(\varepsilon) ) = { O(k  \varepsilon^{-1-d/2} \max_i (T -  \tk ki)^{-d/2}  ) \le } O(k 2^{ kd/2} \varepsilon^{-1-d/2} )$
 to ensure that the all terms in  (\ref{eq:spec:y:er} -- \ref{eq:spec:z:er}) - except those depending on $\bar \cE(Z,k-1,\cdot)$ - 
and also in the correction terms in Proposition \ref{prop:eq:y:z:err:deco:M} are bounded  by $O(\varepsilon)$; observe that we do not need to worry about the terms $1/(T-\tk ki)^{1-\theta}$, because  the sum $\sum_{i=0}^{2^k-1}  \Del k_i/(T-\tk ki)^{1-\theta}$ is  bounded uniformly in $k$.
}

{
\noindent
{\bf $\blacktriangleright$ Iteration to levels $j < k$.}
In the calculations above, it only remains to set parameters such that 
\[\sum_{l = \alpha(i) +1}^{2^{k-1}-1} \bar \cE(Z,k-1,l) \Del{k-1}_j \le O(\Del k_i) \]
for all $i$. 
This is satisfied by setting the precision for the global error $\bar\cE(k-1)$ on level $k-1$ to be less than or equal to $O(\min_i \Del k_i) \le O(2^{-k})$ in the place of $O(\varepsilon)$.
Subsequently, on every level $j \le k-1$ thereafter, we set the precision for the global error  $\bar\cE(j)$  less than or equal to $O(\min_i \Del {j+1}_i) \le O(2^{-j})$ and repeat the first two steps of the procedure above.
For simplicity,  we choose the same basis for every level, although this is  possibly not optimal.
The basis dimension at time $\tk ji$ is  $K(j,i,\varepsilon) := (T-\tk ji)^{-d/2}O(\varepsilon^{-d/2})$
and the number of simulations on level $j<k$ is $M_j = \max_i O(j 2^{jd/2} K(j,2^j-1,\varepsilon) ) = O(j 2^{j + jd/2} \varepsilon^{-d/2}  )$.
}

\noindent
{
{\bf $\blacktriangleright$ Complexity analysis.}
We fix $\varepsilon = O(2^{-k})$,  as this is usually the discretization error between $(y_k ,z_k)$ and the continuous time solution (see \Ht(iv)).
There are two contributions to the computational cost: the cost of simulation of the Markov chain $\X k{}$ and Brownian increments $\DW k{}$, and the cost of the regressions.
The cost of computing the simulations is $O(2^j M_j)$ on level $j$ of the algorithm, therefore the overall simulation cost $\sum_{j=0}^kO(2^j M_j)$.
To compute the cost of the regression, one must first of all remark that there is a closed form formula for regression on indicators (see the partitioning estimate in \cite{gyor:kohl:krzy:walk:02}): for responses $(\psi_m)_{1\le m\le M}$ corresponding to observations $(\phi_m)_{1\le m\le M}$, the precise coefficient of the indicator function denoted by $H$ is given by
}
{
\[
\alpha_H = { \sum_{m=1}^M \psi_m \1_H(\phi_m) \over \sum_{m=1}^M  \1_H(\phi_m) };
\]
therefore, the cost of the regression on each time point is proportional to the cost of sorting the simulations into the indicators, which is proportional to the dimension $d$ times the number of simulations. 
This implies that the cost of the regressions on  level $l$ is equal also  equal to $O(2^l M_l)$.
Therefore, recalling that $\varepsilon = O(2^{-k})$, the overall cost of the algorithm is 
\[
\sum_{j=0}^k O(2^j M_j) \le O(k) \sum_{j=0}^k O(2^{j(1+d)})  = O(\ln(\varepsilon^{-1} +1)\varepsilon^{-2 - d} ).
\]

}

For comparison, we calibrate the basis functions and number of simulations for the LSMDP algorithm described in Remark \ref{rem:best/worst:gen}, using \eqref{eq:er:MDP} 
{in the place of (\ref{eq:spec:y:er} -- \ref{eq:spec:z:er}).}
We choose the same basis functions, and $M_k = O(\varepsilon^{-1} 2^{k} K(k,2^k-1,\varepsilon))$.
Then, setting $\varepsilon = O(2^{-k})$, the overall complexity is 
$
2^k \times M_k  = O( \varepsilon^{-3 - d}  ).
$
We observe that, in comparison to the complexity of the multilevel scheme, one factor in $\ln(\varepsilon^{-1}+1)$ have been replaced by a factor $\varepsilon^{-1}$, which is much larger.
This implies that, in comparison to MDP, the multilevel scheme has a possible efficiency gain of factor $\varepsilon$ (ignoring the log terms). In our setting, is equal to the number of time steps, which is substantial.

\subsection{Proof of Theorems \ref{thm:0:z:er} and \ref{thm:spec:basis}}
\label{section:main:proof}

We state the elementary properties of OLS (Definition \ref{def:ls})  in Proposition \ref{prop:ls:reg:properties} below.
This proposition is in fact the same as \cite[Proposition 4.12]{gobe:turk:13a}, and we refer the reader interested in the proof to that paper.
We are aware that parts (iii) and (iv) of this proposition are given in  high generality, so we  provide some explicit $\sigma$-algebras and functions for the benefit of the reader's intuition following the proposition statement.

\begin{proposition}\label{prop:ls:reg:properties}
With the notation of Definition \ref{def:ls},
suppose that $\LinSpace{}$ is 
finite dimensional and spanned by the functions $\{p_1(.),\dots,p_K(.)\}$.
Let $S^\star$ solve $\OLS(S,\LinSpace{},\nu)$ (resp. $\OLS(S,\LinSpace{},\nu_M)$), according to \eqref{eq:mc:mls} (resp. \eqref{eq:mc:ls:empi}).
The following properties are satisfied:
\begin{enumerate}
\item [\rm (i)] linearity: the mapping $S\mapsto S^\star$ is linear.
\item [\rm(ii)] contraction property: $\|S^\star \|_{\L_2(\cB(\R^l)
,\mu)}\leq \|S\|_{\L_2(\cB(\R^l)
,\mu)}$, where $\mu = \nu$ (resp. $\mu = \nu_M$).
\item [\rm(iii)] conditional expectation solution: 
in the case of the discrete probability
 measure $\nu_M$,
assume additionally that the {sub-}$\sigma$-algebra $\cQ \subset \tilde \cF$ is such that $\big(p_j(\cX^{(1)}),\ldots,p_j(\cX^{(M)}) \big)$ is
$\cQ$-measurable for every $j\in\{1,\ldots,K
\}$. 
Let $S_\cQ (\cdot)$ be any $\tilde \cF \otimes \cB(\R^l)$-measurable, $\R^{l'}$-valued function such that
such that $S_\cQ(\cX^{(m)}):=\tilde\E[S(\cX^{(m)})|\cQ]$ for each $m \in \{1,\dots,M\}$ $\tilde \P$-almost surely.
Then $\tilde \E[S^\star |\cQ](\omega,x)$ solves $\OLS\big( S_\cQ , \LinSpace{}, \nu_M \big)$.
\item  [\rm(iv)]  bounded conditional variance: in the case of the discrete probability measure $\nu_M$, suppose that $S(\omega,x)$ is $\cG \otimes \cB(\R^l)$-measurable, for $\cG \subset \tilde \cF$ independent of $\sigma(\Samp{1:M})$, 
there exists a Borel measurable function $h : \R^l \rightarrow \cE$, for some Euclidean space $\cE$, such that the random variables $\{p_j(\Samp m)  \ : \ m =1,\dots ,M, \ j = 1,\dots, K\}$ are $\cH := \sigma(h(\Samp m) \ : \ m =1,\dots, M)$-measurable, and  there is a finite 
constant $\sigma^2 \geq 0$ that uniformly bounds the conditional variances $\tilde \E\big[|S(\Samp{m})-\tilde \E(S(\Samp{m})| \cG\vee \cH)|^2\ |\ \cG\vee \cH\big]\leq \sigma^2$ $\tilde \P$-a.s. and for all $m \in \{1,\dots,M\}$. Then 
$$\tilde \E\Big[ \| S^\star (\cdot) - \tilde\E[S^\star(\cdot) | \cG\vee \cH] \|_{\L_2(\cB(\R^l),\nu_M)}^2\ \big| \ \cG\vee \cH\Big]  \le \sigma^2 K /M.$$
\end{enumerate}
\end{proposition}

{\bf Intuition for Proposition \ref{prop:ls:reg:properties}.}
The observation $\Samp m{}$ and response $S$ above will be $\Samp{k,m}{} $
 and $\Obsk k{i} {\vecx{x},\vecx{w}} $, respectively, whereas the linear space $\cK$ will be $\LinSpacek k{i}$ and the measure $\nu$ (respectively $\nu_M$) will be $\nu_k$ (respectively $\emeas k$).
For part (iii), we will take $\cQ$ to be  the $\sigma$-algebras $\F M{k,i} $ in Definition \ref{def:sim cnd}; the function  $\tilde\E_\cQ[S(\Samp m)](\cdot)$ will then be equal to $\yn ki$ (respectively $\zn ki$), see below.
For part (iv), we take $\cE =\R^d$ and the Borel function $h : \R^l \to \cE$ to be $h(\Samp m) = \X {k,m} i$, whence the $\sigma$-algebra $\cH$ is $\sigma(\X {k,m}i \ : \ m =1,\ldots,M_k)$. 
We take $\F *{k-1}$ for $\cG$, whence $\cG \vee\cH = \F M{k,i}$.

We now begin the proof of the two theorems.
Recall the $\sigma$-algebras from Definition \ref{def:sim cnd} and the soft truncation function $\cT_r(\cdot)$ in Section \ref{section:notation}.
The Lipschitz continuity (for all $r $) of the function $\cT_{r}(\cdot)$ implies that
\begin{align}
\E[ \normM ki {\yn k{i} - \yn{k,M}{i}} ^2] &  = \E[ \normM ki {\cT_{C_y}(\yn k{i}) - \cT_{C_y}(\psimk k_{Y,i}) } ^2] \le \E[ \normM ki { \yn k{i} - \psimk k_{Y,i} } ^2]
\label{eq:main:y:1}
\\
\E[ \normM ki {\zn k{i} - \zn{k,M}{i}} ^2] & = \E[ \normM ki {\cT_{\Cz{k,i}} (\zn k{i}) - \cT_{\Cz{k,i}}(\psimk k_{Z,i}) } ^2] \le \E[ \normM ki { \zn k{i} - \psimk k_{Z,i} } ^2]
\label{eq:main:z:1}
\end{align}
We introduce the ``fictitious" functions $\psik k_{Y,i} : \R^d \to \R$ and $\psik k_{Z,i} : \R^d \to (\R^q)^\top$ defined by
\begin{align*}
\displaystyle \psik k_{Y,i}(\cdot) \text{ solves }  \OLS ( \Obsk k{Y,i} {\vecx{x},\vecx{w}} , \LinSpacek k{Y,i} , \emeas k)  , \qquad
 \psik k_{Z,i}(\cdot) \text{ solves }  \OLS ( \Obsk k{Z,i} {\vecx{x},\vecx{w}} , \LinSpacek k{Z,i} , \emeas k)  , 
\end{align*}
for functions $ \Obsk k{Y,i}\cdot$ and $ \Obsk k{Z,i}\cdot$ given in \eqref{eq:multilevel:cdnexp:ols} from Algorithm \ref{alg:mlce}; 
the fictitious nature of $\psik k_{Y,i} (\cdot)$ and $\psik k_{Z,i} (\cdot)$  comes from the functions $ \Obsk k{Y,i}\cdot$ and $ \Obsk k{Z,i}\cdot$, which are constructed using the unknown functions $\y k{}\cdot$ and $\z k{} \cdot$, so cannot be computed explicitly.
We will decompose \eqref{eq:main:y:1} and \eqref{eq:main:z:1} using the (random) functions $\CE M_{k,i}[\psik k_{Y,i}](\cdot)$ and $\CE M_{k,i}[\psik k_{Z,i}](\cdot)$, respectively,
but first we make use of Proposition \ref{prop:ls:reg:properties}(iii) to determine that $\CE{M}_{k,i}[\psik k_{Y,i}](\cdot)$ and $\CE{M}_{k,i}[\psik k_{Z,i}](\cdot)$ solve OLS's.
Set $\cQ$ to be the $\sigma$-algebra $\F M{k,i}$.
 $p(\X {k,m}i)$ is $\cQ$-measurable for any $p \in \LinSpace{Y,k,i} \cup \LinSpace{Z,k,i}$. 
Now, since $\z {k-1}j {\X{k-1,m}{j}}$ is $\F M{k,2j}$-measurable for all $j>\alpha(i)$, applying the tower property and the Markov property \HXp \ yields that 
\begin{align}
&\CE{M}_{k,i}[\Obsk k{Y,i} {\cX_m} ]  =  \CE{M}_{k,i}[ \Phi (\X {k,m}{2^k}) ] = \y ki{\X {k,m}i},
\label{eq:cdn:exp:obs:y} \\
&\CE{M}_{k,i}[\Obsk k{Z,i} {\cX_m}]  =  \CE{M}_{k,i}[ {\Phi (\X {k,m}{2^k}) \DW {k,m}i \over \Del k_i }] = \z ki{\X {k,m}i}.
\label{eq:cdn:exp:obs} 
\end{align}
for all $m \in \{ 1, \ldots , M_k\}$,
whence we finally obtain the expression
\begin{align*}
x\in \R^d \mapsto \CE{M}_{k,i}[\psik k_{Y,i}](x) \quad \text{solves} \quad \OLS (  \y ki{x_i} , \LinSpacek k{Y,i} , \emeas k \big), \\
x\in\R^d \mapsto \CE{M}_{k,i}[\psik k_{Z,i}](\cdot) \quad \text{solves} \quad \OLS (  \z ki{x_i} , \LinSpacek k{Z,i} , \emeas k \big).
\end{align*}
Therefore, introducing the random functions $\CE{M}_{k,i}[\psik k_{Y,i}](\cdot)$ and $\CE{M}_{k,i}[\psik k_{Y,i}](\cdot)$ on the right hand side of \eqref{eq:main:y:1} and \eqref{eq:main:z:1}, respectively, and applying Pythagoras' theorem, it follows that
\begin{align}
\E[ \normM ki {\yn k{i} - \yn{k,M}{i}} ^2] 
&\le \E[ \normM ki { \yn k{i} - \CE{M}_{k,i}[\psik k_{Y,i}]}^2] + \E[ \normM ki { (\CE{M}_{k,i}[\psik k_{Y,i}] -\psimk k_{Y,i}) }^2],
\label{eq:main:y:2}
\\
\E[ \normM ki {\zn k{i} - \zn{k,M}{i}} ^2] 
& \le 
\E[ \normM ki { \zn k{i} - \CE{M}_{k,i}[\psik k_{Z,i}]}^2] + \E[ \normM ki { (\CE{M}_{k,i}[\psik k_{Z,i}] -\psimk k_{Z,i}) }^2].
\label{eq:main:z:2}
\end{align}
Moreover, $\E[\normM ki { \zn k{i} - \CE{M}_{k,i}[\psik k_{Z,i}]}^2] \le \Tk k{Z,i}$ and  $\E[\normM ki { \yn k{i} - \CE{M}_{k,i}[\psik k_{Y,i}]}^2] \le \Tk k{Y,i}$, 
and injecting this into inequalities \eqref{eq:main:y:2} and \eqref{eq:main:z:2} yields
\begin{align}
\E[ \normM ki {\yn k{i} - \yn{k,M}{i}} ^2]
\le 
\Tk k{Y,i} + \E[ \normM ki { \CE{M}_{k,i}[\psik k_{Y,i}] -\psimk k_{Y,i} }^2],
\label{eq:MC:Y:1}\\
\E[ \normM ki {\zn k{i} - \zn{k,M}{i}} ^2]
\le 
\Tk k{Z,i} + \E[ \normM ki { \CE{M}_{k,i}[\psik k_{Z,i}] -\psimk k_{Z,i} }^2].
\label{eq:MC:Z:1} 
\end{align}
To treat the second term on the right-hand side of   \eqref{eq:MC:Y:1} (resp. \eqref{eq:MC:Z:1}), we decompose
\begin{align}
\E[ \normM ki { \CE{M}_{k,i}[\psik k_{Y,i}] -\psimk k_{Y,i} }^2] \le 2{\E[ \normM ki {\psik k_{Y,i} -\psimk k_{Y,i} }^2]} + 2{\E[ \normM ki { \CE{M}_{k,i}[\psik k_{Y,i}] -\psik k_{Y,i} }^2]},
\label{eq:main:y:3} 
\\
{\E[ \normM ki { \CE{M}_{k,i}[\psik k_{Z,i}] -\psimk k_{Z,i} }^2]} \le 2{\E[ \normM ki {\psik k_{Z,i} -\psimk k_{Z,i} }^2]} + 2{\E[ \normM ki { \CE{M}_{k,i}[\psik k_{Z,i}] -\psik k_{Z,i} }^2]}.
\label{eq:main:z:3}
\end{align}
We first treat the terms $\E[ \normM ki {\psik k_{Z,i} -\psimk k_{Z,i} }^2]$ and $\E[ \normM ki {\psik k_{Z,i} -\psimk k_{Z,i} }^2]$;
the approach for both terms is identical, so we focus on the upper bound for the latter and only state the result for the former. 
We adopt 
an approach similar to the proof of  Proposition \ref{prop:ls:reg:properties} (iv); see \cite[Appendix A]{gobe:turk:13a} to compare.
First, observe using Proposition \ref{prop:ls:reg:properties} (i) that 
\[
\big[ \psik k_{Z,i} -\psimk k_{Z,i}\big](\cdot) \text{ solves } \OLS ( \Obsk k{Z,i} {\vecx{x},\vecx{w}} -  \ObsMk k{Z,i} {\vecx{x},\vecx{w}} , \LinSpacek k{Z,i} , \emeas k).
\]
Then, since $\LinSpace{Z,k,i}$ is finite dimensional, it has an orthonormal (with respect to the 
norm $ \normM ki {\cdot }^2$) basis $\{\tilde p_1,\ldots, \tilde p_{\tilde K}\}$ with $\tilde K\leq \BasDim{Z,k,i}$.
Using the orthogonality property of $\tilde p$, setting $\alpha^\star := \int \tilde p(x)^\top\{\Obsk k{Z,i}x - \ObsMk k{Z,i}x \}  d \emeas k$, and expanding $|\alpha^\star|^2$ as a summation over the samples yields
\begin{align*}
& \normM ki {\psik k_{Z,i} -\psimk k_{Z,i} }^2=
| \alpha^\star|^2 \nonumber \\
&= \frac{ 1 }{M_{k}^2 }
\sum_{m_1,m_2=1}^{M_{k}}
\trace\Big(\tilde p(\Samp{ m_1}) \tilde p^\top (\Samp {m_2}) \nonumber \\
&\hspace{2cm}(\Obsk k{Z,i}{\Samp {m_{1}} } - \ObsMk k{Z,i}{\Samp {m_{1}}})(\Obsk k{Z,i}{\Samp {m_{2}} } - \ObsMk k{Z,i}{\Samp {m_{2}}})^\top 
\Big) .
\end{align*}
The random variables $\{ \Samp{1}, \dots, \Samp {M_{k}} \}$ are independent, which implies that \linebreak $\{\Obsk k{Z,i}{\Samp {m_{1}} } - \ObsMk k{Z,i}{\Samp {m_{1}}} \ : m =1,\dots, M_{k}\}$ are independent conditionally on $\F M{k,i}$. Thus, taking the conditional expectation $\CE M_{k,i}$ implies that the $(m_1,m_2)$-terms go to 0 for $m_1\neq m_2$.
With matrix $\Sigma^{(m)}:=\CE M_{k,i}\big[(\Obsk k{Z,i}{\Samp {m} } - \ObsMk k{Z,i}{\Samp {m}})(\Obsk k{Z,i}{\Samp {m} } - \ObsMk k{Z,i}{\Samp {m}})^\top\big]$, it follows that
\begin{align}
\CE M_{k,i}\big[ \normM ki {\psik k_{Z,i} -\psimk k_{Z,i} }^2 \big] \ 
= \frac{ 1 }{M_{k}^2 }
\sum_{m=1}^{M_{k}}
\trace\Big([\tilde p \tilde p^\top] (\Samp {m}) \Sigma^{(m)}\Big) 
\leq 
\frac{ 1 }{M_{k}^2 }\sum_{m=1}^{M_{k}}
\trace\big([\tilde p \tilde p^\top] (\Samp {m}) \big) \trace(\Sigma^{(m)}) ,
\label{eq:problem}
\end{align}
where we have used that  $\trace(A B)\leq \trace(A)\trace(B)$ for any symmetric  non-negative definite matrices $A$ and $B$. 
To continue, we require a bound from above on $\E[\trace\big([\tilde p \tilde p^\top] (\Samp {m}) \big) \trace(\Sigma^{(m)})]$.
Two approaches are available depending on the choice of basis: for general basis  (as for Theorem \ref{thm:0:z:er}), we find almost sure upper bounds for $\trace(\Sigma^{(m)})$ that are uniform in $m$; on the other hand, for the special selection of basis in Theorem \ref{thm:spec:basis}, the intrinsic properties of the basis are used to obtain refined bounds.
\begin{lemma}
\label{lem:gen:bds:1}
For any $k\ge 0$, $i\in \{0,\ldots, 2^k-1\}$, and basis functions chosen as in Definition \ref{def:fin dim:space}, 
\begin{align*}
\normM ki {\psik k_{Z,i} -\psimk k_{Z,i} }^2 & \le 
{\BasDimk k{Z,i} \over \Del k_ i M_k }  \Big\{  |\y {k}{i}{\cdot} - \y {k,M}{i}{\cdot} |_\infty^2 + \sum_{j = \alpha(i) +1}^{2^{k-1}-1} | \z {k-1}{j}{\cdot} -\z {k-1,M}{j}{ \cdot} |_\infty^2 \Del{k-1}_j  \Big\} ,\\
\normM ki {\psik k_{Y,i} -\psimk k_{Y,i} }^2 & \le 
{\BasDimk k{Z,i} \over  M_k }   \sum_{j = \alpha(i) +1}^{2^{k-1}-1} | \z {k-1}{j}{\cdot} -\z {k-1,M}{j}{ \cdot} |_\infty^2 \Del{k-1}_j  .
\end{align*}
\end{lemma}
{\bf Proof.}
We treat the terms $\normM ki {\psik k_{Z,i} -\psimk k_{Z,i} }^2$; the proof for the terms $\normM ki {\psik k_{Y,i} -\psimk k_{Y,i} }^2$ is the same and we exclude it.
Recall the estimate \eqref{eq:problem}.
Thanks to the independence of the Brownian increments, one obtains the equality
\begin{align}
\trace(\Sigma^{(m)}) & = \CE M_{k,i}[ |\Obsk k{Z,i}{\Samp m} - \ObsMk k{Z,i}{\Samp m} |^2] \nonumber \\
& = 
\CE M_{k,i}\big[ (\y {k}{i}{\X {k,m}{i}} - \y {k,M}{i}{\X {k,m}{i}} )^2  \big] 
{ \E[|\DW {k,m}i|^2]  \over  (\Del k_i)^2}  
 \nonumber\\ 
&+   \sum_{j = \alpha(i) +1}^{2^{k-1}-1} 
\CE M_{k,i}\big[|\DW {k,m}{i}|^2 
| \z {k-1}{j}{\X {k-1,m}{j}} -\z {k-1,M}{j}{ \X {k-1,m}{j}} |^2   \big]   { \E[| \DW {k-1,m}j\big|^2 ] \over (\Del k _i)^2} \nonumber \\
& \le  
 {1 \over  \Del k_i}  \Big\{  |\y {k}{i}{\cdot} - \y {k,M}{i}{\cdot} |_\infty^2 + \sum_{j = \alpha(i) +1}^{2^{k-1}-1} | \z {k-1}{j}{\cdot} -\z {k-1,M}{j}{ \cdot} |_\infty^2 \Del{k-1}_j \Big\}
\label{eq:R1:0} 
\end{align}
Now, using $\frac{ 1 }{M}\sum_{m=1}^M [\tilde p \tilde p^\top] (\Samp {m})={\rm Id}_{\R^{\tilde K}}$ and $\tilde K\leq \BasDim{Z,k,i}$, one substitutes the bounds of \eqref{eq:R1:0}
into \eqref{eq:problem} in order to obtain the result. \qed

In fact, one can improve on Lemma \ref{lem:gen:bds:1} if one assumes additional structure on the basis functions.

\begin{lemma}
\label{lem:gen:bds:2}
In addition to the general assumptions, assume \HK \ from Theorem \ref{thm:spec:basis}.
For any $k\ge 0$, $i\in \{0,\ldots, 2^k-1\}$, 
\begin{align}
\E[ \normM ki {\psik k_{Z,i} -\psimk k_{Z,i} }^2] & \le 
{\BasDimk k{Z,i} \delta \over \Del k_ i M_k }  \Big\{  \E[ \lawnorm{\y {k}{i}{\cdot} - \y {k,M}{i}{\cdot} }{k,i}^2] 
+ 8 qC_X^2 2^{-k}  T^\theta \theta^{-1} 
\nonumber\\
& \qquad + 2kq \ln(2) \sum_{j = \alpha(i) +1}^{2^{k-1}-1} \E[\lawnorm{ \z {k-1}{j}{\cdot} -\z {k-1,M}{j}{ \cdot} }{k-1,j}^2 ] \Del{k-1}_j  \Big\} ,
\label{eq:specbas:z:var}
\\
\E[ \normM ki {\psik k_{Y,i} -\psimk k_{Y,i} }^2 ] & \le 
{\BasDimk k{Z,i} \delta \over  M_k }   \sum_{j = \alpha(i) +1}^{2^{k-1}-1} \E[\lawnorm{ \z {k-1}{j}{\cdot} -\z {k-1,M}{j}{ \cdot} }{k-1,j}^2 ]\Del{k-1}_j  
\label{eq:specbas:y:var}.
\end{align}
\end{lemma}

{\bf Proof.}
We give the proof for \eqref{eq:specbas:z:var}; the proof for \eqref{eq:specbas:y:var} is analogous (and simpler).
Starting from \eqref{eq:problem}, we apply the  method of \cite{benz:gobe:13}.
For the convenience of the reader, we translate the notation of \cite{benz:gobe:13} to our setting:
the functions $f_j$ are equivalent to our $\tilde p_{j}$, the $j$-th component of the vector $\tilde p$, {whence $\trace\big([\tilde p \tilde p^\top] (\Samp {m}) \big) = \sum_{j=1}^{\BasDimk k{Z,i}} (f_j)^2$}; 
$H^{\alpha^\star}_m$ is equivalent to our $\Obsk k{Z,i}{\Samp {m} }$; $X$ is equivalent to our $\cX$ and $X_m$ is equivalent to our $\Samp m$.
{Assume that $\P_{\X ki}(A_{Z,k,i,j}) \ge \delta/\BasDimk k{Z,i}$ for all $j$.}
Using the  conditioning argument of \cite[case (b) on page 14]{benz:gobe:13}, it follows that
\begin{align*}
&\E[\normM ki {\psik k_{Z,i} -\psimk k_{Z,i} }^2]  
\; \leq \;
\E \Big[
\frac{ 1 }{M_{k}^2 }\sum_{m=1}^{M_{k}}
\trace\big([\tilde p \tilde p^\top] (\Samp {m}) \big) \CE M_{k,i}[ \trace(\Sigma^{(m)})] \Big] \\
& \; { \le {1 \over M_k} \e[ \var(H^{\alpha^\star}(X) | X ) \sum_{j=1}^{\BasDimk k{Z,i}} (f_j)^2 ] 
\qquad \text{(in the equivalent notation of \cite{benz:gobe:13})} }\\
&\; \le \sum_{j=1}^{\BasDimk k{Z,i}} { \E[ |\Obsk k{Z,i}{\cX} - \ObsMk k{Z,i}{\cX} |^2 \1_{\X ki \in A_{Z,k,i,j}}] \over M_k \P(\X {k}i \in A_{Z,k,i,j}) } 
{\le {\BasDimk k{Z,i}\delta \over M_k} \E[ |\Obsk k{Z,i}{\cX} - \ObsMk k{Z,i}{\cX} |^2 ].}
\end{align*}
To complete the proof, we  obtain upper bounds on $\E[ |\Obsk k{Z,i}{\cX} - \ObsMk k{Z,i}{\cX} |^2]$:
\begin{align*}
\E[ |\Obsk k{Z,i}{\cX}&  - \ObsMk k{Z,i}{\cX} |^2] 
 = 
\E[\big[ (\y {k}{i}{\X {k}{i}} - \y {k,M}{i}{\X {k}{i}} )^2  \big] 
{ \E[|\DW {k}i|^2]  \over  (\Del k_i)^2}  
 \nonumber\\ 
&+   \sum_{j = \alpha(i) +1}^{2^{k-1}-1} 
\E[|\DW {k}{i}|^2 
| \z {k-1}{j}{\X {k-1}{j}} -\z {k-1,M}{j}{ \X {k-1}{j}} |^2   ]   { \E[| \DW {k-1}j\big|^2 ] \over (\Del k _i)^2} .
\end{align*}
There is an interdependency issue between $\DW {k}{i}$ and $| \z {k-1}{j}{\X {k-1}{j}} -\z {k-1}{j}{ \X {k-1}{j}} |$ that we now treat;
note that this interdependency does not arise when dealing with $\E[|\Obsk k{Y,i}{\cX}  - \ObsMk k{Y,i}{\cX} |^2]$.
Since  $\DW {k}i$ has $q$ independent components, each with Gaussian distribution with mean 0 and variance $\Del k_i$, these components are each equal in law to $\sqrt{\Del k_i} \cN$, where $\cN$ has a Gaussian distribution with mean 0 and variance 1.
Calculating the expectation by integration-by-parts then using Mill's inequality implies, for any $R>0$, that
\begin{align*}
\E[|\cN|^2 \1_{|\cN| > \sqrt R}] &= 2\big(\P(\cN > \sqrt R)(R + 1) - \frac{\sqrt R e^{-R/2} }{\sqrt{2\pi}} \big)
 \le 2\P(\cN > \sqrt R)(R + 1 - R) \le 2e^{-R/2}.
\end{align*}
Now, using the decomposition $\DW {k}i = \cT_{\sqrt{\Del k_i R }} (\DW {k}i ) + (\DW {k}i - \cT_{\sqrt{\Del k_i R }} (\DW {k}i ) )$ and the almost sure bounds on the $z$ terms from Corollary \ref{cor:bd:z} and Algorithm \ref{alg:mlmc}, it follows that
\begin{align*}
\E[|\DW {k}{i}|^2 &
| \z {k-1}{j}{\X {k-1}{j}}-\z {k-1,M}{j}{ \X {k-1}{j}} |^2   ] \\
& \le
q R \Del k_i \E[ | \z {k-1}{j}{\X {k-1}{j}} -\z {k-1,M}{j}{ \X {k-1}{j}} |^2   ] + { 4q C_X^2 \Del k_i \over (T-\tk {k-1}j)^{1-\theta} } \E[|\cN|^2 \1_{|\cN| > R}] \\
& \le
qR \Del k_i \E[\lawnorm{ \z {k-1}{j}{\cdot} -\z {k-1,M}{j}{ \cdot} }{k-1,j}^2 ] + {8qe^{-R/2} C_X^2 \Del k_i \over (T-\tk {k-1}j)^{1-\theta} } \\
& \le
qR \Del k_i \E[\lawnorm{ \z {k-1}{j}{\cdot} -\z {k-1,M}{j}{ \cdot} }{k-1,j}^2 ] + {8qe^{-R/2} C_X^2 \Del k_i \over (T-\tk {k-1}j)^{1-\theta} } .
\end{align*}
The proof is completed by selecting $R = \ln(2^{2k})$.
\qed

To complete the estimate of \eqref{eq:main:y:3} and \eqref{eq:main:z:3}, it remains only to bound $\E[ \normM ki { \CE{M}_{k,i}[\psik k_{Y,i}] -\psik k_{Y,i} }^2]$ and $\E[ \normM ki { \CE{M}_{k,i}[\psik k_{Z,i}] -\psik k_{Z,i} }^2]$.

\begin{proposition}
\label{prop:bar psi}
In addition to the general assumptions, suppose that either \HXppp \ (from Theorem \ref{thm:0:z:er}) or \HK \ (from Theorem \ref{thm:spec:basis}) is in force.
Then, for all $k\ge0$ and $i\in\{0,\ldots,2^k-1\}$,
\begin{align*}
& \E[ \normM ki { \CE{M}_{k,i}[\psik k_{Y,i}] -\psik k_{Y,i} }^2] \le 
C_1 { 2 \times 2^{-k} \BasDimk k{Y,i} \over   M_k} \left\{   3 C_X^2 + (2+q)\right\} 
 \quad \text{and}\\
& \E[ \normM ki { \CE{M}_{k,i}[\psik k_{Z,i}] -\psik k_{Z,i} }^2] \le 
C_1 {  6  \BasDimk k{Z,i}  (2 +5T^{1-\theta} )C_X^2 \over c_X M_k (T -\tk ki)^{1-\theta} } 
 + 2 C_1   {(2+q) \BasDimk k{Z,i} C_{X} \over c_X M_k}.
\end{align*}
where $C_1 = \delta$ if \HK \ holds and $ C_1 =1$ if \HXppp \ holds.
\end{proposition}

{\bf Proof.}
We will use Proposition \ref{prop:ls:reg:properties}(iv).
For  $\vecx x = (x_0 ,\ldots,x_{2^k}) \in \R^{(2^k + 1)\times d}$, $\vecx {\bar x} = (\bar x_0 ,\ldots,\bar x_{2^{k-1}}) \in \R^{(2^{k-1} +1)\times d}$, 
$\vecx{w} = (w_{0}, \ldots, w_{2^{k} -1}) \in \R^{2^{k}  \times q}$,
define $h(\vecx x ,\vecx {\bar x} , \vecx w) := x_i$; $h$ is a Borel measurable function, and $h(\Samp m) = \X {k,m}i$.
Denote by $\cH$ the  $\sigma$-algebra $ \sigma\big(h_i(\cX_m) \ : \ m =1 , \ldots , M_k \big) $, which is equal to $\sigma(\X {k,m}i \ : \ m =1 , \ldots , M_k ) $,
and by $\cG$ the $\sigma$-algebra $\F *{k-1} \vee  \sigma(\X {k,m}j \ : \ j<i, \ m =1 , \ldots , M_k )$;
 then $\Obsk k{Y,i}\cdot$ and $\Obsk k{Z,i}\cdot$ are $\cG \otimes \cB(\R^l)$-measurable, and $\cG \vee \cH$ is equal to $\F M{k,i}$.
Since $\psimk k_{Y,i}$ (resp. $\psimk k_{Z,i}$) solves $\OLS(\Obsk k{Y,i}{\cdot} , \LinSpacek k{Y,i} , \emeas k)$ 
(resp. $\OLS(\Obsk k{Z,i}{\cdot} , \LinSpacek k{Z,i} , \emeas k)$)
it only remains to find suitable (deterministic) upper bounds for expectation of
$ \Psi_{Y,k,i} := \CE M_{k,i}[| \Obsk k{Y,i}{\Samp m} - \CE M_{k,i} [\Obsk k{Y,i}{\Samp m}] |^2 ]$
 ({resp. }  $\Psi_{Z,k,i} := \CE M_{k,i}[| \Obsk k{Z,i}{\Samp m} - \CE M_{k,i} [\Obsk k{Z,i}{\Samp m}] |^2 ]$ ) 
to allow us to apply Proposition \ref{prop:ls:reg:properties}(iv).
The technique is  similar for both $ \Psi_{Y,k,i}$ and $ \Psi_{Z,k,i}$, so we include the proof for the latter only.
The strategy will be the following: first, we assume that the Markov chains $\X {k,m}{}$ and $\X {k-1,m}{}$ are have deterministic values $x$ at $\tk ki$ and $\bar x$ at $\tk {k-1}{\alpha(i)}$, respectively;
then we decompose the upper bound on  $ \Psi_{Z,k,i}$ by introducing the diffusion processes $\X {\tk ki,x,m}{}$ and $\X {\tk {k-1}{\alpha(i)},\bar x,m}{}$; 
eventually, we fix $x = \X {k,m} i$ and $\bar x = \X{k-1,m}{\alpha(i)}$ (which does not pose difficulties due to the use of the conditional expectation $\CE M_{k,i}[\cdot]$ throughout) to obtain the final bounds.

$\blacktriangleright${\bf Step 1 (fixing the initial value of the Markov chain at $\tk ki$ and $\tk {k-1}{\alpha(i)}$):}
Observe that the random variable $\Obsk k{Z,i}{\Samp m}$ depends on the sample path $\Samp m$ only through the values 
\[
(\X {k,m}i, \ldots, \X {k,m}{2^k} , \X {k-1,m}{\alpha(i)+1},\ldots , \X {k-1,m}{2^{k-1}}, \DW{k,m}i,\ldots, \DW {k,m}{2^k-1}),
\]
i.e., it does not depend on the path $\X {k,m}{}$, $\X {k-1,m}{}$ and $\DW {k,m}{}$ before the time $\tk ki$.
Letting $x,\bar x\in \R^d$,
 we define 
\[
\Samp {m,i}(x,\bar x) :=  (\X {k,m,i,x}i, \ldots, \X {k,m,i,x}{2^k} , \X {k-1,m,\alpha(i),\bar x}{\alpha(i)+1},\ldots , \X {k-1,m,\alpha(i),\bar x}{2^{k-1}}, \DW{k,m}i,\ldots, \DW {k,m}{2^k-1}).
\]
One can then write $\Psi_{Z,k,i} =  \Psi_{Z,k,i}(\X {k,m}i, \X{k-1,m}{\alpha(i)})$, where
\[
 \Psi_{Z,k,i}(x,\bar x) :=  \CE M_{k,i}[| \Obsk k{Z,i}{\Samp {m,i}(x,\bar x) } - \CE M_{k,i} [\Obsk k{Z,i}{\Samp {m,i}(x,\bar x) }] |^2 ].
\]
{
}
$\blacktriangleright${\bf Step 2 (decomposition with intermediate discrete BSDE):} 
Let $\X {t,x,m}{}$ ($m \in \{1,\ldots, M_k \}$) be the simulation of the diffusion started at time $t$ with value $x$ generated with the same path of the Brownian motion as the increments $\DW {k,m}i$.
Recall the discrete BSDE $(\tily k{},\tilz k{})$ from section \ref{section:0:general} and define
\[ \tily{ k,m}j := \E^k_j [\Phi(\X {\tk ki,x,m}T)] \quad \text{and}  \quad \Del k_j \tilz {k,m}j := \E^k_i[(\DW {k, m}i)^\top \Phi(\X {\tk ki,x,m}T)]\] 
for $j \in \{0,\ldots,2^k-1\}$. 
For the coarse grid $\pik {k-1}$,
  define
\begin{align*} 
\tily {k-1,m}j := \E^{k-1}_j [\Phi(\X {\tk {k-1}{\alpha(i)},\bar x,m}T)] \quad \text{and}  \quad \Del {k-1}_j \tilz {k-1,m} j := \E^{k-1}_j[\DW {k-1,m}j \Phi(\X {\tk{k-1}{\alpha(i)},\bar x,m}T)] 
\end{align*}
for $j \in \{0,\ldots,2^{k-1}-1\}$.
We use these processes to decompose $\Obsk k{Z,i}{\Samp m(x,\bar x)}$ into two expressions:
\begin{align*}
\Obsk k{Z,i}{\Samp {m,i}(x,\bar x)} & = 
\left\{
\begin{array}{l}
\displaystyle{\DW {k,m}i \over \Del k_i} \Big\{ \Phi(\X {k,m,i,x}N) - \Phi(\X {\tk ki,x,m}T) - (\y ki{\X {k,m,i,x}i}- \tily {k,m}i) \\
\qquad \quad- \sum_{j=\alpha(i) +1 }^{2^{k-1}-1} ( \z {k-1}j{\X{k-1,m,\alpha(i),\bar x}j} - \tilz{k-1,m}j )\DW{k-1,m}j \Big \} 
\end{array}
\right\} 
\\
& \quad + \left\{  {\DW {k,m}i \over \Del k_i} \big\{ \Phi(\X {\tk ki,x,m}T) - \tily {k,m}i - \sum_{j=\alpha(i) +1 }^{2^{k-1}-1} \tilz{k-1,m}j \DW{k-1,m}j \big\}  \right\} \\
& =: A_1(x,\bar x) + A_2(x,\bar x).
\end{align*}
The trivial inequality $(x + y)^2 \le 2 x^2 +2 y^2$ for all real $x$ and $y$  then yields
\[
\Psi_{Z,k,i}(x,\bar x) \le 
2\CE M_{k,i}[A_1(x,\bar x)^2] + 2\CE M_{k,i}[A_2(x,\bar x)^2].
\]

$\blacktriangleright$ {\bf Step 3 (bound on 
$\CE M_{k,i} [A_1(x,\bar x)^2]$).}
Using  the Cauchy-Schwarz inequality, we have
\begin{align*}
&\CE M_{k,i} [ |\z {k-1}j{\X{k-1,m,\alpha(i),\bar x}j} - \tilz{k-1,m}j |^2] \\
& = \frac 1{(\Del {k-1}_j)^2} \CE M_{k-1,i} [ | \CE M_{k,j} [ (\DW {k-1,m}{j})^\top (\y {k-1}{j+1}{\X{k-1,m,\alpha(i),\bar x}{j+1}} - \tily{k-1,m}{j+1} )  ]|^2] \\
&\le \frac q{\Del {k-1}_j} \Big \{ \CE M_{k,i} [ |\y {k-1}{j+1}{\X{k-1,m,\alpha(i),\bar x}{j+1}} - \tily{k-1,m}{j+1} |^2]
-   \CE M_{k,i} [ |\CE M_{k-1,j}[\y {k-1}{j+1}{\X{k-1,m,\alpha(i),\bar x}{j+1}} - \tily{k-1,m}{j+1} ] |^2]  \Big \} .
\end{align*}
Observe that
\(
\CE M_{k-1,j}[\y {k-1}{j+1}{\X{k-1,\alpha(i),\bar x,m}{j+1}} - \tily{k-1,m}{j+1} ]  = \y {k-1}{j}{\X{k-1,\alpha(i),\bar x,m}{j}} - \tily{k-1,m}{j}. 
\)
Then, a shift of  summation {indicies} gives
\begin{align*}
\sum_{j=\alpha(i) +1 }^{2^{k-1}-1} & \CE M_{k,i} [ |\z {k-1}j{\X{k-1,m,\alpha(i),\bar x}j} - \tilz{k-1,m}j |^2]\Del {k-1}_{j} \\
&\le \sum_{j=\alpha(i) +1 }^{2^{k-1}-1} q \Big \{ \CE M_{k,i} [ |\y {k-1}{j+1}{\X{k-1,m,i,\bar x}{j+1}} - \tily{k-1}{j+1} |^2]
-   \CE M_{k,i} [ |\CE M_{k-1,j}[\y {k-1}{j+1}{\X{k-1,m,i,\bar x}{j+1}} - \tily{k-1}{j+1} ] |^2]  \Big \} \\
& \le \sum_{j=\alpha(i) +1 }^{2^{k-1}-1} q \Big \{ \underbrace{  \CE M_{k,i} [ |\y {k-1}{j}{\X{k-1, m,i,\bar x}{j}} - \tily{k-1,m}{j} |^2]
-   \CE M_{k,i} [ |\CE M_{k-1,j}[\y {k-1}{j+1}{\X{k-1,m,i,\bar x}{j+1}} - \tily{k-1,m}{j+1} ] |^2]     }_{=\; 0} \Big \} \\
& \qquad +q\CE M_{k,i}[|\Phi(\X {k-1,m,\alpha(i),\bar x}N) - \Phi(\X {\tk {k-1}{\alpha(i)},\bar x,m}T)|^2] .
\end{align*}
Therefore,
using the independence of the Brownian increments,
 one obtains the upper bound
\begin{align}
\CE M_{k,i}[A_1 &(x,\bar x)^2]   = 
\displaystyle{1 \over \Del k_i} \Big\{ \CE M_{k,i}[|\Phi(\X {k,m,i,x}N) - \Phi(\X {\tk ki,x,m}T)|^2] 
+ \CE M_{k,i}[|\y ki{\X {k,m, i,x}i}- \tily ki|^2] 
\nonumber\\
& \qquad \quad + \sum_{j=\alpha(i) +1 }^{2^{k-1}-1} \CE M_{k,i} [ |\z {k-1}j{\X{k-1,m,i,\bar x}j} - \tilz{k-1,m}j |^2]\Del {k-1}_{j} \Big \} 
\nonumber\\
& \le \displaystyle{1 \over \Del k_i} \Big\{ 2 \CE M_{k,i}[|\Phi(\X {k,m,i,x}N) - \Phi(\X {\tk ki,x,m}T)|^2] + q\CE M_{k,i}[|\Phi(\X {k-1,m,\alpha(i),\bar x}N) - \Phi(\X {m,\tk {k-1}{\alpha(i)},\bar x}T)|^2] \Big\} 
\label{eq:bd:a1}
\end{align}

It follows from assumption on the Markov chains \HXp(i) and the assumption on the time-grids \Ht(iii)  that the terms in parenthesis in \eqref{eq:bd:a1} can be bounded by $C_{X} \max_{i} \{2^{-k} + q 2^{-(k-1)}\}$.
{Hence,} 
\[
\CE M_{k,i}[A_1(x,\bar x)^2] \le {(2+q) C_{X} 2^{-k} \over \Del k_{i} } \le {(2+q) C_{X} \over c_X} .
\]

$\blacktriangleright$ {\bf Step 4 (bound on 
$\CE M_{k,i} [A_2(x,\bar x)^2]$).}
Using  equality \eqref{eq:disc:bsde:0} and  Lemma \ref{lem:disc:bsde:0}, we have
\[
 \sum_{j = \alpha(i) +1}^{2^{k-1}-1}\tilz {k,m}j \DW{k-1,m}j 
= \Phi (\X {\tk {k-1}{\alpha(i)},\bar x,m}T ) - \tily {k-1,m}{\alpha(i) +1}   - \sum_{j = \alpha(i) +1}^{2^{k-1}-1}\DL {k-1,m}{j}\]
where $\DL {k-1,m}j := \int_{\tk {k-1}j}^{\tk {k-1}{j+1}} (\zn {\tk {k-1}{\alpha(i)},\bar x,m}t  - \tilz {k-1,m}j )^\top dW^{(m)}_t$ {and} $\zn {\tk {k-1}{\alpha(i)},\bar x,m}{}$ {is} 
the process given in \HX(iii) with $\X{\tk {k-1}{\alpha(i)},\bar x,m}{}$ in the place of $\X{\tk {k-1}{\alpha(i)},\bar x}{}$.
Substituting this into the definition of $A_{2}(x,\bar x)$, it follows that
\[
A_{2}(x,\bar x) = {(\DW ki )^\top \over \Del k_i} \Big\{  \Phi(\X {\tk ki,x,m}T) - \Phi (\X {\tk {k-1}{\alpha(i)},\bar x,m }T ) - (\tily {k,m}i - \tily{k-1,m}{\alpha(i)+1}) 
+\sum_{j = \alpha(i) +1}^{2^{k-1}-1}\DL {k-1,m}{j}\Big\} .
\]
Now we square and take expectations, and treat the terms in $\Delta L$, $\Phi$ and $\tilde y$ individually.

To treat the terms in $\Delta L$,
 we apply property \Ht(iv) to obtain  {that}
\(
\sum_{j = \alpha(i) +1 } ^{2^{k-1}}\CE M_{k,i} [|\DL {k-1,m}{j}|^2]
\)
{is bounded by}  $ 2C_X 2^{-k}$; this upper bound is independent of the starting value $\bar x$ of $\X {\tk {k-1}{\alpha(i)},\bar x,m }{}$.
In order to treat the terms in{volving} $\Phi$, we apply assumption \HX(ii) to obtain
\begin{align}
\CE M_{k,i}[| \Phi(\X {\tk ki,x,m}T) - \Phi (\X {\tk {k-1}{\alpha(i)},\bar x,m}T ) |^2]  \le C_X | x -  \X {\tk {k-1}{\alpha(i)},\bar x,m}{\tk ki} |^2.
\label{eq:Phi:bd}
\end{align} 
\begin{remark}
\label{rem:gen:hold}
In \eqref{eq:Phi:bd}, one can see the impact of {condition} \HX(ii); 
it is needed to obtain the upper bound $O(2^{-k})$ of the terms in $\Phi$, once we put $x = \X {k,m}{i}$, $\bar x = \X {k-1,m}{\alpha(i)}$, and  take expectations.
\end{remark}

Finally, to treat the terms in $\tilde y$, we use that $\tily{k,m}i - \tily{k-1,m}{\alpha(i)+1}$ is equal to
\begin{equation*}
\label{eq:tily:bd:1}
  \tily{k,m}i \pm \tily {k,m}{2(\alpha(i)+1)} - \tily{k-1,m}{\alpha(i)+1} = (\tily {k,m}i - \tily{k,m}{2(\alpha(i)+1)}) + \CE M_{k,2(\alpha(i) +1)}[\Phi(\X {\tk ki,x,m}T) - \Phi (\X {\tk {k-1}{\alpha(i)},\bar x,m}T ) ].
\end{equation*}
The terms in $\Phi$ are treated as in \eqref{eq:Phi:bd}.
We further expand $(\tily ki - \tily{k}{2(\alpha(i)+1)})$ using \eqref{eq:disc:bsde:0}
\begin{align*}
\tily{k,m}i - \tily{k,m}{2(\alpha(i)+1)} & = \sum_{j=i}^{2(\alpha(i) +1 )} \{ \tilz kj \DW {k,m}j + \DL {k,m}j\}.
\end{align*}
Squaring and taking conditional expectations, we obtain from \Ht(iv) and Corollary \ref{cor:bd:z} that
\begin{align*}
\CE M_{k,i}[|\tily {k,m}i - \tily{k,m}{2(\alpha(i)+1)} |^2] & = \sum_{j=i}^{2(\alpha(i) +1 )} \CE M_{k,i}[ | \tilz kj|^2 | \DW kj|^2 + |\DL kj|^2] 
 \le  {2 C_X^2 2^{-k} \over (T -\tk ki)^{1-\theta} } + C_X2^{-k} .
\end{align*}
To conclude {\bf Step 4}, we combine the above upper bounds to obtain
\begin{align*}
\Del k_i \CE M_{k,i} [A_2(x, \bar x)^2] \le 3 \times (   3C_X | x -  \X {\tk {k-1}{\alpha(i)},\bar x,m}{\tk ki} |^2 + {2 C_X^2 2^{-k} \over (T -\tk ki)^{1-\theta} } + 2 C_X2^{-k} ).
\end{align*}
Therefore, plugging $x = \X {k,m}i$ and $\bar x = \X{k-1,m}{\alpha(i)}$, we obtain
\begin{align}
 \Del k_i \CE M_{k,i} & [A_2(\X {k,m}i,\X{k-1,m}{\alpha(i)})^2] 
 \le 
9 C_X  | \X {k,m}{i}  -\X{k-1,m}{\alpha(i)}|^2 + {6 C_X^2 2^{-k} \over (T -\tk ki)^{1-\theta} } + 6 C_X2^{-k}
\label{eq:equal:X's}
\end{align}

{\bf Concluding the proof.}
The proof of the proposition under \HXppp \ is now completed by observing that $| \X {k,m}{i}  -\X{k-1,m}{\alpha(i)}|^2 = 0$ in \eqref{eq:equal:X's}, piecing together the estimates obtained in {\bf Steps 1-4} on $\CE M_{k,i}[| \Obsk k{Z,i}{\Samp m} - \CE M_{k,i} [\Obsk k{Z,i}{\Samp m}] |^2 ]$ to find that there is a deterministic bound, and applying Proposition \ref{prop:ls:reg:properties}(iv).
On the other hand, if \HK \ were in force, we again (as in proof of Lemma \ref{lem:gen:bds:2})
use the conditioning arguments of \cite[case (b) on page 14]{benz:gobe:13} in order to obtain 
\begin{align*}
\E[ \normM ki { \CE{M}_{k,i}[\psik k_{Z,i}] -\psik k_{Z,i} }^2] \le {  \delta \BasDimk k{Z,i} \E \big[\CE M_{k,i}[| \Obsk k{Z,i}{\Samp 1} - \CE M_{k,i} [\Obsk k{Z,i}{\Samp 1}] |^2 ] \big] \over M_k}.
\end{align*}

{By combining}
 the estimates obtained in {\bf Steps 1-4} on
 $\CE M_{k,i}[| \Obsk k{Z,i}{\Samp 1} - \CE M_{k,i} [\Obsk k{Z,i}{\Samp 1}] |^2 ]$, 
substitut{ing} them into the above inequality, 
 {one sees} that
{ the expectation}
 \( \E[ \normM ki { \CE{M}_{k,i}[\psik k_{Z,i}] -\psik k_{Z,i} }^2]\) 
 {is bounded by }
\begin{align}
   6 \delta { \BasDimk k{Z,i}  \over M_k} & \left(3 C_X { \E [ | \X {k}{i}  -\X{k-1}{\alpha(i)}|^2  ] \over \Del k_i} + {(2 +2T^{1-\theta} )C_X^2 \over  c_X (T -\tk ki)^{1-\theta} } \right) \nonumber
 + 2  \delta {(2+q) \BasDimk k{Z,i} C_{X} \over c_X M_k} \nonumber
\\
& \le 
 {  6 \delta \BasDimk k{Z,i}  (2 +5T^{1-\theta} )C_X^2 \over c_X M_k (T -\tk ki)^{1-\theta} } 
 + 2  \delta {(2+q) \BasDimk k{Z,i} C_{X} \over c_X M_k}
 \label{eq:coup}
\end{align}
where we have used \HXp(iii) \ in the last inequality, and \Ht(iii) for the bound $\Del k_i \ge c_X 2^{-k}$.
\qed

\begin{remark}
\label{rem:coup}
We see in equation \eqref{eq:coup} the impact of assumption \HXp(iii);
were a lower rate of convergence assumed, the overall rate of convergence with respect to $k$ of the upper bound in Theorem \ref{thm:spec:basis} would be lower.
\end{remark}

The proof 
is completed by estimating the terms in (\ref{eq:main:y:3} -- \ref{eq:main:z:3}) using Lemma \ref{lem:gen:bds:1} (resp. Lemma \ref{lem:gen:bds:2}) and Proposition \ref{prop:bar psi}, and substituting the estimates into (\ref{eq:MC:Y:1} -- \ref{eq:MC:Z:1}).

\subsection{{Computational examples}}
\label{section:num:ml}

The computational examples in this section illustrate and
 compare the actual errors and efficiency of different simulation schemes to the  BSDE with zero  generator \eqref{eq:lin} 
 to support the results of the theoretical analysis based on error estimates.
We consider cases of BSDE with analytically known solutions in order to investigate the actual {global mean squared errors (MSE)}
 of the approximate solutions for the respective approximation schemes.
{The MSEs 
are computed by Monte Carlo on a fine time grid in the same way as the global error in  in \eqref{eq:glob:er}.
The overall MSE is the sum of the MSEs with respect to the $Y$ and the $Z$ components, 
corresponding to the first (with squared maximum over time) respectively second (time-weighted)  summand.}

\subsubsection{Sine payoff}{ At first let us consider a $q=1$ dimensional example with
$T=1$, $X=W$, and terminal condition $\Phi(x) = \sin(x)$.
For a regression basis, we took Hermite polynomials up to degree $7$, adjusted for time so that $\{1, p_1(t,W_t), \ldots, p_7(t,W_t)\}$ is orthonormal for all time $t>0$, i.e. $K = 8$.
We run the multilevel (ML) scheme with $M_k = 40 \times K \times 2^k $ simulations at {final} level $k$, while
at any lower level $j<j+1\le k$ the number of simulations  $M_{j}=2 M_{j+1}$ doubles.
The overall complexity for ML up to level $k$ is therefore $\cC_{ML} = O(k \times 2^{2k})$.
}
\graphicspath{{./Plots/}}
\begin{figure}[h]
\includegraphics[width=13cm]{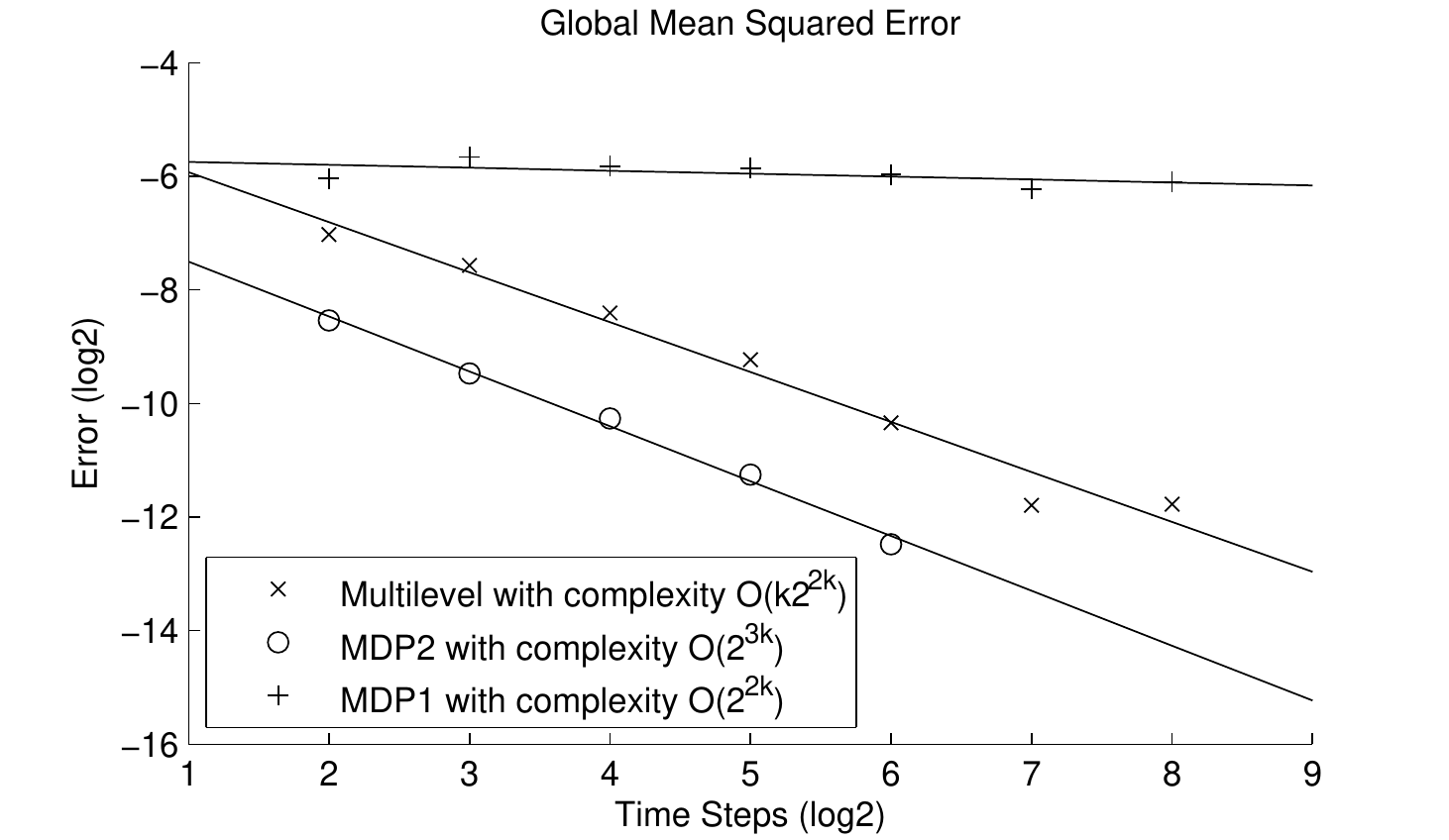}
\caption{Sine payoff: $\log_2$ {Mean squared error} versus $\log_2(N)$ for ML and MDP1,2
\label{fig:mlmdp:sine}}
\end{figure}
{For comparison, we run two instances of the MDP scheme: At first (MDP1) with $M_k = 40 \times K \times 2^k $ simulations, 
and then (MDP2) with a much higher number $M_k = 40 \times K \times 2^{2k} $  of simulations.
The complexity in the first case is $\cC_{MDP,1} = O( 2^{2k})$ and $\cC_{MDP,2} = O(2^{3k})$ in the second.
}
Figure \ref{fig:mlmdp:sine} shows the log of the global mean squared 
errors (MSE) of the MDP1, MDP2 and and the ML scheme vs. 
the log of the number of time steps $k=\log_2 N$.
The respective regression lines are (ML) $-0.88 x - 5.0$, (MDP1) $-0.05 x - 5.7$, respectively (MDP2) $-0.97-6.5$.
{
This example supports the results of Theorems \ref{thm:0:z:er} and \ref{thm:spec:basis} and the subsequent complexity analysis in Section \ref{section:ML:err an};
{indeed, one sees} that one needs to have $2^k$-times as many simulations for MDP as for multilevel to achieve {a}  convergence rate {of about}  $-1$.
Moreover, {one sees that results} the {computational results indicate, that efficiency gains as stated  in} Theorem \ref{thm:spec:basis} 
may  {be obtained, beyond the assumptions of the theorem, for} 
  a wider class of basis functions than allowed by \HK, as the basis functions used in this example do not satisfy this condition.
}

\subsubsection{{A multi-dimensional example}}

{
Let the forward process be a Brownian motion $X=W$ in dimension $q=3$
and consider the  terminal condition $\Phi(X_T)=\prod_1^q X^i_T$ with $T=1$.  
{This is beyond the assumptions used in the complexity analysis of Section \ref{section:ML:err an},
{ the boundedness assumption \Hg, and the Lipschitz assumption \HX(iv)}.}
}
{
This example is to compare the the {MSE} separately in the contributions of the $Y$- and the $Z$-part
  of the multilevel (ML) and the MDP scheme. 
Moreover, we compare two different sets of regression bases. 
The first regression basis (`indicator') consists of  indicator functions on {equiprobable} hypercubes of a
partition of $\R^3$ into $K=8^3=512$ sets. 
The second  regression basis (`linear')  
{consists of}
{ functions,}
 {each being affine within one hypercube} of a partition of ${\R}^3$ into $K=5^3$ sets {and vanishing}  {outside}.
The linear basis contains $4\times5^3=500$ of regression functions{, so that both bases have about the same size}.
Number simulations is $M=2\times10^6$ for both schemes.
{For ML, the same number of simulations is used at each level.}

\begin{table}[h]
\centering
\caption{\textbf{MSE for $Y$ and $Z$} in dimension 3}
\label{table:multid}
\small
\begin{tabular}{lrrrrrrrrrr}
$\text{log}_{2}(N)$
&\textbf{2}&\textbf{3}&\textbf{4}&\textbf{5}&\textbf{6}&\textbf{7}\\
\hline
ML Y (linear)&0.1528&0.1266&0.1215&0.1194&0.1183&0.1190\\
ML Z (linear)&0.0334&0.0184&0.0160&0.0157&0.0166&0.0185\\
MDP Y (linear)&0.1578&0.1316&0.1253 &0.1236&0.1231   &0.1222\\
MDP Z (linear)&0.0358 &0.0245 &0.0301 &0.0462 &0.0786    &0.1441\\
\hline
ML Y (indicator)&0.5815&0.5454&0.5356&0.5331&0.5310&0.5306\\
ML Z (indicator)&0.1509&0.1219&0.1148&0.1135&0.1154&0.1210\\
MDP Y (indicator)&0.5865 &0.5465 &0.5351 &0.5318 &0.5297   &0.5297\\
MDP Z (indicator)&0.1514 &0.1253 &0.1230 &0.1330 &0.1550   &0.2044\\
\end{tabular}
\end{table}

The table of global mean squared errors, Table \ref{table:multid}, shows that the  multilevel scheme  achieves lower 
errors for the $Z$-part for finer time grids, whereas errors for the $Y$-part are similar. 
The multilevel scheme shows a higher reduction of error in $Z$ for the linear basis.
{Error reduction by factors beyond $1/2$ for $k=\log_2 N\ge 4$ are  significant,
 even when noting that the computational cost for MDP with $N=2^{k+1}$ time steps 
are basically equal to that for ML with $N=2^{k}$ steps  at final level $k$.
 Compared  to MDP,  errors (in $Z$) for the multilevel scheme begin to increase at a later stage $k=\log_2 N$ 
and increase at a much milder rate, regardless of the choice of the basis;}
this is best seen by comparing Figure~\ref{graph:multidim:a} 
with Figure~\ref{graph:multidim:b}. 
This effect 
{fits with the results of} Theorems \ref{thm:0:z:er} and \ref{thm:spec:basis}, 
which state that the error of multilevel scheme is 
less affected by the number of time points than the MDP scheme \eqref{eq:er:MDP}; 
the error of the multilevel scheme should increase only logarithmically with {$N$.}
 and indeed Figure~\ref{graph:multidim:a} {shows an error curve increasing only mildly at large $k$. 
Note that for a fixed number of simulations, as here, it is inevitable that statistical errors increase and take over at some stage;
such simply means that more simulations would be required for larger $k=\log_2 N$.
The advantages of the linear over the indicator basis can understood in the sense that in this example the  {bias} from $L^2$ projection on the function 
space spanned by this basis is smaller, whereas the indicator basis would require a finer partition to achieve the same.} 
This example
 shows that efficiency gains from multilevel can be realized in actual computations;
and moreover indicates that efficiency gains may be expected in a more general context 
beyond the specific assumptions  required in Section~\ref{section:ML:err an}. 

\begin{figure}[h]
\centering
\includegraphics[width=11cm]{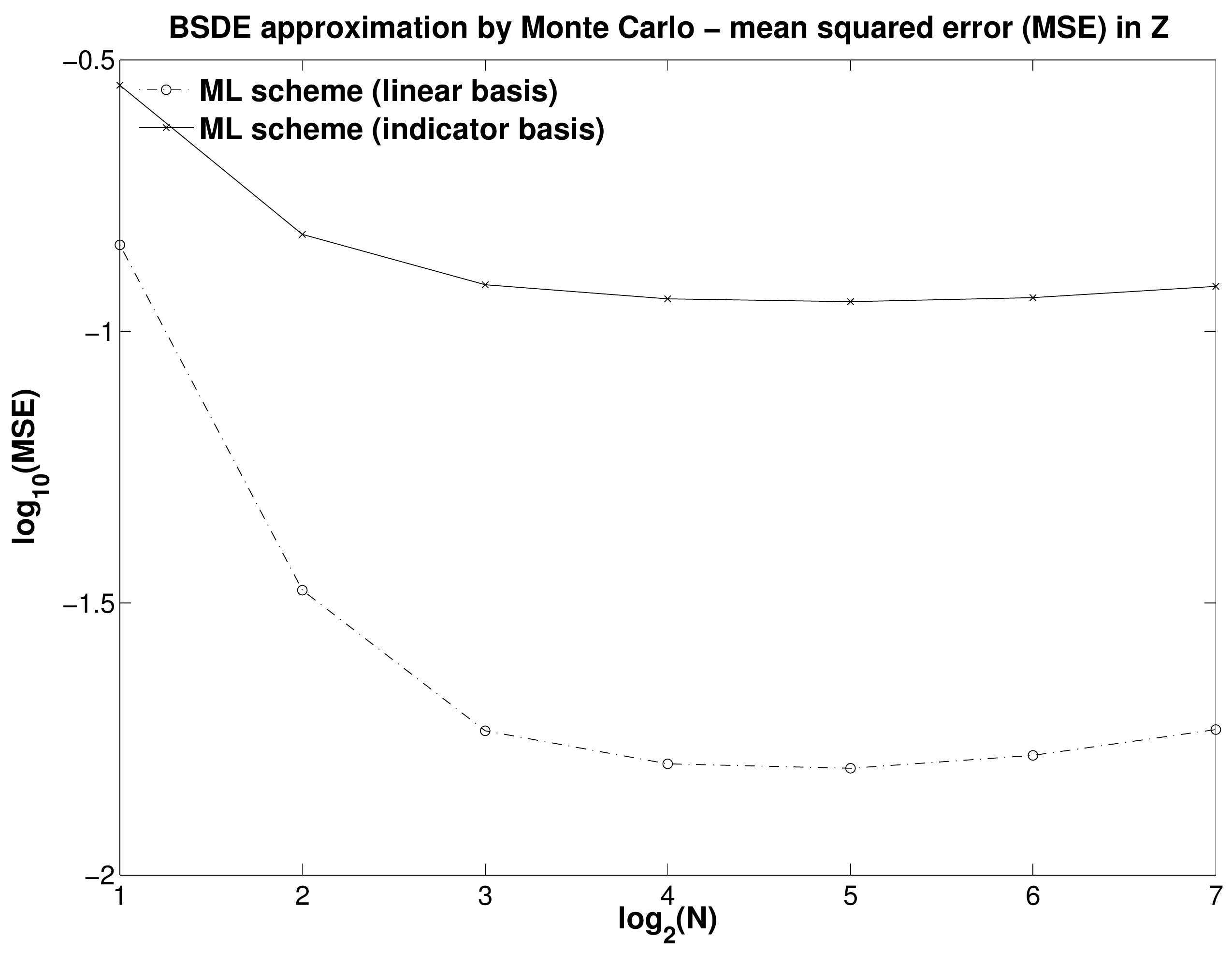}
\caption{MSE in $Z$ for multilevel in dimension 3}
\label{graph:multidim:a}
\end{figure}

\begin{figure}[h]
\centering
\includegraphics[width=11cm]{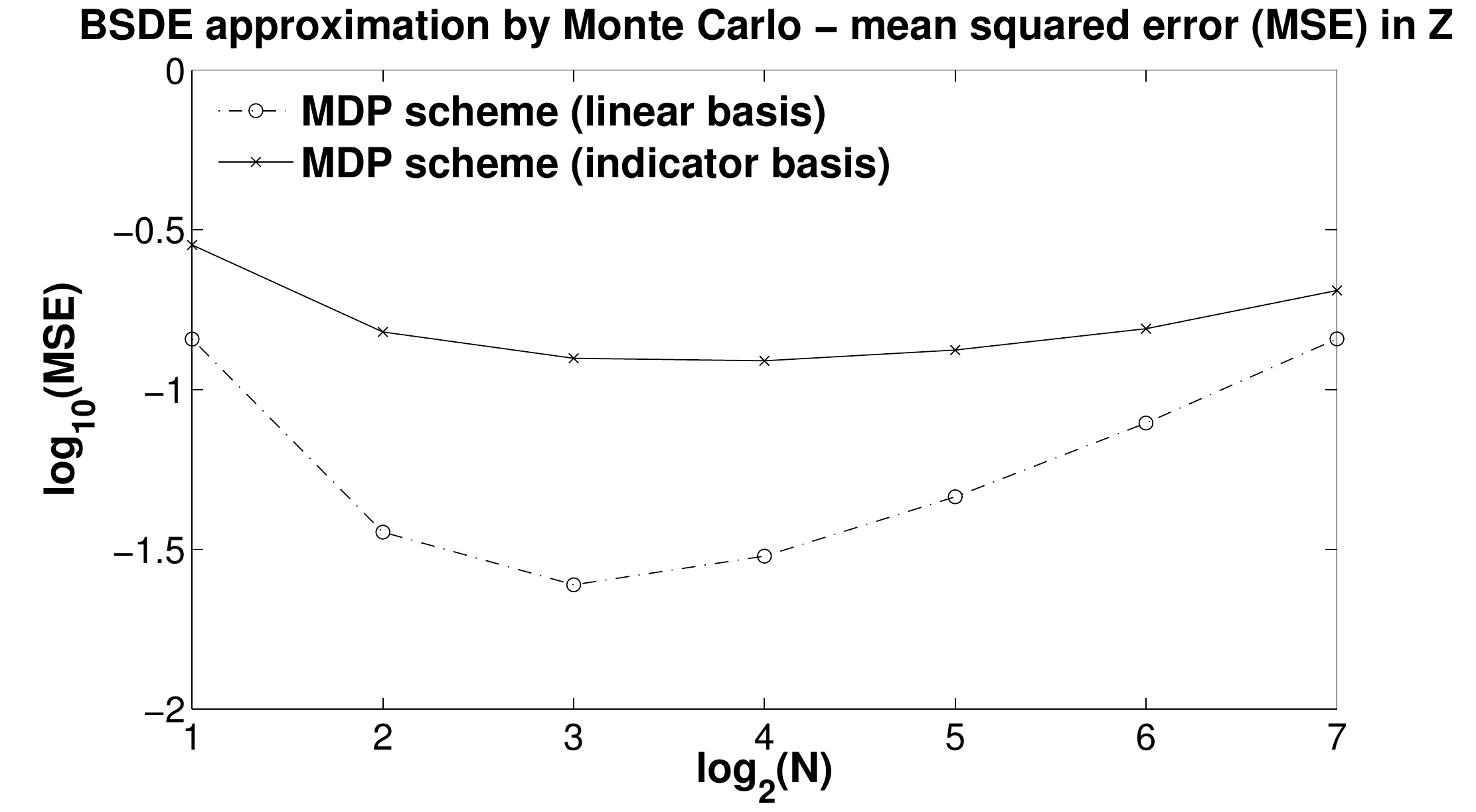}
\caption{MSE error in $Z$ for MDP in dimension 3}
\label{graph:multidim:b}
\end{figure}

\section{Completing the splitting algorithm}
\label{section:nonzero}
{
Fix $k >0$.}
In this section, approximate the second part of the split system \eqref{eq:nonl}, namely the functions $\bar y^{(k)}_i : \R^d \to \R$ and $\bar z^{(k)}_i : \R^d \to (\R^q)^\top$ in Lemma \ref{lem:markov:bsde}.
Omitting the superscript $(k)$ {to ease notation}
in what follows, we recall that these functions satisfy 
\begin{align*}
\bar y_i(X_i) & := \E_i\Big[   \sum_{j=i+1}^{2^k-1} f_j \big(X_j , \y k{j+1}{\X k{j+1}} + \bar y_{j+1}(\X k{j+1}) , \z k j{\X kj} + \bar z_j(\X kj)\big) \Delta_j \Big] , \\
 \Delta_i \times \bar z_i(X_i) & := \E_i\Big[\DW ki  \Big( \sum_{j=i+1}^{2^k-1} f_j \big(\X kj ,{ y_{j+1}}(\X k{j+1}) + \bar y_{j+1}(\X k{j+1}) ,{ z_j}(\X kj) + \bar z_j(\X kj)\big) \Delta_j\Big) \Big] .
\end{align*}

{Let} the functions $(y_j(\cdot),z_j(\cdot))_{0\le j \le N-1}$ {be} approximated using the multilevel algorithm, and 
{ $\pi$ denote} the time-grid $\pik k$ of the {highest} level of the multilevel algorithm.
{We use} 
 least-squares multistep dynamical programming (LSMDP) from\cite{gobe:turk:13a} for the discrete BSDE with zero terminal condition and random driver
\begin{equation}
\label{eq:proxy:driver}
\fM_j(y,z) :=  f_j \big(\X kj , \y {k,M}{j+1}{\X k{j+1}} + y, \z {k,M}{j+1}{\X k{j+1}} + z \big)
\end{equation}
to approximate $\bar y_i(\cdot)$ and $\bar z_i(\cdot)$.
{We maintain the superscript $k$ in notation $\X k{}$, $\DW k{}$, and $\tk k{}$  to remind that the time-grid in use is $\pik k$, although LSMDP does not make use of 
{earlier (coarser)}
time-grids.
The driver has two sources of randomness: the random functions $(\yn {k,M}{}, \zn {k,M}{})$, which depend on the samples used in the multilevel algorithm, {and}  the Markov chain $\X k{}$.
{ The notation}
\[
\fM_j(x_1,x_2,y,z) :=  f_j \big(x_1 , \y {k,M}{j+1}{x_2} + y, \z {k,M}{j}{x_1} + z \big),
\]
{will be helpful in the sequel. }
We briefly recall the LSMDP algorithm for the convenience of the reader.
 {Like the algorithm in Section \ref{section:MC}}, LSMDP is a least-squares Monte Carlo algorithm; 
the difference  in the choice of basis functions and the generation of simulations compared to the multilevel algorithm of Section \ref{section:MC} is threefold: firstly, since {there is no use of multiple levels}, only simulations of the Markov chain $\X k{}$ are generated; secondly, independent clouds of simulations are generated for every time-point, which means that the empirical measure for each time-point is independent of the empirical measure used at any other time-point; thirdly, the choice of basis functions is different to that used for the multilevel scheme.
We formalize  this in the following definitions.
\begin{definition}[Finite dimensional approximation spaces]
\label{def:fin-dim:space:2}
For 
$i\in \{0,\dots, 2^k-1\}$, we
 finite functional linear spaces 
{ $\LinSpace{Y,i}$ and $ \LinSpace{Z,i}$ }  of dimensions $\BasDim{Y,i}$ and $\BasDim{Z,i}$, given by
\begin{equation*}
\begin{cases}
  \LinSpace{Y,i} := {\rm span}\{\pl{Y,i}1,\ldots, \pl{Y,i}{\BasDim{Y,i} } \}, \text{ for } \pl{Y,i}l :\R^d \rightarrow \R \text{ s.t. }
\e[|\plx{Y,i}l{X_i}|^2] < +\infty,\\[2mm]
 \LinSpace{Z,i} := {\rm span}\{\pl{Z,i}1,\ldots, \pl{Z,i}{\BasDim{Z,i} } \}, \text{ for } \pl{Z,i}l :\R^d \rightarrow {(\R^q)^\top} \text{ s.t. }
\e[|\plx{Z,i}l{X_i}|^2] < +\infty.
\end{cases}
\end{equation*}
\end{definition}
We {suppress} the subscript $k$ in the notation of the basis functions and linear spaces to distinguish them from those in Definition \ref{def:fin dim:space}.
The functions $\bar y_i(\cdot)$ and  $\bar z_{i}(\cdot)$  will be approximated in the linear spaces $\LinSpace{Y,i}$ and  $\LinSpace{Z,i}$, respectively.
 will be approximated in $\LinSpace{Z,i}$.}
{We define } 
\begin{equation*}
{T^{Y}_{1,i}  := \inf_{ \phi \in \LinSpace{Y,i} }\e \Big[  |\phi(X_i) - \bar y_i(X_i)|^2 \Big]}
\quad \text{ and } \quad
{T^{Z}_{1,i}  := \inf_{ \phi \in \LinSpace{Z,i} }\e \Big[  |\phi(X_i) - \bar z_i(X_i)|^2 \Big];}
\end{equation*}
as in Definition \ref{def:fin dim:space}, these are the best approximation errors {possible with the} 
{choosen} basis functions.
\begin{definition}[Simulations and empirical measures]
\label{def:sims and empi:2}
For $i\in\{0,\ldots,2^k-1\}$, generate $M_i\geq 1$ independent copies  $\cC_{k,i} := \{( \DW{k,i,m}{i},\X{k,i,m}{} ) \; : \  m = 1,\dots,M_i \}$
 of $(\DW ki, \X k{})$: $\cC_{k,i}$ forms the {cloud of simulations} used for the regressions at time $i$.
We assume that the  clouds of simulations $(\cC_{k,i}:0\leq i < N)$ are independently generated, and are also independently generated from the clouds $\{ \cC_k \ : 0 \le k \le k\}$ of Definition \ref{def:sims and empi:2} used for the multilevel algorithm. 
{Let} $\emeasi k{i,M}$ {denote} the empirical measure of the $\cC_{k,i}$-simulations, 
i.e.\ 
\begin{equation*}
\emeasi k{i,M} = \frac 1{M_i} \sum_{m=1}^{M_i} \delta_{(\DW{k,i,m}i, \X {k,i,m}{i},\ldots,\X{k,i,m}{2^k} )}.
\end{equation*}
\end{definition}
We use the additional subscript $i$ in the notation for the clouds of simulations $\cC_{k,i}$ and the empirical measure $\emeasi k{i,M}$ for the LSMDP algorithm to distinguish them from those used for the multilevel algorithm and to specify the time-point.
As in Section \ref{section:MC}, 
{we enlage the probability space, while continuing to denote it for simplicity by}  $(\Omega, \cF, \P)$,
{to contain} also the simulations used for the LSMDP  \emph{and} the multilevel algorithms;
recall that the prototype processes {$W^{(k)}$ and $\X k{}$ } 
 are independent of all 
{simulation clouds}.

\begin{alg}
\label{def:lsmdp}
Recall the the linear spaces $\LinSpace{Y,i}$ and  $\LinSpace{Z,i}$
from 
Definition \ref{def:fin-dim:space:2}, 
 the empirical measures $\{\emeasi k{i,M} : i = 0,\dots,2^k-1\}$
from 
Definition \ref{def:sims and empi:2},
 the almost sure bounds from \Hf(iii),
 {the definition of the } truncation functions $\cT_L(\cdot)$ from Section \ref{section:notation},  { and OLS from Definition~\ref{def:ls}}.

Set $\byM{2^k}{\cdot} := 0$. 
For each $i=2^k-1,2^k-2,\dots,0$, set the random functions $\byM i\cdot$ and $\bzM i\cdot$ recursively as follows: 
 Define
$\byM i\cdot := \cT_{C_{y,i}}\big(\psim_{Y,i}(\cdot) \big)$
and $\bzM i\cdot = \cT_{{C_{z, i}}
}\big( \psim_{Z,i}(\cdot) \big)$, where
$C_{y,i} := C_X (T- \tk ki)^{(\thetaL+\theta)/2}$, $C_{z,i} := C_X (T- \tk ki)^{(\thetaL + \theta)/2} /\Del k_i$ and
\begin{equation}
\label{eq:PsiM}  
\left\{\begin{array}{l}
\displaystyle \psim_{Y,i}(\cdot)  \quad \text{solves} \quad
\OLS\left(\ObsM {Y,i}{\vecx{x}} \ ,  \ \LinSpace{Y,i} \ , \ \emeasi k{i,M}\right) 
 \\
\displaystyle \qquad \text{for} \quad \ObsM{Y,i}{\vecx{x}} :=  \sum_{j=i}^{2^k-1} \fM_k\big(x_j, x_{j+1}, \byM {j+1}{x_{j+1}}, \bzM{j}{x_j} \big ) \Del k_j, \quad\text{and}
\\
 \displaystyle \psim_{Z,i}(\cdot) \quad \text{solves} \quad
\OLS\left( \ \ObsM{Z,i}{w,\vecx{x}} \ , \ \LinSpace{Z,i} \ , \ \emeasi k{i,M}\right)\\  \\
\displaystyle \quad \text{for} \quad \ObsM{Z,i}{w,\vecx{x}} := 
\frac 1{\Del k_i}\ObsM{Y,i+1}{\vecx{x}} \ w^\top,\qquad\text{for $w \in \R^q$, $\vecx{x} = (x_0,\dots,x_{2^k}) \in (\R^d)^{2^k+1}$. }
\end{array}\right. 
\end{equation}
\end{alg}

We now come to the main result of this section, which is the error analysis of the LSMDP algorithm.
\begin{theorem}[{Error for the LSMDP scheme}]
\label{thm:mc:err decomp}
Recall the constants $C_{y,i}$ and $C_{z,i}$ from Algorithm \ref{def:lsmdp}.
For each $j \in \{0,\dots,2^k-1\}$, define
\[
  \cE(j) := 
 T^{Y}_{1,j}
 +  {T^{Z}_{1,j}} + C_S^2\Big(\frac{ 3\BasDim{Y,j} }{M_j } + 2q\frac{\BasDim{Z,j} }{\Del k_j M_j} \Big)
 + 800 \Big(C_{y,j}^2(\BasDim{Y,j}+1)+\Cz j^2(\BasDim{Z,j}+1)q \Big) \frac{\log(3M_j)}{M_j}.
\]
where \[{C_S := \sum_{i=0}^{2^k-1} \Big\{C_f + {L_f (C_{y,i} + C_y + C_{z,i} + \Cz {k,i}) \over (T-\tk ki)^{(1-\thetaL)/2} } \Big\} \Del k_i.} \]
Recall $C_{\pik k} $ from \Htp \ and assume that  $k$ is sufficiently large so that $C_{\pik k} L^2_f (\cpib \vee 1) \leq ({384}  (2q + (1+T)e^{T/2} )(1+T))^{-1}$, 
and that the parameters of the multilevel algorithm are such that the global error is estimated by $\bar \cE(k) \le \varepsilon$ for some $\varepsilon >0$ (for definition of $\bar \cE(k)$, see equation \eqref{eq:glob:er} and subsequent remarks).
Then, for all 
$0\leq i \leq 2^k-1$, 
\begin{align}
& \e\Big[ \frac 1{M_i}\sum_{m=1}^{M_i} |   \byM i{\X {k,i,m}i} - \bar y_i(\X {k,i,m}i) |^2 \Big] 
\le 
{T^{Y}_{1,i}}  +  \frac{ 3 C_S^2 \BasDim{Y,i} }{M_i} 
+ C_\Gamma (1+T) \varepsilon
+ C_\Gamma \sum_{j=i}^{2^k-1} \cE(j){\Del k_j},
\label{eq:empi:y} \\
& \sum_{j=i}^{2^k-1}  \e \Big[ \frac 1{M_j}\sum_{m=1}^{M_j}  | \bzM j{\X {k,j,m}j}-\bar z_j(\X {k,j,m}j) |^2 \Big] \Del k_j  \le
C_\Gamma (1+T) \varepsilon +  C_\Gamma \sum_{j=i}^{2^k-1} \cE(j) \Del k_j,
\label{eq:empi:z}
\end{align}
where $C_\Gamma = {8} \exp \big( {384} (\cpib \vee 1) (2q + (1+T)e^{T/2} ) (1+T)T^{\thetaL} L_f^2 / \thetaL \big)$.
\end{theorem}
The proof of the above theorem is analogous to the proof of \cite[Theorem 4.11]{gobe:turk:13a}. 
Indeed, the proof of \cite[Theorem 4.11]{gobe:turk:13a} relies only on conditioning arguments, a-priori estimates, concentration of measure inequalities, and elementary properties of ordinary least-squares regression (Proposition \ref{prop:ls:reg:properties});
the a-priori estimates \cite[Proposition 3.2]{gobe:turk:13a} admit randomness in the driver, and the properties of ordinary least-squares regression are universal, therefore these arguments  
{require no alterations for} our setting. The concentration of measure result is provided in Proposition \ref{prop:conc:meas:nonl} below.
There are three minor adaptations to the analysis, which we now detail for the convenience of the reader.
Firstly, one must augment the $\sigma$-algebras used in the conditioning arguments   by adding $\sigma(\cC_0, \ldots , \cC_k)$, the $\sigma$-algebra of the simulations used in the multilevel algorithm {for the first part of the split system}.
Secondly, after the application of the a-priori estimates,  one must estimate
\begin{equation}
\label{eq:diff:drivers}
\E[| f_j(\y k{j+1}{\X k{j+1}} + \bar y_{j+1}(\X k{j+1}) , \z kj{\X kj} + \bar z_j(\X kj)  ) - \fM_j ( \byM {j+1}{\X k{j+1}}, \bzM j{\X kj})|^2]
\end{equation}
whereas, in  the {respective} computation in  \cite[equation (33)]{gobe:turk:13a}, one only needed to estimate
\(
\E[| f_j(y_{j+1}(\X k{j+1}) , z_j(\X kj) ) - f_j(\yM {j+1}{\X k{j+1}} , \zM j{\X kj} )|^2].
\)
Recalling that 
\[
\fM_j(y,z) :=  f_j (\X kj , \y {k,M}{j+1}{\X k{j+1}} + y, \z {k,M}{j+1}{\X k{j+1}} + z ),
\]
one uses the Lipschitz continuity of $f_j(y,z)$ and the hypothesis that  the approximation of $(\yn k{},\zn k{})$ by $(\yn {k,M}{} , \zn {k,M}{})$ produces a global error less {than} or equal to $\varepsilon$ to estimate the error due to multilevel.
Finally, one replaces the constant $C_{(4.7)}$ in \cite[Theorem 4.11]{gobe:turk:13a} by $C_S$; the explicit value of $C_S$ is obtain exactly as the explicit value of $C_{(4.7)}$ in \cite[Lemma 4.7]{gobe:turk:13a}, only using the almost absolute bounds of $\byM i\cdot$ and $\bzM i\cdot$.

\section{Conclusion: Comparison of the  schemes with and without splitting and multilevel}
\label{section:conclusion}
Using the results of Sections \ref{section:MC} and \ref{section:nonzero}, we are now in a position to compare our algorithm (splitting combined with multilevel)
to the  least-squares multistep dynamical programming (LSMDP) scheme  with neither.
For simplicity, we will assume that $\theta = \thetaL = 1$, meaning that the terminal condition $\Phi(\cdot)$ is Lipschitz continuous (but not necessarily differentiable) and that the driver is uniformly Lipschitz continuous in $(y,z)$.
For the remainder of this section, we write $g(y) = O(y)$ if there exists a constant $C$ independent of $k$ and $y$ such that $g(y)/y \to C$ as $y \to 0$.
For given precision level $\varepsilon>0$, it is our goal to set the basis functions and the number of simulations for each time-point of the approximation of $(\bar y_i,\bar z_i)$ so that the global error {satisfies}
\begin{equation}
\bar\cE(M) := \max_{0\le i \le 2^k-1} \e[ |   \byM i{\X ki} - \bar y_i(\X ki) |^2 ] + \sum_{i=0}^{2^k} \e[ |   \bzM i{\X ki} - \bar z_i(\X ki) |^2 ] \Del k_i \le O(\varepsilon).
\label{er:nonl}
\end{equation}
{To apply}   Theorem \ref{thm:mc:err decomp}, 
{we first} provide a concentration of measure result.
\begin{proposition}\label{prop:conc:meas:nonl} 
For each $k\in\{0,\ldots,\kappa\}$ and $i \in \{0,\dots , 2^k-1\}$, we have
\begin{align}
  &\e[ |   \byM i{\X ki} - \bar y_i(\X ki) |^2 ]
 \le 2  \e \Big[ \frac 1{M_i}\sum_{m=1}^{M_i}  | \byM i{\X {k,i,m}i}-\bar y_i(\X {k,i,m}i) |^2 \Big] 
+  {2028 (\BasDim {Y,i}+1) C_{y,i}^2 \log(3M_i) \over M_i},
\nonumber
\\
& \e[ |   \bzM i{\X ki} - \bar z_i(\X ki) |^2 ] \le 2 \e \Big[ \frac 1{M_i}\sum_{m=1}^{M_i}  | \bzM i {\X {k,i,m}i}-\bar z_i (\X {k,i,m}i) |^2 \Big]  
+  {{2028 (\BasDim {Z,i}+1) q \Cz{i}^2 \log(3M_i) } \over  M_i};
\nonumber
\end{align}
we recall that $C_{y,i} = C_X^2(T-\tk ki)^{(1+\theta)/2}$ and $\Cz{i} = C_X^2(T-\tk ki)^{\theta/2}$.
\end{proposition}
Just as Proposition \ref{prop:eq:y:z:err:deco:M}, Proposition \ref{prop:conc:meas:nonl} is analogous to \cite[Proposition 4.10]{gobe:turk:13a}.
The second 
{terms} on the right hand side of the inequalities  are  correction terms which can be interpreted as { interdependency } errors due to the change of the inner measure.
We see the interdependency errors have  the same dependence on $\BasDim {\cdot,i}$ and $M_i$ as the last term in $\cE(i)$ in Theorem \ref{thm:mc:err decomp}.
{Hence,} to ensure \eqref{er:nonl},
 it is  sufficient to set the numerical parameters so that the local error terms {satisfy} $\cE(i) \le O(\varepsilon)$ for every $i \in \{ 0,\ldots, 2^k -1\}$.
Using \HXpp(ii), we  can {replace} $T^Y_{1,i}$ and $T^Z_{1,i}$ {by} 
\[T^Y_{2,i} := \inf_{\phi\in \LinSpace{Y,i} }\e[ |  \phi(X_{\tk ki}) - U(\tk ki , X_{\tk ki}) |^2 ] \ \text{ and } \ T^Z_{2,i}:= \inf_{\phi \in\LinSpace{Z,i}}\e[ |   \phi(X_{\tk ki}) - V(\tk ki, X_{\tk ki}) |^2 ],\] 
respectively, in the local error term $\cE(i)$, 
and choose basis functions
{such} that $T^Y_{2,i}$ and $T^Z_{2,i}$ are dominated by $O(\varepsilon)$ for every $i$.
Thanks to the 
Lipschitz continuity in \Hf(iii),
 it is sufficient to use (for every time-point and both for $Y$ and $Z$) a basis of 
 functions 
{on a partition of disjoint hypercubes  with} diameter $O(\sqrt{\varepsilon})$.
This basis is infinite dimensional, but one can make a simple truncation to get around this problem. 
We assume additionally (as in \cite[Section 4.4]{gobe:turk:13a}) that, for each $i \in \{ 0,\dots, 2^k-1\}$, $\X ki$ has exponential moments, so that we may set the basis in the region outside $[-R,R]^d$ to zero for $R = \ln (\varepsilon^{-1} + 1)$; 
this truncation induces an error $O(\varepsilon)$, which is admissible.
The dimension of the hypercube basis $\BasDim {l,i}$ is therefore, uniformly in  $l \in \{Y,Z\}$ and $i \in \{ 0, \ldots ,2^k-1\}$, equal to $O(\varepsilon^{-d/2} \ln(\varepsilon^{-1}+1)^d)$.
It follows that we must choose the number of simulations $M_i$ to be equal, uniformly in $i$, to $O(2^{k}\varepsilon^{-1-d/2} \ln(\varepsilon^{-1} +1)^d)$.

It remains only to compute the complexity of the scheme.
There are two contributions to the computational cost: the cost of simulation of the Markov chain $\X k{}$ and Brownian incremends $\DW k{}$, and the cost of the regressions.
The cost of simulation is equal to $O(2^{3k}\varepsilon^{-1-d/2} \ln(\varepsilon^{-1}+1)^d)$;
the additional factor $2^k$ comes from re-simulation the paths of the Markov chain $\X k{}$ at every time-step.
Since we are using the partitioning estimate to compute the regression coefficient, the regression cost is equal to $\sum_{j=0}^{2^k -1} O( M_i) = O(2^{2k}\varepsilon^{-1-d/2} \ln(\varepsilon^{-1}+1)^d)$;
see Section \ref{section:ML:err an}  for details on the partitioning estimate.
Therefore, recalling that $2^k = \varepsilon^{-1}$, the overall complexity is equal to equal to $O(\varepsilon^{-4-d/2} \ln(\varepsilon^{-1}+1)^d)$. 
Therefore, using {the assumptions and} computations  of Section \ref{section:ML:err an}, it follows that the overall complexity of the splitting scheme with multilevel, i.e.\ Algorithms \ref{def:lsmdp} and \ref{alg:mlmc} together, is
\begin{equation}
\label{eq:complexity:full}
O(\varepsilon^{-2 - d} \ln(\varepsilon^{-1}+1)) + O(\varepsilon^{-4 - d/2}\ln(\varepsilon^{-1}+1)^d).
\end{equation}

{
We now calibrate the LSMDP algorithm with no splitting and no multilevel, which we recall  below for completeness in Algorithm \ref{def:lsmdp:full}. 
We then compute the complexity of this algorithm in order to provide a suitable comparison to an established algorithm and determine the possible gains of the splitting algorithm with multilevel.
\begin{alg}
\label{def:lsmdp:full}
Recall the the linear spaces $\LinSpace{Y,i}$ and  $\LinSpace{Z,i}$
from 
Definition \ref{def:fin-dim:space:2}, 
 the empirical measures $\{\emeasi k{i,M} : i = 0,\dots,2^k-1\}$
from 
Definition \ref{def:sims and empi:2},
 the  bounds from \Hf(iii),
  and the truncation function $\cT_L(\cdot)$ from Section \ref{section:notation}.
Set $\yM{2^k}{\cdot} := \Phi(\cdot)$. 
For each $i=2^k-1,2^k-2,\dots,0$, set the random functions $\yM i\cdot$ and $\zM i\cdot$ recursively as follows: 
 Define
$\yM i\cdot := \cT_{C_{y}}\big(\psim_{Y,i}(\cdot) \big)$
and $\zM i\cdot = \cT_{{C_{z, i}}
}\big( \psim_{Z,i}(\cdot) \big)$, where
$C_{y} := C_X $, $C_{z,i} := \Cz {k,i}$ and
\begin{equation*}
\left\{\begin{array}{l}
\displaystyle \psim_{Y,i}(\cdot)  \quad \text{solves} \quad
\OLS(\ObsM {Y,i}{\vecx{x}} \ ,  \ \LinSpace{Y,i} \ , \ \emeasi k{i,M}) 
 \\
\displaystyle \qquad \text{for} \quad \ObsM{Y,i}{\vecx{x}} :=  \Phi(x_N) + \sum_{j=i}^{2^k-1} f_k\big(x_j, \yM {j+1}{x_{j+1}}, \zM{j}{x_j} \big ) \Del k_j, \quad \text{and}
\\
 \displaystyle \psim_{Z,i}(\cdot) \quad \text{solves} \quad
\OLS( \ \ObsM{Z,i}{w,\vecx{x}} \ , \ \LinSpace{Z,i} \ , \ \emeasi k{i,M})\\ 
\displaystyle \qquad \text{for} \quad \ObsM{Z,i}{w,\vecx{x}} := 
\frac 1{\Del k_i}\ObsM{Y,i+1}{\vecx{x}} \ w^\top,\qquad {\text{for $w \in \R^q$, $\vecx{x} = (x_0,\dots,x_{2^k}) \in (\R^d)^{2^k+1}$.}}
\end{array}\right.
\end{equation*}
\end{alg}
}
{
The error  of this algorithm is studied in \cite[Theorem 4.11]{gobe:turk:13a}.
With this algorithm, we are directly approximating the continuous time function $v(t,\cdot) + V(t,\cdot)$.
The complexity analysis for Algorithm \ref{def:lsmdp:full} is the same as that for Algorithm \ref{def:lsmdp} above, however we must take into account the additional weight due to the time-dependency of the Lipschitz coefficient:
the Lipschitz constant for $v(t,\cdot) + V(t,\cdot)$ is equal to  $O((T-t)^{-1/2})$ for all $t \in [0,T)$ - see assumptions \HX(iv) and \Hf(iii).
Therefore, we choose a hypercube basis for each time-point $i\in\{0,\ldots,2^k-1\}$ whose cubes have diameter $\sqrt{T-\tk ki} O(\sqrt{\varepsilon})$.
Therefore, the overall complexity of Algorithm \ref{def:lsmdp:full} is
\[O(\varepsilon^{-3 - d/2} \ln(\varepsilon^{-1}+1)^d)\sum_{i=0}^{2^k-1}(T-\tk ki)^{d/2} \le O(\varepsilon^{-4 - d} \ln(\varepsilon^{-1}+1)^d) .\]
Compared with the two terms in \eqref{eq:complexity:full},  the complexity of Algorithm \ref{def:lsmdp:full} dominates: 
if $d <4$, \eqref{eq:complexity:full} is dominated by $O(\varepsilon^{-4 - d/2}\ln(\varepsilon^{-1}+1)^d))$, whereas for $d\ge 4$, \eqref{eq:complexity:full} is dominated by $O(\varepsilon^{-2 - d} \ln(\varepsilon^{-1}+1))$. 
Therefore, 
one gaines two orders in $\varepsilon$ in high dimension $d \ge 4$ thanks to the use of splitting algorithm with multilevel.
If we were to use a splitting method but no multilevel (i.e., LSMDP as in Remark \ref{rem:best/worst:gen}) in the {zero driver} part, 
{the} gain compared to pure LSMDP 
 would {still be substantial, being of} order one in $\varepsilon$,
however an additional order is to be gained by multilevel (for $d \ge 4$), cf. Section~\ref{section:ML:err an},

}

Finally, 
{we support the theory for the non-linear generator 
by computational results for an} example from finance.
To this end, consider a $d=2$-dimensional  forward process  $X=(S,H)$ for correlated geometric Brownian motions 
\(dS = S \sigma^S dW^1\) with \(S_0=1,\) and 
\begin{align}
dH&= H \left( \gamma dt + \sigma^H(\rho dW^1 + \sqrt{1-\rho^2} dW^2)\right), \quad H_0=1,
\end{align} 
with parameters $ \sigma^S=\sigma^H=0.5, \rho=0.6$ and $\gamma=0.1$. 
 Considering $S$ and $H$ 
 as the (discounted) price processes
of a liquidly tradable risky asset and  of a non-tradable asset,  the so called no-good-deal
{valuation bound} $Y$  for an option $\Phi(X_T)=(H_T-S_T)^+$ to  exchange at maturity $T$ one traded asset $S_T$ into one non-traded asset $H_T$
is described by the non-linear BSDE
\begin{align}
 dY_t &= -h|Z^{(2)}_t|dt + Z_t dW= -|Z^{(2)}_t| dt + Z^{(1)}_t dW^{(1)}_t + Z^{(2)}_t dW^{(2)}_t \,,\quad
Y_T = \Phi(X_T), 
\end{align}
where $Z = (Z^{(1)},Z^{(2)})$;
{for a good-deal constraint that we take as $ h=0.2$, see \cite{Becherer09,kentiab14}}.
The BSDE  has an explicit solution in terms of a Margrabe-type formula, see \cite{kentiab14}, and 
a  corresponding good-deal hedging strategy can be obtained from $Z$.

The regression basis for ML and MDP is given by indicator functions on the hypercubes of a
partition of $\R^2$ into $K=50^2$ sets. Number of simulations is $M=2*10^6$ for both schemes;
{for multilevel (ML), the same number of simulations is used at every level.}
Results on mean squared errors are reported in Table~\ref{tab:bsp4}.
They show substantial error reduction by multilevel (ML) in combination with the splitting scheme, 
in particular for the MSE in $Z$ for finer time grids (larger $k=\log_2 N$),
confirming insights as 
{before} also for the present example with non-zero generator. 
  
\begin{table}[h]
\centering
\caption{\textbf{MSE in $Y$ and $Z$} with non-zero generator}
\label{tab:bsp4}
\small
\begin{tabular}{lrrrrrrrrrr}
$\text{log}_{2}(N)$
&\textbf{1}&\textbf{2}&\textbf{3}&\textbf{4}&\textbf{5}\\
\hline
MDP Y&0.1372&0.0795&0.0515&0.0379&0.0322\\
MDP Z&0.0161&0.0089&0.0092&0.0143&0.0253\\
\hline
ML Y&0.1371&0.0791&0.0510&0.0373&0.0314\\
ML Z&0.0156&0.0068&0.0039&0.0032&0.0031\\
\end{tabular}
\end{table}

\noindent
{\bf \large Acknowledgements:}
{We like to thank Emmanuel Gobet for advice on  this paper and 
the thesis \cite{turk:13}
where {a first multilevel scheme has been} introduced, 
and  for pointing out
the special basis used in
 Theorem \ref{thm:spec:basis} in particular.
We 
thank Axel Mosch 
and Klebert Kenita 
 {for help with the examples}.
}

\end{document}